\documentclass[draft,preprint,12pt,numbers,sort&compress]{elsarticle}
\usepackage{}
\usepackage{mathrsfs}
\usepackage{amssymb}
\usepackage{amsthm}
\usepackage{mathrsfs}
\usepackage[centertags]{amsmath}
\usepackage{amsfonts}

\usepackage{color}

\input amssym.def
\input amssym.tex

\date{}

\newtheorem{Theorem}{Theorem}[section]

\newtheorem{Lemma}{Lemma}[section]

\newcommand\R{\mbox{\bf R}}

\newcommand\Z{\mbox{\bf Z}}
\newcommand\z{\mbox{\bf z}}
\newcommand\T{\mbox{\bf T}}
\newcommand\SR{\mbox{\scriptsize\bf R}}

\newcommand{\definition}{{\lower .5ex
  \hbox{$\>\>\stackrel{\triangle}{=}\>\>$} }}
\newcommand\supp{\mathop{\rm supp}}

\begin{document}

\baselineskip=22pt
\thispagestyle{empty}

\mbox{}
\bigskip

\begin{center}
{\Large \bf The  Cauchy problem for the  generalized   KdV equation with rough data and random data}\\[1ex]

{Wei Yan\footnote{ Email:011133@htu.edu.cn}$^a$,
\quad Xiangqian Yan\footnote{Email:yanxiangqian213@126.com}$^a$,
\quad  Jinqiao  Duan
\footnote{ Email: duan@iit.edu}$^{b}$,     Jianhua Huang\footnote{ Email:jhhuang32@nudt.edu.cn}$^c$}\\[1ex]

{$^a$College of Mathematics and Information Science, Henan Normal University,}\\
{Xinxiang, Henan 453007,   China}\\[1ex]

{$^b$ Department of Applied Mathematics, Illinois Institute of Technology,}\\[1ex]
{Chicago, IL 60616, USA  }\\[1ex]
{$^c$College of Science, National University of Defense  Technology,}\\
{ Changsha, Hunan 410073,  China}\\[1ex]

\end{center}

\bigskip
\bigskip

\noindent{\bf Abstract.}    In this paper,
 we consider the Cauchy problem for the generalized
KdV equation
\begin{eqnarray*}
      u_{t}+\partial_{x}^{3}u+\frac{1}{k+1}(u^{k+1})_{x}=0,k\geq5
    \end{eqnarray*}
    with rough data and random data.
     Firstly, by using the Fourier restriction norm method and the
     frequency truncated technique,
  we prove that $u(x,t)\longrightarrow u(x,0)$ as $t\longrightarrow0$
  for a.e. $x\in \R$ with $u(x,0)\in H^{s}(\R)(s>\frac{1}{2}-\frac{2}{k},k\geq8).$
   Secondly, we  prove that
   $u(x,t)\longrightarrow  e^{-t\partial_{x}^{3}}u(x,0)$ as $t\longrightarrow0$
  for a.e. $x\in \R$ with $u(x,0)\in H^{s}(\R)(s>\frac{1}{2}-\frac{2}{k},k\geq8).$
   Thirdly, by using the Fourier restriction norm method,
  we prove that  $
\lim\limits_{t\longrightarrow 0}\left\|u(x,t)-e^{-t\partial_{x}^{3}}u(x,0)\right\|_{L_{x}^{\infty}}=0$
 with   $u(x,0)\in H^{s}(\R)(s>\frac{1}{2}-\frac{2}{k+1},k\geq5)$.
 Fourthly,
by using the Fourier restriction norm method, Strichartz estimates,
Wiener randomization of the initial data,
   probabilistic  Strichartz estimates,    we  establish
 the probabilistic
well-posedness in $H^{s}(\R)\left(s>{\rm max}
\left\{\frac{1}{k+1}\left(\frac{1}{2}-\frac{2}{k}\right),
\frac{1}{6}-\frac{2}{k}\right\}\right)$ with random data.
   Our result  improves the result of
   Hwang,  Kwak (Proc. Amer. Math. Soc.  146(2018), 267-280.).
 Fifthly,
  we prove that $\forall \epsilon>0,$ $\forall \omega \in \Omega_{T},$
  $\lim\limits_{t\longrightarrow0}\left\|u(x,t)-e^{-t\partial_{x}^{3}}u^{\omega}(x,0)\right\|_{L_{x}^{\infty}}=0$
   with $u(x,0)\in H^{s}(\R)(s>\frac{1}{6},k\geq6),$ where ${\rm P}(\Omega_{T})\geq 1- C_{1}{\rm exp}
\left(-\frac{C}{T^{\frac{\epsilon}{48k}}\|u(x,0)\|_{H^{s}}^{2}}\right)$
    and $u^{\omega}(x,0)$  is the randomization of   $u(x,0)$.
   Finally, we prove that  $\forall \epsilon>0,$ $\forall \omega \in \Omega_{T},$
  $\lim\limits_{t\longrightarrow0}\left\|u(x,t)-u^{\omega}(x,0)\right\|_{L_{x}^{\infty}}=0$
    with  ${\rm P}(\Omega_{T})\geq 1- C_{1}{\rm exp}
\left(-\frac{C}{T^{\frac{\epsilon}{48k}}\|u(x,0)\|_{H^{s}}^{2}}\right)$
 and $u(x,0)\in H^{s}(\R)(s>\frac{1}{6},k\geq6).$


\bigskip

\leftskip 0 true cm \rightskip 0 true cm

\newpage{}

\begin{center}

{\Large \bf The  Cauchy problem for the  generalized   KdV equation
 with rough data and random data}\\[1ex]
{Wei Yan\footnote{ Email:011133@htu.edu.cn}$^a$,
\quad Xiangqian Yan\footnote{Email:yanxiangqian213@126.com}$^a$,
\quad  Jinqiao  Duan
\footnote{ Email: duan@iit.edu}$^{b}$,
   Jianhua Huang\footnote{ Email:jhhuang32@nudt.edu.cn}$^c$}\\[1ex]

\end{center}
  In this paper,
 we consider the Cauchy problem for the generalized
KdV equation
\begin{eqnarray*}
      u_{t}+\partial_{x}^{3}u+\frac{1}{k+1}(u^{k+1})_{x}=0,k\geq5
    \end{eqnarray*}
    with rough data and random data.
     Firstly, by using the Fourier restriction norm method and the
     frequency truncated technique,
  we prove that $u(x,t)\longrightarrow u(x,0)$ as $t\longrightarrow0$
  for a.e. $x\in \R$ with $u(x,0)\in H^{s}(\R)(s>\frac{1}{2}-\frac{2}{k},k\geq8).$
   Secondly, we  prove that
   $u(x,t)\longrightarrow  e^{-t\partial_{x}^{3}}u(x,0)$ as $t\longrightarrow0$
  for a.e. $x\in \R$ with $u(x,0)\in H^{s}(\R)(s>\frac{1}{2}-\frac{2}{k},k\geq8).$
   Thirdly, by using the Fourier restriction norm method,
  we prove that  $
\lim\limits_{t\longrightarrow 0}\left\|u(x,t)-e^{-t\partial_{x}^{3}}u(x,0)\right\|_{L_{x}^{\infty}}=0$
 with   $u(x,0)\in H^{s}(\R)(s>\frac{1}{2}-\frac{2}{k+1},k\geq5)$.
 Fourthly,
by using the Fourier restriction norm method, Strichartz estimates,
Wiener randomization of the initial data,
   probabilistic  Strichartz estimates,    we  establish
 the probabilistic
well-posedness in $H^{s}(\R)\left(s>{\rm max}
\left\{\frac{1}{k+1}\left(\frac{1}{2}-\frac{2}{k}\right),
\frac{1}{6}-\frac{2}{k}\right\}\right)$ with random data.
   Our result  improves the result of
   Hwang,  Kwak (Proc. Amer. Math. Soc.  146(2018), 267-280.).
 Fifthly,
  we prove that $\forall \epsilon>0,$ $\forall \omega \in \Omega_{T},$
  $\lim\limits_{t\longrightarrow0}\left\|u(x,t)-e^{-t\partial_{x}^{3}}u^{\omega}(x,0)\right\|_{L_{x}^{\infty}}=0$
   with $u(x,0)\in H^{s}(\R)(s>\frac{1}{6},k\geq6),$ where ${\rm P}(\Omega_{T})\geq 1- C_{1}{\rm exp}
\left(-\frac{C}{T^{\frac{\epsilon}{48k}}\|u(x,0)\|_{H^{s}}^{2}}\right)$
    and $u^{\omega}(x,0)$  is the randomization of   $u(x,0)$.
   Finally, we prove that  $\forall \epsilon>0,$ $\forall \omega \in \Omega_{T},$
  $\lim\limits_{t\longrightarrow0}\left\|u(x,t)-u^{\omega}(x,0)\right\|_{L_{x}^{\infty}}=0$
    with  ${\rm P}(\Omega_{T})\geq 1- C_{1}{\rm exp}
\left(-\frac{C}{T^{\frac{\epsilon}{48k}}\|u(x,0)\|_{H^{s}}^{2}}\right)$
 and $u(x,0)\in H^{s}(\R)(s>\frac{1}{6},k\geq6).$

\bigskip

{\large\bf 1. Introduction}
\bigskip

\setcounter{Theorem}{0} \setcounter{Lemma}{0}

\setcounter{section}{1}

In this paper, we consider the Cauchy problem for the generalized KdV equation
\begin{eqnarray}
     && u_{t}+\partial_{x}^{3}u+\frac{1}{k+1}(u^{k+1})_{x}=0,k\geq5\label{1.01},\\
     &&u(x,0)=f(x)\label{1.02}
 \end{eqnarray}
 with rough data and random data.

Carleson \cite{Carleson} initiated the study on   the pointwise convergence  problem of  one dimensional Schr\"odinger equation in $H^{s}(\R)$, $s\geq 1/4$.
   Dahlberg and Kenig \cite{DK} showed that $s\geq \frac{1}{4}$ is
the necessary condition for the pointwise
convergence problem of the Schr\"odinger equation in any dimension. Dahlberg,  Kenig \cite{DK} and Kenig  et al. \cite{KPV1991,KPVCPAM}
  have showed that $s\geq \frac{1}{4}$
  is the necessary and sufficient condition for
  the pointwise convergence problem
of  KdV equation in $H^{s}(\R)$.    Bourgain \cite{Bourgain2016} presented counterexamples about the pointwise convergence problem of Schr\"odinger
 equation showing
 that the pointwise convergence problem requires  $s\geq\frac{n}{2(n+1)}$. Luc$\grave{a}$ and Rogers \cite{LR2019} presented a new proof of Proposition 1 of \cite{Bourgain2016}.
  Du et al. \cite{DGL} proved
 that when   $s>\frac{1}{3},$
 the pointwise convergence problem
  of two dimensional
  Schr\"odinger
 equation in $H^{s}(\R)$ can hold.
 Du and Zhang \cite{DZ}
  proved that when  $s>\frac{n}{2(n+1)},n\geq3,$   the
  pointwise convergence problem of $n$ dimensional Schr\"odinger
 equation in $H^{s}(\R)$ can hold.
Rogers and Villarroya \cite{RV} have proved that when $s>{\rm max}\left\{n(\frac{1}{2}-\frac{1}{q}),
\frac{n+1}{4}-\frac{n-1}{2q},\frac{1}{2}\right\}(q\geq1)$, the pointwise convergence problem of   the wave equation
$\frac{1}{2}\left[e^{it\sqrt{-\Delta}}+e^{-it\sqrt{-\Delta}}\right]
f\longrightarrow f$ in  $H^{s}(\R^{n})$ can hold.
Some  authors have investigated  the pointwise
 convergence problem of the  Schr\"odinger equation  in higher  dimension and
 other dispersive equations
\cite{Bourgain1992,Bourgain1995,Bourgain2013,CLV,Cowling,D,DG,DGZ,DK,GS,MVV,Lee,
LR2015,LR2017,MYZ2015,MZZ2015,KPV1991,KPVCPAM,RVV,Shao,S,Tao,TV,Vega}.
Some authors have investigated
 the pointwise
 convergence problem of Schr\"odinger equation on the torus $\mathrm{\T}^{n}$ \cite{MV,WZ,CLS,EL,WYY}.

  Lebowitz et al.\cite{LRS}  studied the
 statistical mechanics of the
nonlinear Schr\"odinger equation with random data. Bourgain \cite{B1994CMP, B1994Duke, B1996CMP} initiated the study
on  the
problem of constructing the invariant measures for the Sch\"ordinger equation and  Zakharov system with random data.
By establishing the Khintchine-type  inequality,   Burq and Tzvetkov \cite{BT2008-L,BT2008-G} established
the probabilistic well-posedness of the wave equation on compact manifolds in supercritical regime with random data.
The idea  of  randomized initial data are effective in constructing the invariant measures and
    establishing
 local and global  strong solution for a large set of initial
  data
 in supercritical regime   with random  data
 \cite{B1994CMP, B1994Duke, B1996CMP,B-2014, BB2014, BT2007, BT2014,BOP,BT2008-L,BT2008-G,CG,CO,DC,DJEMS,de2013, de2014, FS,HO15,HO,
 LMCPDE,LRS,NPS,NO,OSIAM,OhCMP,Oh,OhDIE,OP,OQ2013,OR,ORT, P,Poiret12,PRT,R, ST, T, TT,TV,TV2013,TV2014,
 T2009,ZF2011,ZF2012}.

 Richards \cite{R} studied the  invariance of the Gibbs measure for the
  periodic
 quartic gKdV. Oh  et al. \cite{OR}  studied  the  invariant Gibbs measure
  for the
generalized KdV equations. Hwang and  Kwak \cite{HK} established
  the probabilistic well-posedness of the generalized KdV
in $H^{s}(\R)$ with  $(s>{\rm max}\left\{\frac{1}{k+1}\left(\frac{1}{2}-\frac{2}{k}\right),
\frac{1}{4}-\frac{2}{k}\right\},k\geq5)$. Oh \cite{OhCMP}  proved the invariance of the mean 0 white noise for the periodic KdV.
Killip et al.\cite{KMV} considered the Korteweg-de Vries equation with white noise initial
data, posed on the whole real line, and proved the almost sure existence of
solutions and  the solutions obey the group property and
follow a white noise law at all times and  presented  a new proof of the existence of solutions and the
invariance of white noise measure in the torus setting. Recently, Compaan et al. \cite{CLS}
 studied the  question of  pointwise almost everywhere limit
nonlinear Schr\"odinger flows to the initial data in both continuous and periodic settings.  Linares
and Ramos \cite{LR}  established the pointwise convergence results for the flow of the generalized
 Zaharov-Kuznetsov equation in any dimension.

In this paper,  we investigate (\ref{1.01}) with rough data and random data.
Firstly,  by using the Fourier restriction norm method introduced in   \cite{Bourgain-S,Bourgain-GAFA93}
 developed in \cite{KPV-Duke,KPV1996} and
the frequency truncated technique as well as the idea of (4) in \cite{CLS},
  we prove that $u(x,t)\longrightarrow u(x,0)$ as $t\longrightarrow0$
  for a.e. $x\in \R$ with $u(x,0)\in H^{s}(\R)(s>\frac{1}{2}-\frac{2}{k},k\geq8),$
  which is just Theorem 1.1.
  Secondly,  by using Theorem 1.1 and  the  idea of Lemma2.3 in
  \cite{D},
  we  prove that
   $u(x,t)\longrightarrow  e^{-t\partial_{x}^{3}}f$ as $t\longrightarrow0$
  for a.e. $x\in \R$ with $u(x,0)\in H^{s}(\R)(s>\frac{1}{2}-\frac{2}{k},k\geq8),$
  which is just Theorem 1.2.
   Thirdly, by using the idea of \cite{Compaan},
  we prove that  $
\lim\limits_{t\longrightarrow 0}\left\|u(x,t)-e^{-t\partial_{x}^{3}}f\right\|_{L_{x}^{\infty}}=0$
  with   $f\in H^{s}(\R)(s>\frac{1}{2}-\frac{2}{k+1},k\geq5)$,
which is just Theorem 1.3.
Fourthly,
by using the Fourier restriction norm method, Strichartz estimates,
Wiener randomization of the initial data,
   probabilistic  Strichartz estimates,    we  establish
 the probabilistic
well-posedness in $H^{s}(\R)\left(s>{\rm max}
\left\{\frac{1}{k+1}\left(\frac{1}{2}-\frac{2}{k}\right),
\frac{1}{6}-\frac{2}{k}\right\}\right)$ with random data,  which is just Theorem 1.5.
 Fifthly,
  we prove that
  \begin{eqnarray*}
 \forall \epsilon>0,\forall \omega \in \Omega_{T},\lim\limits_{t\longrightarrow0}\left\|u(x,t)-e^{-t\partial_{x}^{3}}f^{\omega}(x)\right\|_{L_{x}^{\infty}}=0
  \end{eqnarray*}
   with $f(x)\in H^{s}(\R)(s>\frac{1}{6},k\geq6),$ where ${\rm P}(\Omega_{T})\geq 1- C_{1}{\rm exp}
\left(-\frac{C}{T^{\frac{\epsilon}{48k}}\|f\|_{H^{s}}^{2}}\right)$
    and $f^{\omega}(x)$  is the randomization of   $f(x)$,  which is just Theorem 1.6.
   Finally, we prove that $\forall \epsilon>0,$ $\forall \omega \in \Omega_{T},$
  $\lim\limits_{t\longrightarrow0}\left\|u(x,t)-f^{\omega}(x)\right\|_{L_{x}^{\infty}}=0$
    with  ${\rm P}(\Omega_{T})\geq 1- C_{1}{\rm exp}
\left(-\frac{C}{T^{\frac{\epsilon}{48k}}\|f\|_{H^{s}}^{2}}\right)$ and $u(x,0)\in H^{s}(\R)(s>\frac{1}{6},k\geq6),$
which is just Theorem 1.7.

We give some notations before presenting the main results.
$a\simeq b$ means that there exists two positive constants $C_{1},C_{2}$ such that $C_{1}|a|\leq |b|\leq C_{2}|a|.$
We define $\langle\cdot\rangle=1+|\cdot|.$
Let $\phi,\psi_{1}$ be smooth jump function such that $\phi(\xi)=1$
 for $|\xi|\leq 1$ and $\phi(\xi)=0$  for $|\xi|>2$ as well as $\psi_{1}(\xi)=1$
 for $|\xi|\leq 2$ and $\psi_{1}(\xi)=0$  for $|\xi|>4$.
 Then, we define
for every dyadic integer $N\in 2^{\z}$,
\begin{eqnarray*}
&&\mathscr{F}_{x}P_{N}f(\xi)=\left[\phi\left(\frac{\xi}{N}\right)-
\phi\left(\frac{2\xi}{N}\right)\right]\mathscr{F}_{x}f(\xi),\\
&&\mathscr{F}_{x}P_{\leq N}f(\xi)=\phi\left(\frac{\xi}{N}\right)
\mathscr{F}_{x}f(\xi),\\
&&\mathscr{F}_{x}P_{>N}f(\xi)=\left[1-
\phi\left(\frac{\xi}{N}\right)\right]\mathscr{F}_{x}f(\xi),\\
&&\mathscr{F}_{x}f(\xi)=\frac{1}{\sqrt{2\pi}}\int_{\SR}e^{-ix \xi}f(x)dx,\\
&&\mathscr{F}_{x}^{-1}f(\xi)=\frac{1}{\sqrt{2\pi}}\int_{\SR}e^{ix \xi}f(x)dx,\\
&&\mathscr{F}u(\xi,\tau)=\frac{1}{2\pi}\int_{\SR^{2}}e^{-ix\xi-it\tau}u(x,t)dxdt,\\
&&\mathscr{F}^{-1}u(\xi,\tau)=\frac{1}{2\pi}\int_{\SR^{2}}e^{ix\xi+it\tau}u(x,t)dxdt,\\
&&D_{x}^{\alpha}u_{0}=\frac{1}{\sqrt{2\pi}}\int_{\SR} |\xi|^{\alpha}e^{ix\xi}
\mathscr{F}_{x}{u}_{0}(\xi)d\xi,\\
&&J_{x}^{\alpha}u_{0}=\frac{1}{\sqrt{2\pi}}\int_{\SR} \langle\xi\rangle^{\alpha}e^{ix\xi}
\mathscr{F}_{x}{u}_{0}(\xi)d\xi,
\end{eqnarray*}
\begin{eqnarray*}
&&J_{t}^{\alpha}U(t)u_{0}=\frac{1}{\sqrt{2\pi}}\int_{\SR}
\left\langle\xi^{3}\right\rangle^{\alpha}e^{ix\xi+it\xi^{3}}
\mathscr{F}_{x}{u}_{0}(\xi)d\xi,\\
&&\|f\|_{L_{t}^{p}L_{x}^{q}}=\left(\int_{\SR}
\left(\int_{\SR}|f(x,t)|^{q}dx\right)^{\frac{p}{q}}dt\right)^{\frac{1}{p}},\\
&&\|a_{k}\|_{\ell_{k}^q}=\left[\sum_{k}|a_{k}|^{q}\right]^{\frac{1}{q}},\\
&&\|P_{N}f\|_{\ell_{N}^{r}L_{t}^{p}L_{x}^{q}}=\left[\sum_{N}\left(\int_{\SR}
\left(\int_{\SR}|P_{N}f(x,t)|^{q}dx\right)^{\frac{p}{q}}dt\right)^{\frac{r}{p}}\right]^{\frac{1}{r}},\\
&&\|P_{N}f\|_{L_{t}^{p}L_{x}^{q}\ell_{N}^{r}}=
\left\|\left[\sum_{N}|P_{N}f|^{r}\right]^{\frac{1}{r}}\right\|_{L_{t}^{p}L_{x}^{q}},\\
&&\|f\|_{L_{xt}^{p}}=\|f\|_{L_{x}^{p}L_{t}^{p}}.
\end{eqnarray*}
$H^{s}(\R)=\left\{f\in \mathscr{S}^{'}(\R):\|f \|_{H^{s}(\SR)}=
 \|\langle\xi\rangle ^{s}\mathscr{F}_{x}{f}\|_{L_{\xi}^{2}(\SR)}<\infty\right\}$.
 The space $X_{s,b}(\R^{2})$ is defined as follows:
\begin{eqnarray*}
X_{s,b}(\R^{2})=\left\{u\in \mathscr{S}^{'}(\R^{2}):\|u\|_{X_{s,b}}
=\left[\int_{\SR^{2}}\langle \xi\rangle^{2s}\langle \sigma\rangle^{2b}
|\mathscr{F}u(\xi,\tau)|^{2}d\xi\right]^{\frac{1}{2}}<\infty\right\}.
\end{eqnarray*}
Here, $\langle \sigma\rangle =1+|\tau-\xi^{3}|.$

Now we present the pointwise convergence problem and uniform convergence problem related to (\ref{1.01})-(\ref{1.02}), which are just Theorems 1.1-1.3.

\begin{Theorem} \label{Theorem1}
Let  $k\geq8$ and $f\in H^{s}(\R),s>\frac{1}{2}-\frac{2}{k}$. Then,   we have
\begin{eqnarray*}
u(x,t)\longrightarrow f
\end{eqnarray*}
as $t\longrightarrow 0$ for a.e. $x\in \R$.
\end{Theorem}

\noindent {\bf Remark1:}We present the outline of proof of Theorem 1.1. From (3.9) of  \cite{KPVCPAM}, we have
\begin{eqnarray}
\|U(t)f\|_{L_{x}^{4}L_{t}^{\infty}}\leq C\|f\|_{H^{\frac{1}{4}}(\SR)}\label{1.03}.
\end{eqnarray}
Combining (\ref{1.03}) with a standard proof of \cite{KPV-Duke}, we have
\begin{eqnarray}
\|u\|_{L_{x}^{4}L_{t}^{\infty}}\leq C\|u\|_{X_{\frac{1}{4},\frac{1}{2}+\epsilon}}\label{1.04}.
\end{eqnarray}
We consider
\begin{eqnarray}
&&(u^{N})_{t}+(u^{N})_{xxx}+\frac{1}{k+1}\partial_{x}P_{\leq N}((u^{N})^{k+1})=0\label{1.05},\\
&&
u^{N}(x,0)=P_{\leq N}u(x,0)=P_{\leq N}f.\label{1.06}
\end{eqnarray}
Here,  $P_{\leq N}u(x,0)=\frac{1}{2\pi}\int_{\SR}\phi\left(\frac{\xi}{N}\right)e^{ix\xi}\mathscr{F}_{x}u(\xi,0)d\xi
=\frac{1}{2\pi}\int_{\SR}\phi\left(\frac{\xi}{N}\right)e^{ix\xi}\mathscr{F}_{x}f(\xi)d\xi.$
We define $u:=\lim\limits_{N\longrightarrow \infty}u^{N}.$
For $f\in H^{s}(\R)(s>\frac{1}{2}-\frac{2}{k},k\geq8),$  by using the Fourier restriction
 norm method introduced in   \cite{Bourgain-S,Bourgain-GAFA93}
 developed in \cite{KPV-Duke,KPV1996} and
the frequency truncated technique as well as the idea of (4) in \cite{CLS}, we firstly prove
\begin{eqnarray}
\lim\limits_{N\longrightarrow \infty}
\left\|\sup\limits_{0\leq t\leq \delta}\left|u-u^{N}\right|\right\|_{L_{x}^{4}}=0,\label{1.07}
\end{eqnarray}
which is just Lemma 4.1 in this paper.
Since $P_{\leq N}f\in H^{\infty}(\R)$, then, we have
\begin{eqnarray}
\lim\limits_{t\longrightarrow0}u^{N}=P_{\leq N}f\label{1.08}
\end{eqnarray}
for a.e. $x\in \R$. Since
\begin{eqnarray}
\left|u-f\right|\leq |u-u^{N}|+|u^{N}-P_{\leq N}f|+|P_{>N}f|,\label{1.09}
\end{eqnarray}
by using (\ref{1.09}), we have
\begin{eqnarray}
&&\lim\limits_{t\longrightarrow0}\sup\left|u-f\right|\leq \lim\limits_{t\longrightarrow0}\sup|u-u^{N}|
+\lim\limits_{t\longrightarrow0}\sup|u^{N}-P_{\leq N}f|+\lim\limits_{t\longrightarrow0}\sup|P_{>N}f|\nonumber\\
&&\leq \lim\limits_{t\longrightarrow0}\sup|u-u^{N}|+\lim\limits_{t\longrightarrow0}\sup|P_{>N}f|.\label{1.010}
\end{eqnarray}
By using the Chebyshev inequality and the Sobolev embeddings theorem $H^{\frac{1}{4}}(\R)\hookrightarrow L^{4}(\R)$
 as well as Lemma 4.1, since $f\in H^{s}(\R)$,  we have
\begin{eqnarray}
&&\left|\left\{x\in B_{1}: \lim\limits_{t\longrightarrow0}\sup\left|u-f\right|>\alpha\right\}\right|\leq
 \left|\left\{x\in B_{1}: \lim\limits_{t\longrightarrow0}\sup\left|u-u^{N}\right|>\frac{\alpha}{2}\right\}\right|
\nonumber\\&&\qquad\qquad+\left|\left\{x\in B_{1}: \lim\limits_{t\longrightarrow0}\sup\left|P_{>N}f\right|>\frac{\alpha}{2}\right\}\right|\nonumber\\
&&\leq C\alpha^{-4}\left\|\sup\limits_{0\leq t\leq \delta}\left|u-u^{N}\right|\right\|_{L_{x}^{4}}+C\alpha^{-4}\|P_{>N}f\|_{L_{x}^{4}}\nonumber\\
&&\leq C\alpha^{-4}\left\|\sup\limits_{0\leq t\leq \delta}\left|u-u^{N}\right|\right\|_{L_{x}^{4}}+C\alpha^{-4}\|P_{>N}f\|_{H_{x}^{\frac{1}{4}}}\nonumber\\
&&\leq C\alpha^{-4}\left\|\sup\limits_{0\leq t\leq \delta}\left|u-u^{N}\right|\right\|_{L_{x}^{4}}+C\alpha^{-4}\|P_{>N}f\|_{H_{x}^{s}}\label{1.011}.
\end{eqnarray}
Here $B_{1}\subset \R.$
From Lemma 4.1 and $f\in H^{s}(\R),$ we have
\begin{eqnarray}
\left|\left\{x\in B_{1}: \lim\limits_{t\longrightarrow0}\sup\left|u-f\right|>\alpha\right\}\right|\leq \left|\left\{x\in B_{1}:
\lim\limits_{t\longrightarrow0}\sup\left|u-u^{N}\right|>\frac{\alpha}{2}\right\}\right|=0.\label{1.012}
\end{eqnarray}

\begin{Theorem} \label{Theorem2}
Let  $k\geq8$ and $f\in H^{s}(\R),s>\frac{1}{2}-\frac{2}{k}$. Then,   we have
\begin{eqnarray*}
u(x,t)\longrightarrow U(t)f.
\end{eqnarray*}
as $t\longrightarrow 0$ for a.e. $x\in \R$.
\end{Theorem}
\noindent {\bf Remark2:} We present the outline of proof of Theorem 1.2.
 By using (\ref{1.012}) and a proof similar to
  Lemma 2.3 of \cite{D}, since $f\in H^{s}(\R)(s>\frac{1}{2}-\frac{2}{k}\geq \frac{1}{4},k\geq8)$ which yields
\begin{eqnarray}
\|U(t)f\|_{L_{x}^{4}L_{t}^{\infty}}\leq C\|f\|_{H^{\frac{1}{4}}}\leq C\|f\|_{H^{s}}\label{1.013},
\end{eqnarray}
we have
\begin{eqnarray}
\lim\limits_{t\longrightarrow0}\sup\left|u-U(t)f\right|\leq \lim\limits_{t\longrightarrow0}
\sup|u-f|+\lim\limits_{t\longrightarrow0}\sup|f-U(t)f|=0.\label{1.014}
\end{eqnarray}

\begin{Theorem} \label{Theorem3}
Let  $k\geq5$ and $f\in H^{s}(\R),s>\frac{1}{2}-\frac{2}{k+1}$. Then,   we have
\begin{eqnarray*}
\lim\limits_{t\longrightarrow 0}\left\|u(x,t)-U(t)f\right\|_{L_{x}^{\infty}}=0.
\end{eqnarray*}
\end{Theorem}
\noindent {\bf Remark3:}We present the outline of proof of Theorem 1.3.
We firstly prove  that  there exists a unique solution $u\in X_{s,b}(b>\frac{1}{2})$ to (\ref{1.01})-(\ref{1.02})
 for $f\in H^{s}(\R)(s>\frac{1}{2}-\frac{2}{k},k\geq5)$, which is just Lemma 6.1.
Then, by using  Lemma 6.1 and Lemma 3.2 and
$X_{\frac{1}{2}+\epsilon,b}\hookrightarrow C([0,T];H^{\frac{1}{2}+\epsilon})$ as well as (\ref{6.004}),  we have
\begin{eqnarray}
&&\left\|u-\phi(t)U(t)f\right\|_{L_{x}^{\infty}}=\left\|u-U(t)f\right\|_{L_{x}^{\infty}}
\leq C\left\|u-U(t)f\right\|_{C([0,T];H_{x}^{\frac{1}{2}+\epsilon})}\nonumber\\&&= C\left\|u-U(t)f\right\|_{X_{\frac{1}{2}+\epsilon,b}}\nonumber\\&&=
C\left\|\phi\left(\frac{t}{T}\right)\int_{\SR}U(t-\tau)\partial_{x}\left[(u)^{k+1}\right]\right\|_{X_{\frac{1}{2}+\epsilon,b}}
\nonumber\\&&\leq CT^{\epsilon}\|\partial_{x}\left[u^{k+1}\right]\|_{X_{s,b_{1}}}\nonumber\\&&\leq  CT^{\epsilon}\|u\|_{X_{s,b}}^{k+1}
\leq CT^{\epsilon}\|f\|_{H^{s}}^{k+1}<\infty.\label{1.015}
\end{eqnarray}
From (\ref{1.015}), we have
\begin{eqnarray}
\lim\limits_{t\longrightarrow 0}\left\|u-U(t)f\right\|_{L_{x}^{\infty}}=0.\label{1.016}
\end{eqnarray}

To present Theorems 1.4-1.7,  now we introduce the randomization procedure for the initial data, which can
 be seen in \cite{BOP,BOP2015,LMCPDE,ZF2011}.
 Let $\psi\in\mathcal{S}(\R)$ be an even, non-negative jump function
with $supp(\psi)\subseteq [-1,1]$ and such that for $\xi\in\R$,
 \begin{eqnarray}
\sum_{k\in
\z}\psi(\xi-k)=1.\label{1.017}
\end{eqnarray}
For every $k\in
\Z$, we define the function $\psi(D-k)f:\R\rightarrow\mathbb{C}$ by
\begin{eqnarray*}
(\psi(D-k)f)(x)=\mathscr{F}^{-1}_{\xi}\big(\psi(\xi-k)\mathscr{F}_xf\big)(x), ~x\in \R.
\end{eqnarray*}
Note that these projections satisfy
 a unit-scale Bernstein inequality,
 namely that,
 for any $p_1,p_2$, which satisfies that
$2 \leq p_1 \leq p_2 \leq\infty$,
there exists a $C= C(p_1, p_2)>0$ such that for any $f\in L^2(\R)$
and $k\in \Z$,
 \begin{eqnarray}
\left\|\psi(D-k)f\right\|_{L_{x}^{p_2}(\SR)}\leq C
\left\|\psi(D-k)f\right\|_{L_{x}^{p_1}(\SR)}\leq C
\left\|\psi(D-k)f\right\|_{L_{x}^{2}(\SR)}.\label{1.018}
\end{eqnarray}

Let $\{g_k\}_{k\in Z}$ be a sequence of independent, zero-mean, complex-valued
 Gaussian random variables
 on a probability space $(\Omega,\mathcal{A}, \mathbb{P})$,
 where the real and imaginary parts of $g_k$ are independent
  and endowed with probability
distributions $\mu_k^1$ and $\mu_k^2$ respectively. The probability
distributions $\mu_k^1$ and $\mu_k^2$ satisfies
 \begin{eqnarray}
\Big|\int_{-\infty}^{+\infty}e^{\gamma x}d\mu_k^j(x)\Big|\leq e^{c\gamma^2},\quad\text{for all}\quad
 \gamma\in \R, k\in  \Z, j=1, 2.\label{1.019}
\end{eqnarray}
Thereafter for a given $f\in H^{s}(\R)$, we define its
randomization by
\begin{eqnarray}
f^\omega:=\sum_{k\in \z}g_k(\omega)\psi(D-k)f.\label{1.020}
\end{eqnarray}
We define
\begin{eqnarray*}
\|f_1\|_{L_{\omega}^{p}(\Omega)}=\left[\int_{\Omega}|f_{1}(\omega)|^{p}d\mathbb{P}(\omega)\right]^{\frac{1}{p}}.
\end{eqnarray*}
If $f\in H^{s}(\R)$,
then  $f^\omega$
is almost surely
in $H^{s}(\R)$ and $\|\|f^\omega\|_{H^{s}}\|_{L_{\omega}^{2}}\simeq\|f\|_{H^{s}}$,
 see Lemma 2.2 in \cite{BOP}. This randomization improved the
 integrability of $f$, see Lemma 2.3 of \cite{BOP}. Such results for random
 Fourier series are known as Paley-Zygmund's theorem \cite{PZ}.

\begin{Theorem}\label{Thm4}
Let $k\geq 5,s\geq \frac{1}{2}-\frac{2}{k}+88\epsilon$ and $s_{1}\geq{\rm max}
\left\{\frac{s}{k+1},s-\frac{1}{3},
\frac{s-\frac{1}{3}}{k}\right\}+88\epsilon$, $z=\phi(t)S(t)f^{\omega}$
 and   $b=\frac{1}{2}+\frac{\epsilon}{24}$,   we denote by $f^{\omega}$ the randomization of $f$ as
defined in (\ref{1.06}).
  Then, we have
\begin{eqnarray}
\left\|\partial_{x}\left[(v+z)^{k+1}\right]\right\|_{X_{s,-\frac{1}{2}+\frac{\epsilon}{12}}}\leq C
\left[\|v\|_{X_{s,b}}^{k+1}+\lambda^{k+1}\right]\label{1.021}
\end{eqnarray}
outside a set of probability
$C_{1}{\rm exp}
\left(-\frac{C}{T^{\frac{\epsilon}{48k}}\|f\|_{H^{s}}^{2}}\right).$
\end{Theorem}
\noindent {\bf Remark4:} We present the outline of proof of Theorem 1.4.
We use   Lemmas 2.2, 2.3, 2.6-2.8, 2.13-2.15 to  establish Theorem 1.4. The establishment of Lemma 3.3 is most difficult in  proving Theorem 1.4.
Thus, we use Lemmas 2.7-2.8, 2.14 introduced in this paper to overcome the difficulty.
 Hwang and Kwak
\cite{HK} proved that when $k\geq5,
{\rm max}\left\{\frac{1}{k+1}(\frac{1}{2}-\frac{2}{k}),
\frac{1}{4}-\frac{2}{k}\right\}<s_{1}<\frac{1}{2}-\frac{2}{k},$
$\frac{1}{2}-\frac{3}{2(k-1)}<s<{\rm min}\left\{(k+1)s_{1},s_{1}+\frac{1}{4}\right\}$,
(\ref{1.07}) is
valid. The main tools used in \cite{HK} to  establish multilinear
estimates are Proposition 2.4 called
probabilistic local smoothing estimate and Lemma 3.5 of \cite{HK}.
Thus, we improve the result of  Hwang and Kwak
\cite{HK}.

\begin{Theorem}\label{Thm5}
Let $k\geq5,s_{1}>{\rm max}\left\{\frac{1}{k+1}\left(\frac{1}{2}-\frac{2}{k}\right),
\frac{1}{6}-\frac{2}{k}\right\}$, $f\in H^{s_{1}}(\R)$ and $f^{\omega}$ be its
randomization defined in (\ref{1.020}) and
\begin{eqnarray}
\left|\int_{\SR}e^{\gamma x}d\mu_{n}^{(j)}(x)\right|\leq e^{C\gamma^{2}}\label{1.022}
\end{eqnarray}
for all $\gamma\in \R,n\in Z,j=1,2.$ Then  (\ref{1.01}) on $\R$ is probabilistically   well-posed with random  data
$f^{\omega}$. More  precisely, there exist $C,C_{1}>0$ and
$s>\frac{1}{2}-\frac{2}{k}$
such that for each $T\ll1$,
there exists a set $\Omega_{T}\subset \Omega$ with the following properties:

\noindent {\rm (i)}:$P\left(\Omega\setminus\Omega_{T}\right)\leq C_{1}{\rm exp}
\left(-\frac{C}{T^{\frac{\epsilon}{48k}}\|f\|_{H^{s_{1}}}^{2}}\right).$

\noindent {\rm (ii)}: For each $\omega\in \Omega_{T},$ there exists a unique
solution to (\ref{1.01})
with $u(x,0)=f^{\omega}$ in the class
$
S(t)f^{\omega}+C([-T,T]:H^{s}(\R))\subset C([-T,T]:H^{s_{1}}(\R)).
$

\end{Theorem}

\noindent {\bf Remark5:} Combining Theorem 1.4 with  the fixed point argument, we derive Theorem 1.3.
In Theorem  1.3  of \cite{HK},  Hwang and Kwak
 proved that when $f \in H^{s_{1}}(\R)$,
${\rm max}\left\{\frac{1}{k+1}(\frac{1}{2}-\frac{2}{k}),
\frac{1}{4}-\frac{2}{k}\right\}<s_{1}<\frac{1}{2}-\frac{2}{k},$
$\frac{1}{2}-\frac{3}{2(k-1)}<s<{\rm min}\left\{(k+1)s_{1},s_{1}+\frac{1}{4}\right\}$,
 (\ref{1.01}) is
probabilistically well-posed with randomized initial data
 $f^{\omega}.$ Thus, the result of Theorem 1.5 improves the result
 of Theorem  1.3  of \cite{HK}.

\begin{Theorem} \label{Theorem6}
Let  $k\geq6$ and $f\in H^{s_{1}}(\R),s_{1}>\frac{1}{6}$, we denote by $f^{\omega}$ the randomization of $f$ as
defined in (\ref{1.020}). Then,  $\forall \omega\in \Omega_{T}$,  we have
\begin{eqnarray*}
\lim\limits_{t\longrightarrow0}\left\|u(x,t)-U(t)f^{\omega}(x)\right\|_{L_{x}^{\infty}}=0.
\end{eqnarray*}
Here, ${\rm P}(\Omega_{T})\geq 1- C_{1}{\rm exp}
\left(-\frac{C}{T^{\frac{\epsilon}{48k}}\|f\|_{H^{s_{1}}}^{2}}\right)$.
\end{Theorem}
\noindent {\bf Remark6:} We present the outline of proof of Theorem 1.6.Let $v=u-U(t)f^{\omega}.$ Then,
 $\forall \omega \in \Omega_{T}$,  by using  Theorems 1.4, 1.5  and  $X_{\frac{1}{2}+\epsilon,b}\hookrightarrow C([0,T];H^{\frac{1}{2}+\epsilon})$,  we have
\begin{eqnarray}
&&\left\|u-U(t)f^{\omega}\right\|_{L_{x}^{\infty}}=\left\|v\right\|_{L_{x}^{\infty}}
\leq C\left\|v\right\|_{C([0,T];H_{x}^{\frac{1}{2}+\epsilon})}\leq
C\left\|v\right\|_{X_{\frac{1}{2}+\epsilon,b}}\nonumber\\&&\leq C
\left\|\int_{\SR}U(t-\tau)\partial_{x}\left[(v+U(t)f^{\omega})^{k+1}\right]\right\|_{X_{\frac{1}{2}+\epsilon,b}}
\leq C\left[\|v\|_{X_{s,b}}^{k+1}+\lambda^{k+1}\right]\nonumber\\
&&\leq C\left[1+\lambda^{k+1}\right].\label{1.023}
\end{eqnarray}
From (\ref{1.023}) and $v(0)=0$, we have
\begin{eqnarray}
\lim\limits_{t\longrightarrow 0}\|v(t)\|_{H^{s}}=\lim\limits_{t\longrightarrow 0}\left[\|v(t)\|_{H^{s}}-\|v(0)\|_{H^{s}}\right]=0.\label{1.024}
\end{eqnarray}
From (\ref{1.024}) and $H^{s}(\R)\hookrightarrow L^{\infty}(\R)(s>\frac{1}{2})$, we have
\begin{eqnarray}
\lim\limits_{t\longrightarrow 0}\left\|u-U(t)f^{\omega}\right\|_{L_{x}^{\infty}}
=\lim\limits_{t\longrightarrow 0}\left\|v\right\|_{L_{x}^{\infty}}
\leq C\lim\limits_{t\longrightarrow 0}\|v(t)\|_{H^{s}}=0.\label{1.025}
\end{eqnarray}
From (\ref{1.025}), we have
\begin{eqnarray}
\lim\limits_{t\longrightarrow 0}\left\|u-U(t)f^{\omega}\right\|_{L_{x}^{\infty}}=0.\label{1.026}
\end{eqnarray}
From Theorems 1.3 and 1.6, we know that the uniform convergence problem of  (\ref{1.01})
 with  random data requires less regularity of the initial data than
the uniform convergence problem of (\ref{1.01})  with  rough data.

\begin{Theorem} \label{Theorem7}
Let  $k\geq6$ and $f\in H^{s_{1}}(\R),s_{1}>\frac{1}{6}$, we denote by $f^{\omega}$ the randomization of $f$ as
defined in (\ref{1.020}). Then,  $\forall \omega \in \Omega_{T}$,  we have
\begin{eqnarray*}
\lim\limits_{t\longrightarrow0}\left\|u(x,t)-f^{\omega}(x)\right\|_{L_{x}^{\infty}}=0.
\end{eqnarray*}
Here, ${\rm P}(\Omega_{T})\geq 1- C_{1}{\rm exp}
\left(-\frac{C}{T^{\frac{\epsilon}{48k}}\|f\|_{H^{s_{1}}}^{2}}\right)$.
\end{Theorem}
\noindent {\bf Remark7:}
Combining Theorem 1.5 with  Theorem 1.1 of \cite{YDLY},  $\forall \omega \in \Omega_{T},$ we have
\begin{eqnarray}
\lim\limits_{t\longrightarrow0}\sup\left|u-f^{\omega}\right|\leq \lim\limits_{t\longrightarrow0}
\sup|u-U(t)f^{\omega}|+\lim\limits_{t\longrightarrow0}\sup|f^{\omega}-U(t)f^{\omega}|=0.\label{1.027}
\end{eqnarray}

The rest of the paper is arranged as follows. In Section 2,  we give some
preliminaries. In Section 3, we show some multilinear estimates. In Section 4, we prove
 Theorem 1.1. In Section 5, we prove
 Theorem 1.2. In Section 6, we prove
 Theorem 1.3. In Section 7, we prove
 Theorem 1.4. In Section 8, we prove
 Theorem 1.5. In Section 9, we prove
 Theorem 1.6. In Section 10, we prove
 Theorem 1.7.

\bigskip

 \noindent{\large\bf 2. Preliminaries }

\setcounter{equation}{0}

\setcounter{Theorem}{0}

\setcounter{Lemma}{0}

\setcounter{section}{2}

In this section, we present the  probabilistic Strichartz estimates,
 linear estimates and Strichartz estimates used in this paper.

\begin{Lemma}\label{Lemma2.1}
Let $T\in (0,1)$ and $s\in \R$ and $-\frac{1}{2}<b^{\prime}
\leq0\leq b\leq b^{\prime}+1$.
Then, for $h\in X_{s,b^{\prime}},$  we have
\begin{eqnarray}
&&\left\|\phi(t)U(t)f\right\|_{X_{s,\frac{1}{2}+\epsilon}}
\leq C\|f\|_{H^{s}},\label{2.01}\\
&&\left\|\phi\left(\frac{t}{T}\right)\int_{0}^{t}U(t-\tau)h(\tau)
d\tau\right\|_{X_{s,b}}\leq C
T^{1+b^{\prime}-b}\|h\|_{X_{s,b^{\prime}}}.\label{2.02}
\end{eqnarray}
\end{Lemma}

For the proof of Lemma 2.1, we refer the readers to
 \cite{Bourgain-GAFA93,G}.

\begin{Lemma}\label{Lemma2.2}
\noindent Let $b=\frac{1}{2}+\frac{\epsilon}{24}$.
Then,  for $0\leq s\leq \frac{1}{2},$  we have
\begin{eqnarray}
\left\|I^{s}(u_{1},u_{2})\right\|
  _{L_{xt}^{2}}\leq C\prod_{j=1}^{2}\| u_{j}\| _{X_{0,
  \frac{3+2s}{4}b}},\label{2.03}
\end{eqnarray}
where
\begin{eqnarray*}
\mathscr{F}I^{s}(u_{1},u_{2})(\xi,\tau)
&=&\int_{\!\!\!\mbox{\scriptsize $
\begin{array}{l}
\xi=\xi_{1}+\xi_{2}\\
\tau=\tau_{1}+\tau_{2}
\end{array}
$}}
\left|\xi_{1}^{2}-\xi_{2}^{2}\right|^{s}
\mathscr{F}{u_{1}}(\xi_{1},\tau_{1})\mathscr{F}{u_{2}}
(\xi_{2},\tau_{2})\,d\xi_{1}d\tau_{1}.
\end{eqnarray*}

\end{Lemma}

For the proof of Lemma 2.2, we refer the readers to \cite{G}.

\begin{Lemma}\label{Lemma2.3}
Let $q\geq8$ and $s=\frac{1}{8}-\frac{1}{q}$ and
 $b=\frac{1}{2}+\frac{\epsilon}{24}$
 and $0\leq\epsilon\leq 10^{-3}$. Then, we have
\begin{eqnarray}
&&\|U(t)u_{0}\|_{L_{xt}^{8}}\leq C\|u_{0}\|_{L_{x}^{2}},\label{2.04}\\
&&\|u\|_{L_{xt}^{q}}\leq C\|u\|_{X_{4s,b}},\label{2.05}\\
&&\|D_{x}^{\frac{1}{6}}u\|_{L_{xt}^{6}}\leq C\|u\|_{X_{0,b}},\label{2.06}\\
&&\|u\|_{L_{xt}^{8}}\leq C\|u\|_{X_{0,b}}\label{2.07},\\
&&\|u\|_{L_{xt}^{\frac{8}{1+\epsilon}}}\leq C\|u\|_{X_{0,\frac{3-\epsilon}{3}b}}
\leq C\|u\|_{X_{0,\frac{1}{2}-\frac{\epsilon}{12}}}\label{2.08},\\
&&\|D_{x}^{\frac{3-\epsilon}{18}}u\|_{L_{xt}^{\frac{6}{1+\epsilon}}}\leq C
\|u\|_{X_{0,\frac{3-\epsilon}{3}b}}\leq C\|u\|_{X_{0,\frac{1}{2}-\frac{\epsilon}{12}}}
\label{2.09},\\
&&\|D_{x}^{\frac{1-8\epsilon}{6}}u\|_{L_{xt}^{\frac{6}{1-2\epsilon}}}\leq C
\|u\|_{X_{0,b}}\label{2.010},\\
&&\|D_{x}^{1-2\epsilon}u\|_{L_{x}^{\frac{1}{\epsilon}}L_{t}^{2}}\leq C
\|u\|_{X_{0,(1-2\epsilon)b}}\label{2.011}.
\end{eqnarray}
\end{Lemma}
\noindent {\bf Proof.}(\ref{2.04}), (\ref{2.06}),  (\ref{2.07}) can be seen in
 Lemma 2.1 of \cite{KPV-Duke}.
  By using the Sobolev embeddings theorem and (\ref{2.04}), we have
\begin{eqnarray}
&&\|U(t)u_{0}\|_{L_{xt}^{q}}\leq C\|D_{x}^{s}D_{t}^{s}U(t)u_{0}\|_{L_{xt}^{8}}
=C\left\|\int_{\SR}e^{ix\xi+it\xi^{3}}|\xi|^{4s}\mathscr{F}_{x}u_{0}(\xi)d\xi\right\|_{L_{xt}^{8}}
\nonumber\\&&\leq C\|u_{0}\|_{H^{4s}}
\label{2.012}.
\end{eqnarray}
Obviously, changing variable $\tau=\lambda+\xi^{3}$, we have
 \begin{eqnarray}
&&u(x,t)=\frac{1}{2\pi}\int_{\SR^{2}}e^{ix\xi+it\tau}\mathscr{F}u(\xi,\tau)d\xi d\tau\nonumber\\
&&=\frac{1}{2\pi}\int_{\SR^{2}}e^{ix\xi+it(\lambda+\xi^{3})}\mathscr{F}
u(\xi,\lambda+\xi^{3}d\xi d\lambda\nonumber\\
&&=\frac{1}{2\pi}\int_{\SR}e^{it\lambda}\left(\int_{\SR}e^{ix\xi+it\xi^{3}}
\mathscr{F}u(\xi,\lambda+\xi^{3}) d\xi\right)d\lambda\label{2.013}.
 \end{eqnarray}
By using  (\ref{2.013}),  (\ref{2.012}),  Minkowski's inequality,
 taking $b>\frac{1}{2}$  and   variable substitution $\lambda+\xi^{3}=\tau,$ we have
\begin{eqnarray}
&&\|u\|_{L_{xt}^{q}}\leq C\int_{\SR}\left\|\left(\int_{\SR}e^{ix\xi+it\xi^{3}}
\mathscr{F}u(\xi,\lambda+\xi^{3})d\xi\right)\right\|_{L_{xt}^{q}} d\lambda\nonumber\\&&\leq
 C\int_{\SR}\left\|\mathscr{F}u(\xi,\lambda+\xi^{3})\right\|_{H^{4s}}d\lambda\nonumber\\&&\leq C
\left[\int_{\SR}(1+|\lambda|)^{2b}\left\|\mathscr{F}u(\xi,\lambda+\xi^{3})\right\|_{H^{4s}}^{2}
d\lambda\right]^{\frac{1}{2}}
\left[\int_{\SR}(1+|\lambda|)^{-2b}d\lambda\right]^{\frac{1}{2}}\nonumber\\
&&\leq C\left[\int_{\SR}(1+|\tau-\xi^{3}|)^{2b}\left\|\mathscr{F}u(\xi,\tau)
\right\|_{H^{4s}}^{2}d\tau\right]^{\frac{1}{2}}
=\|u\|_{X_{4s,b}}.\label{2.014}
\end{eqnarray}
Thus, we derive (\ref{2.05}).
 Interpolating (\ref{2.06}) with
\begin{eqnarray*}
\|u\|_{L_{xt}^{2}}\leq C\|u\|_{X_{0,0}}
\end{eqnarray*}
yields (\ref{2.09}). Interpolating (\ref{2.06}) with (\ref{2.07})
yields (\ref{2.010}).
   (\ref{2.011}) can be seen in Lemma 3.3 of  \cite{HK}.

We have completed the proof of
Lemma 2.3.
\bigskip

\begin{Lemma}\label{lem2.4}
Let  $\phi_{j}(j=1,2)\in C^{\infty}(\R)$,  $\supp \phi_{2}\subset (a,b)$ and  $\phi_{1}^{\prime}(\xi)\neq 0$ for all $\xi\in [a,b]$.  Then,
we have
\begin{align*}
\left|\int_{a}^{b}e^{i\lambda\phi_{1}(\xi)}\phi_{2}(\xi)d\xi\right|
\leq \frac{C}{|\lambda|^{k}}
\end{align*}
for all $k\geq 0$. Here $C=C(\phi_{1},\phi_{2},k)$.

\end{Lemma}

Lemma 2.4 can be seen in \cite{Stein-new}.
\begin{Lemma}\label{lem2.5}
Let  $\phi_{4}\in C^{\infty}(\R)$ and  $\phi_{3} \in C^{3}(\R)$and  $\supp \phi_{3} \subset (a,b)$ and
 $\left|\phi_{3}^{(3)}(\xi)\right|\geq 1$. Then, we have
\begin{align*}
\left|\int_{a}^{b}e^{i\lambda\phi_{3}(\xi)}\phi_{4}(\xi)d\xi\right|\leq
 \frac{C}{|\lambda|^{\frac{1}{3}}}\left(|\phi_{4}(b)|+\int_{a}^{b}
 |\phi_{4}^{'}(\xi)|d\xi\right).
\end{align*}
\end{Lemma}

Lemma 2.4 can be seen in \cite{Stein-new}.

\begin{Lemma}\label{lem2.6}We define
\begin{align*}
 K(x,t):=\int_{\SR}e^{it\xi^{3}+ix\xi}\psi_{1}\left(\frac{\xi}{4N}\right)d\xi.
 \end{align*}
 Then, for $\gamma\geq7$, we have
\begin{align}\label{2.015}
 \|K\|_{L_{x}^{\frac{\gamma}{2}}L_{t}^{\infty}}\leq C N^{\frac{\gamma-2}{\gamma}}.
 \end{align}
\end{Lemma}
\noindent{\bf Proof.} Obviously  $\R_{x}\times \R_{t}\subset \bigcup\limits_{j=1}^{3}\Omega_{j}$, where
\begin{eqnarray*}
&&\Omega_{1}:=\left\{(x,t)\in \R\times \R: |x|\leq \frac{1}{N}\right\},\\
&&\Omega_{2}:= \left\{(x,t)\in \R\times \R: |x|> \frac{1}{N},|x|\geq 1536N^{2}t\right\},\\
&&\Omega_{3}:= \left\{(x,t)\in \R\times \R: |x|>\frac{1}{N},|x|< 1536N^{2}t\right\}.
\end{eqnarray*}
We define $\Omega_{x,i}:=\left\{t\in \R|(x,t)\in \Omega_{i}\right\}$ for a fixed $x\in \R$.
Let $\eta=\frac{\xi}{4N}$. Thus,  we obtain
\begin{align*}
 K(x,t)=4N\int_{\SR}e^{64iN^{3}\eta^{3}t+4ixN\eta}\psi_{1}(\eta)d\eta.
 \end{align*}
From the definition of $\Omega_{x1},$  we have
 \begin{align}\label{2.016}
\left[\int_{|x|<\frac{1}{N}}\left[\sup\limits_{t\in \Omega_{x,2}}|K(x,t)|
\right]^{\frac{\gamma}{2}}dx\right]^{\frac{2}{\gamma}}\leq C
N\left[\int_{|x|\leq \frac{1}{N}}\left(\int_{\SR}\psi_{1}(\eta)d\eta\right)^{\frac{\gamma}{2}}dx\right]^{\frac{2}{\gamma}}
\leq C N^{\frac{\gamma-2}{\gamma}}.
\end{align}
 Obviously,  $64iN^{3}\eta^{3}t+4ixN\eta
=4ixN(\eta+\frac{16N^{2}\eta^{3}t}{x}):=
4ixN\phi_{5}(\eta),\phi_{5}(\eta)=\eta+\frac{16N^{2}\eta^{3}t}{x}$ yields
\begin{align*}
 |\phi_{5}^{\prime}(\eta)|=\left|1+\frac{48tN^{2}\eta^{2}}{x}\right|\geq 1-768\left|
 \frac{tN^{2}}{x}\right|\geq 1-\frac{1}{2}=\frac{1}{2}
 \end{align*}
uniformly with respect to  $(x,t)\in \Omega_{2}$. Therefore, $\phi_{5}^{\prime}\neq0$
for any  $(x,t)\in \Omega_{2}$. By Lemma 2.4, we have  $|K(x,t)|\leq C N(N|x|)^{-2}=N^{-1}x^{-2}$.
 By using  a  direct  computation, we have
 \begin{align}\label{2.017}
\left[\int_{|x|>\frac{1}{\gamma}}\left[\sup\limits_{t\in \Omega_{x,2}}|K(x,t)|
\right]^{\frac{\gamma}{2}}dx\right]^{\frac{2}{\gamma}}\leq C
 \left[\int _{|x|>\frac{1}{\gamma}}N^{-\frac{\gamma}{2}}|x|^{-\gamma}dx\right]^{\frac{2}{\gamma}}
 \leq CN^{-1}N^{2-\frac{2}{\gamma}}=N^{\frac{\gamma-2}{\gamma}}.
\end{align}
 We define
$64iN^{3}\eta^{3}t+4ixN\eta=4itN^{3}(16\eta^{3}+\frac{x\eta}{N^{2}t}):
=4itN^{3}\phi_{6}(\eta)$. By using a  direct computation, we  have  $|\phi_{6}^{(3)}(\eta)|\geq 1$.
By using  Lemma 2.5, we have
\begin{align}\label{2.018}
 |K(x,t)|=CN\left|\int_{\SR}e^{64iN^{3}\eta^{3}t+4ixN\eta}\psi_{1}(\eta)d\eta\right|
 \leq CN(N^{3}t)^{-\frac{1}{3}}
 \leq \frac{C}{t^{\frac{1}{3}}}\leq C\frac{N^{\frac{2}{3}}}{|x|^{\frac{1}{3}}}.
 \end{align}
 From (\ref{2.018}), for $\gamma \geq 7,$  we have
 \begin{align}\label{2.019}
\left[\int_{|x|>\frac{1}{N}}\left[\sup\limits_{t\in \Omega_{x,3}}
|K(x,t)|^{\frac{\gamma}{2}}\right]dx\right]^{\frac{2}{\gamma}}\leq C
N^{\frac{2}{3}} \left[\int_{|x|>\frac{1}{N}}|x|^{-\frac{\gamma}{6}}
dx\right]^{\frac{2}{\gamma}}\leq C
 N^{\frac{2}{3}}N^{\frac{1}{3}-\frac{2}{\gamma}}=N^{1-\frac{2}{\gamma}}.
\end{align}
Combining (\ref{2.016}), (\ref{2.017}) with (\ref{2.019}), we have
\begin{align*}
\|K\|_{L_{x}^{\frac{\gamma}{2}}L_{t}^{\infty}}\leq C
 N^{\frac{\gamma-2}{\gamma}}.
\end{align*}

We have completed the proof of  Lemma 2.6.

\noindent {\bf Remark3:} We follow the idea of Proposition 2.5 of \cite{Porn},
 we derive Lemma 2.6.

\begin{Lemma}\label{lem2.7}(Maximal function estimate.)
 Let $f\in L_{x}^{2}(\R)$  and  $\gamma\geq 4$ and  $supp(|\mathscr{F}_{x}{f}|)
 \subseteq [N,4N]$,
  where $N\in 2^{Z}$, we have
\begin{align}\label{2.020}
\|U(t)f(x)\|_{L_{x}^{\gamma}L_{t}^{\infty}}\leq C N^{\frac{\gamma-2}{2\gamma}}
\|f\|_{L_{x}^{2}},
\end{align}
where
\begin{eqnarray}
U(t)f(x)=\int_{\SR}e^{ix\xi+it\xi^{3}}\mathscr{F}_{x}f(\xi)d\xi.\label{2.021}
\end{eqnarray}
In particular, for $b>\frac{1}{2},$  $\gamma\geq 4$ and  $supp(|\mathscr{F}_{x}{u}|)
\subseteq [N,4N]$,
where $N\in 2^{Z}$,      we have
\begin{eqnarray}
&&\|u\|_{L_{x}^{\gamma}L_{t}^{\infty}}\leq C N^{\frac{\gamma-2}{2\gamma}}
\|u\|_{X_{0,b}}.\label{2.022}
\end{eqnarray}

\end{Lemma}
\noindent{\bf Proof.}The maximal function inequality
\begin{eqnarray}
\|U(t)f(x)\|_{L_{x}^{4}L_{t}^{\infty}}\leq C N^{\frac{1}{4}}
\|f\|_{L_{x}^{2}}\label{2.023}
\end{eqnarray}
can be seen in \cite{KPVCPAM}. We firstly consider the case $\gamma \geq 7.$
  Let $Tf=U(t)f$. Obviously, $T$  is an operator which  maps   $L_{x}^{2}$ into $L_{x}^{\gamma}L_{t}^{\infty}$.
Thus, we have $T^{*}F=\int_{\SR}e^{-t\partial_{x}^{3}}F dt$. With the aid of  $TT^{*}$ argument, we obtain
 that  to obtain (\ref{2.020}), it suffices  to prove
 \begin{align}\label{2.025}
 \left\|\int_{\SR}U(t-s)Fds\right\|_{L_{x}^{\gamma}L_{t}^{\infty}}\leq C
  N^{\frac{\gamma-2}{\gamma}}\|F\|_{L_{x}^{\frac{\gamma}{\gamma-1}}L_{t}^{1}}.
 \end{align}
Here  $F\in L_{x}^{2}L_{t}^{1}(\R\times\R)$ has  the same frequency support as $u$.
 Thus, it suffices to prove (\ref{2.025}).
 Obviously,
 \begin{align*}
 \mathscr{F}_{x}^{-1}\left(e^{i(t-s)\xi^{3}}\mathscr{F}_{x}{F}(\xi,s)\right)
 &=C\int_{\SR}e^{i(t-s)\xi^{3}+ix\xi}\mathscr{F}_{x}{F}(\xi,s)d\xi\\
 &=\mathscr{F}_{x}^{-1}\left(e^{i(t-s)\xi^{3}}\psi_{1}\left(\frac{\xi}{4N}\right)\right)*F(x,s).
 \end{align*}
By using a direct computation, we have that  the term on the left hand side of (\ref{2.025}) equals to
 \begin{align*}
 \int_{\SR}\mathscr{F}_{x}^{-1}\left(e^{i(t-s)\xi^{3}}\psi_{1}\left(\frac{\xi}{4N}\right)\right)*F(x,s)ds
 &=\mathscr{F}_{x}^{-1}\left(e^{it\xi^{3}}\psi_{1}\left(\frac{\xi}{4N}\right)\right)*F(x,t)=CK*F.
 \end{align*}
 Here $*$ denotes convolution with respect  to  the space-time variable    and
\begin{align}\label{2.026}
 K(x,t)=\int_{\SR}e^{it\xi^{3}+ix\xi}\psi_{1}\left(\frac{\xi}{4N}\right)d\xi.
 \end{align}
 By using Young's inequality and Lemma 2.6, we obtain
 \begin{align*}
\|K*F\|_{L_{x}^{\gamma}L_{t}^{\infty}}\leq \|K\|_{L_{x}^{\frac{\gamma}{2}}L_{t}^{\infty}}
\|F\|_{L_{x}^{\frac{\gamma}{\gamma-1}}L_{t}^{1}}\leq C
  N^{\frac{\gamma-2}{\gamma}}\|F\|_{L_{x}^{\frac{\gamma}{\gamma-1}}L_{t}^{1}}.
 \end{align*}
 Thus, when $\gamma\geq7,$ we have
 \begin{align}\label{2.027}
\|U(t)f(x)\|_{L_{x}^{\gamma}L_{t}^{\infty}}\leq C N^{\frac{\gamma-2}{2\gamma}}
\|f\|_{L_{x}^{2}}.
\end{align}
Interpolating (\ref{2.023}) with (\ref{2.027}) yields
\begin{align}\label{2.028}
\|U(t)f(x)\|_{L_{x}^{\gamma}L_{t}^{\infty}}\leq C N^{\frac{\gamma-2}{2\gamma}}
\|f\|_{L_{x}^{2}},4\leq \gamma \leq 7.
\end{align}
Obviously, changing variable $\tau=\lambda+\xi^{3}$, we have
 \begin{eqnarray}
&&u(x,t)=\frac{1}{2\pi}\int_{\SR^{2}}e^{ix\xi+it\tau}\mathscr{F}u(\xi,\tau)
d\xi d\tau\nonumber\\
&&=\frac{1}{2\pi}\int_{\SR^{2}}e^{ix\xi+it(\lambda+\xi^{3})}
\mathscr{F}u(\xi,\lambda+\xi^{3})d\xi d\lambda\nonumber\\
&&=\frac{1}{2\pi}\int_{\SR}e^{it\lambda}\left(\int_{\SR}e^{ix\xi+it\xi^{3}}
\mathscr{F}u(\xi,\lambda+\xi^{3}) d\xi\right)d\lambda\label{2.029}.
 \end{eqnarray}
By using  (\ref{2.028}),  (\ref{2.029}),  Minkowski's inequality
 and taking $b>\frac{1}{2}$ and variable substitution $\lambda+\xi^{3}=\tau,$ we have
\begin{eqnarray}
&&\|u\|_{L_{x}^{\gamma}L_{t}^{\infty}}\leq C\int_{\SR}\left\|\left(\int_{\SR}e^{ix\xi+it\Phi(\xi)}
\mathscr{F}u(\xi,\lambda+\xi^{3})d\xi\right)
\right\|_{L_{x}^{\gamma}L_{t}^{\infty}} d\lambda\nonumber\\&&\leq CN^{\frac{\gamma-2}{2\gamma}}
\int_{\SR}\left\|\mathscr{F}u(\xi,\lambda+\xi^{3})\right\|_{L^{2}}
d\lambda\nonumber\\&&\leq CN^{\frac{\gamma-2}{2\gamma}}
\left[\int_{\SR}(1+|\lambda|)^{2b}\left\|\mathscr{F}u(\xi,\lambda+\xi^{3})\right\|_{L^{2}}^{2}
d\lambda\right]^{\frac{1}{2}}
\left[\int_{\SR}(1+|\lambda|)^{-2b}d\lambda\right]^{\frac{1}{2}}\nonumber\\
&&\leq CN^{\frac{\gamma-2}{2\gamma}}\left[\int_{\SR}(1+|\tau-\xi^{3}|)^{2b}\left\|\mathscr{F}u(\xi,\tau)
\right\|_{L^{2}}^{2}d\tau\right]^{\frac{1}{2}}
=CN^{\frac{\gamma-2}{2\gamma}}\|u\|_{X_{0,b}}.\label{2.030}
\end{eqnarray}

We have completed the proof of  Lemma 2.7.

\begin{Lemma}\label{lem2.8}
 Let $f\in L_{x}^{2}(\R)$  and  $\gamma\geq 4$ and  $supp(|\mathscr{F}_{x}{u}|)\subseteq [N,4N]$,
  where $N\in 2^{Z}$.  For $b>\frac{1}{2},$    we have
\begin{align}\label{2.031}
\|u\|_{L_{xt}^{\infty}}\leq C N^{\frac{1}{2}}\|u\|_{X_{0,b}}.
\end{align}
\end{Lemma}
\noindent{\bf Proof.} Obviously, since $supp(|\mathscr{F}_{x}{u}|)\subseteq [N,4N]$, by using
 the Cauchy-Schwarz inequality and the Minkowski's inequality  as well as the Plancherel identity, we have
\begin{eqnarray}
&&\left|u(x)\right|=\left|\int_{N}^{4N}e^{ix\xi}\mathscr{F}_{x}u(\xi,t)d\xi\right|
\leq \int_{N}^{4N}\left|\mathscr{F}_{x}u(\xi,t)\right|d\xi\leq CN^{\frac{1}{2}}\left[\int_{N}^{4N}
\left|\mathscr{F}_{x}u(\xi,t)\right|^{2}d\xi\right]^{\frac{1}{2}}\nonumber\\
&&=CN^{\frac{1}{2}}\|\mathscr{F}_{x}u(\xi,t)\|_{L_{\xi}^{2}}=CN^{\frac{1}{2}}\|u\|_{L_{x}^{2}}\label{2.032}.
\end{eqnarray}
Thus, from (\ref{2.032}) and Lemma 2.3 of \cite{KPV-Duke}, we have
\begin{eqnarray}
&&\left\|u\right\|_{L_{xt}^{\infty}}\leq CN^{\frac{1}{2}}\|u\|_{L_{t}^{\infty}L_{x}^{2}}\leq
 CN^{\frac{1}{2}}\|u\|_{L_{x}^{2}L_{t}^{\infty}}\leq CN^{\frac{1}{2}}\|u\|_{X_{0,b}}\label{2.033}.
\end{eqnarray}

This completes the proof of Lemma 2.8.

\begin{Lemma}\label{lem2.9}Suppose $f\in L^{p}(\R),1<p<\infty.$ Then, there exist two
 positive constants $C_{1},C_{2}$ such that
\begin{eqnarray}
C_{1}\|f\|_{L^{p}}\leq \left\|\left[\sum\limits_{N}|P_{N}f|^{2}\right]^{\frac{1}{2}}\right\|_{L^{p}}
=\left\|\left\|P_{N}f\right\|_{\ell_{N}^{2}}\right\|_{L^{p}}\leq C_{2}\|f\|_{L^{p}}\label{2.034}
\end{eqnarray}
and
\begin{eqnarray}
\left\|\left\|P_{N}g\right\|_{L_{t}^{p}L_{x}^{q}}\right\|_{\ell_{N}^{r}}\leq C\|g\|_{L_{t}^{p}L_{x}^{q}}.\label{2.035}
\end{eqnarray}
Here $r\geq{\rm max}\left\{p,q\right\}$ and $g\in L_{t}^{p}L_{x}^{q}.$
\end{Lemma}

\noindent {\bf Proof.} (\ref{2.034})  can be seen in Theorem 1.5
 of page 104 in
 \cite{Stein}. By using the Minkowski's inequality and
  $\ell^{2}\hookrightarrow \ell^{r}(r\geq2)$ as well as (\ref{2.034}),
   we have
\begin{eqnarray*}
\left\|P_{N}g\right\|_{\ell_{N}^{r}L_{t}^{p}L_{x}^{q}}\leq
 C\|P_{N}g\|_{L_{t}^{p}L_{x}^{q}\ell_{N}^{r}}\leq C
\|P_{N}g\|_{L_{t}^{p}L_{x}^{q}\ell_{N}^{2}}\leq C\|g\|_{L_{t}^{p}L_{x}^{q}}.
\end{eqnarray*}

This completes the proof of Lemma 2.9.

\begin{Lemma}\label{lem2.10}

Let $g_{n}(\omega)(1\leq n\leq \infty)$ be a sequence of real, 0-mean, independent
random variables with associated sequence of distributions  Assume
that $\mu_{n}(1\leq n\leq \infty)$ satisfy the property (\ref{1.05}).
Then there exists $\alpha> 0$ such that for every $\lambda>0$ every
 sequence $(c_{n})(1\leq n\leq \infty)\in l^{2}(\Z)$
 of real numbers,
\begin{eqnarray*}
\mathbb{P}\left(\left\{\omega\in \Omega:
\left|\sum\limits_{k=1}^{\infty}c_{k}g_{k}(\omega)\right|>\lambda\right\}\right)
\leq 2e^{-\frac{\alpha\lambda^{2}}{\sum\limits_{k=1}^{\infty}c_{k}^{2}}}.
\end{eqnarray*}
As a consequence there exists $C>0$ such that for every $p\geq2,$ every
$\{c_k\}\in l^{2}(\Z)$
\begin{eqnarray}
\left\|\sum_{k\in\z}g_k(\omega)c_k\right\|_{L_{\omega}^p(\Omega)}
\leq C\sqrt{p}\left\|c_k\right\|_{l^{2}(\z)}.\label{2.036}
\end{eqnarray}
\end{Lemma}

For the proof of Lemma 2.10, we refer the readers to Lemma 3.1 of \cite{BT2008-L}.

\begin{Lemma}\label{Lemma2.11}
Let $f\in H^{s}(\R)$ and  $f^{\omega}=\sum\limits_{n\in Z}g_{n}(\omega) \psi(D-n)f$.
 Then,
 there exist $C_{1}>0, C>0$ such that
\begin{eqnarray}
{\rm P} \left(\Omega_{1}\right)\leq C_{1}e^{-C\lambda^{2}\|\phi\|_{H^{s}}^{-2}}\label{2.037}
\end{eqnarray}
for all $\lambda>0$.
Here, $\Omega_{1}=\left\{\omega:\|f^{\omega}\|_{H^{s}}>\lambda\right\}.$
\end{Lemma}

For the proof of Lemma 2.11, we refer the readers to Lemma 2.2 of \cite{BOP}.

\begin{Lemma}\label{lem2.12}

 Let $f\in H^{s}(\R)$ and we denote by $f^{\omega}$
 the randomization of $f$ as defined in (\ref{1.020}) and $b>\frac{1}{2}.$  Then,
 there exist $C_{1}>0, C>0$ such that
\begin{eqnarray}
\mathbb{P}\left(\Omega_{2}\right)\leq C_{1}e^{-C\lambda^{2}\|f\|_{H^{s}}^{-2}}\label{2.038}
\end{eqnarray}
 for all $\lambda>0.$
Here
$
\Omega_{2}=\left\{\omega\in \Omega:
\left\|\phi(t)U(t)f^{\omega}\right\|_{X_{s,b}}>\lambda\right\}.
$

\end{Lemma}
\noindent {\bf Proof.} From Lemma 2.1, we have
\begin{eqnarray}
\left\|\phi(t)S(t)f^{\omega}\right\|_{X_{s,b}}\leq C\|f^{\omega}\|_{H^{s}}\label{2.039}.
\end{eqnarray}
Thus, from (\ref{2.037}) and Lemma 2.10, we have
\begin{eqnarray*}
\mathbb{P}\left(\Omega_{2}\right)\leq \mathbb{P}\left(\omega\in \Omega:
C\left\|f^{\omega}\right\|_{H^{s}}>\lambda\right)\leq C_{1}e^{-C\lambda^{2}\|f\|_{H^{s}}^{-2}}.
\end{eqnarray*}

This completes the proof of Lemma 2.12.$\hfill\Box$

\begin{Lemma}\label{lem2.13}(Improved local-in-time Strichartz estimate.)
 Let $f\in L^{2}(\R)$ and we denote by $f^{\omega}$
 the randomization of $f$ as defined in (\ref{1.020}).
   Then,   for $0<T\leq1,$   $2\leq q,r<\infty,$
    there exist $C>0, C_{1}>0$ such that
\begin{eqnarray}
\mathbb{P}\left(\Omega_{3}\right)\leq C_{1}
e^{-C\left(\frac{\lambda}{T^{\frac{1}{q}}\|f\|_{L^{2}}}\right)^{2}}\label{2.040}
\end{eqnarray}
for all $ \lambda>0.$
Here
$
\Omega_{3}=\left\{\omega\in \Omega:
\left\|U(t)f^{\omega}\right\|_{L_{t}^{q}L_{x}^{r}([0,T]\times \SR)}>\lambda\right\}.
$

\end{Lemma}

Lemma 2.13 can be seen in Proposition 1.3 of \cite{BOP}.

\begin{Lemma}\label{lem2.14}
 Let $f\in H^{8\epsilon}(\R)$
 and we denote by $f^{\omega}$ the randomization of $f$ as defined in
 (\ref{1.020}).  Then,
     there exist  $C>0$ and $C_{1}>0$ such that
\begin{eqnarray}
\mathbb{P}\left(\Omega_{4}\right)\leq C_{1} e^{-C\left(\frac{\lambda}
{\|f\|_{H^{8\epsilon}}}\right)^{2}} \label{2.041}
\end{eqnarray}
$\forall \lambda>0.$
Here
\begin{eqnarray}
\Omega_{4}=\left\{\omega\in \Omega: \left\|U(t)f^{\omega}\right\|_{L_{xt}^{\infty}
(\SR\times [0,T])}>\lambda\right\}.\label{2.042}
\end{eqnarray}
Here, $C,C_{1}$ is independent of $x,t$.
\end{Lemma}
\noindent {\bf Proof.} When $p\geq\frac{1}{\epsilon},$ by using
  the Sobolev embeddings theorem,
 Lemma  2.9 and the Minkowski's inequality,  we have
\begin{eqnarray}
&&\hspace{-1.2cm}\left\|\left\|U(t)f^{\omega}\right\|_{L_{xt}^{\infty}(\SR\times [0,T])}\right\|_{L_{\omega}^{p}(\Omega)}
=\left\|\left\|\sum\limits_{k\in \z}g_{k}(\omega)U(t)\psi(D-k)f\right\|_{L_{xt}^{\infty}(\SR\times [0,T])}\right\|_{L_{\omega}^{p}(\Omega)}\nonumber\\&&
\hspace{-1.2cm}=\left\|\left\|\sum\limits_{k\in \z}g_{k}(\omega)J_{x}^{2\epsilon}
J_{t}^{2\epsilon}U(t)\psi(D-k)f\right\|_{L_{xt}^{\frac{1}{\epsilon}}{(\SR\times [0,T])}}
\right\|_{L_{\omega}^{p}(\Omega)}
\nonumber\\
&&\hspace{-1.2cm}\leq C\left\|\left\|\sum\limits_{k\in \z}g_{k}(\omega)J_{x}^{2\epsilon}J_{t}^{2\epsilon}
U(t)\psi(D-k)f\right\|_{L_{\omega}^{p}(\Omega)}
\right\|_{L_{xt}^{\frac{1}{\epsilon}}{(\SR\times [0,T])}}\nonumber\\
&&\hspace{-1.2cm}\leq C\sqrt{p}\left\|\left[\sum\limits_{k\in \z}|J_{x}^{2\epsilon}
J_{t}^{2\epsilon}U(t)\psi(D-k)f|^{2}\right]^{\frac{1}{2}}
\right\|_{L_{xt}^{\frac{1}{\epsilon}}{(\SR\times [0,T])}}\nonumber\\&&\leq C
\sqrt{p}\left\|J_{x}^{2\epsilon}J_{t}^{2\epsilon}U(t)\psi(D-k)f
\right\|_{\ell_{k}^{2}L_{xt}^{\frac{1}{\epsilon}}{(\SR\times [0,T])}}\nonumber\\&&
\hspace{-1.2cm}\leq C\sqrt{p}\left\|\left\|J_{x}^{2\epsilon}J_{t}^{2\epsilon}U(t)
\psi(D-k)f\right\|_{L_{x}^{\frac{1}{\epsilon}}}\right\|_{\ell_{k}^{2}L_{t}^{\frac{1}{\epsilon}}([0,T])}
\nonumber\\&&\hspace{-1.2cm}\leq C\sqrt{p}\left\|\left\|J_{x}^{2\epsilon}J_{t}^{2\epsilon}U(t)
\psi(D-k)f\right\|_{L_{x}^{2}}\right\|_{\ell_{k}^{2}L_{t}^{\frac{1}{\epsilon}}([0,T])}\nonumber\\
&&\hspace{-1.2cm}\leq C\sqrt{p}\left\|\left\|\psi(D-k)f\right\|_{H^{8\epsilon}}
\right\|_{\ell_{k}^{2}L_{t}^{\frac{1}{\epsilon}}([0,T])}\leq CT^{\epsilon}\sqrt{p}\left\|\psi(D-k)f\right\|_{\ell_{k}^{2}H^{8\epsilon}}\nonumber\\&&\hspace{-1.2cm}
\leq C\sqrt{p}\|f\|_{H^{8\epsilon}}.\label{2.043}
\end{eqnarray}
Thus, by Chebyshev inequality, from (\ref{2.043}), we have
\begin{eqnarray}
&&\mathbb{P}\left(\Omega_{4}\right)\leq \int_{\Omega_{4}}\left[
\frac{\left\|U(t)f\right\|_{L_{xt}^{\infty}}}{\lambda}\right]^{p}d\mathbb{P}(\omega)
\leq \left(\frac{C\sqrt{p}\|f\|_{H^{8\epsilon}}}{\lambda}\right)^{p}\label{2.044}.
\end{eqnarray}
Take
\begin{eqnarray}
p=\left(\frac{\lambda}{Ce\|f\|_{H^{8\epsilon}}}\right)^{2}\label{2.045}.
\end{eqnarray}
If $p\geq2$,
 we have
\begin{eqnarray}
\mathbb{P}(\Omega_{4})\leq e^{-p}=
e^{-\left(\frac{\lambda}{Ce\|f\|_{H^{8\epsilon}}}\right)^{2}}.\label{2.046}
\end{eqnarray}
If $p\leq2,$
  we have
\begin{eqnarray}
&&\mathbb{P}\left(\Omega_{4}\right)\leq e^{2}e^{-2}=C_{1}e^{-\left(\frac{\lambda}
{Ce\|f\|_{H^{8\epsilon}}}\right)^{2}}\label{2.047}.
\end{eqnarray}
Here $C_{1}=e^{2}.$

This completes the proof of Lemma 2.14.

\begin{Lemma}\label{lem2.15}
  Let $f\in H^{6\epsilon}(\R)$
 and we denote by $f^{\omega}$ the randomization of $f$ as defined in
 (\ref{1.020}).  Then,
     there exist  $C>0$ and $C_{1}>0$ such that
\begin{eqnarray*}
\mathbb{P}\left(\Omega_{5}\right)\leq C_{1} e^{-C\left(\frac{\lambda}
{\|f\|_{H^{6\epsilon}}}\right)^{2}}\leq\epsilon
\end{eqnarray*}
 $\forall \lambda>0.$
Here
$
\Omega_{5}=\left\{\omega\in \Omega: \left\|U(t)f^{\omega}\right\|_{L_{t}^{\infty}L_{x}^{\frac{1}{\epsilon}}(\SR\times [0,T])}>\lambda\right\}.
$
Here, $C,C_{1}$ is independent of $x,t$.
\end{Lemma}

Lemma 2.15 can be proved similarly to Lemma 2.14.

\bigskip
\bigskip

\noindent{\large\bf 3. Multilinear estimates}

\setcounter{equation}{0}

 \setcounter{Theorem}{0}

\setcounter{Lemma}{0}

 \setcounter{section}{3}
 In this section, we present some crucial multilinear estimates which play an
  important role in establishing Theorems 1.1-1.7.

\begin{Lemma}\label{Lemma3.1}
Let $s\geq\frac{1}{2}-\frac{2}{k}+88\epsilon,k\geq5$ and $b=\frac{1}{2}+\frac{\epsilon}{24},b_{1}=\frac{1}{2}-\frac{\epsilon}{12}$.
Then, we have
\begin{eqnarray}
\left\|\partial_{x}(\prod\limits_{j=1}^{k+1}u_{j})\right\|_{X_{s,b_{1}}}\leq C\prod_{j=1}^{k+1}\|u_{j}\|_{X_{s,b}}.\label{3.01}
\end{eqnarray}
\end{Lemma}
\noindent {\bf Proof.}To prove (\ref{3.01}), by duality,  it suffices to prove
\begin{eqnarray}
\left|\int_{\SR^{2}}J^{s}\partial_{x}\left(\prod_{j=1}^{k+1}u_{j}\right)\bar{h}dxdt\right|\leq C\|h\|_{X_{0,b_{1}}}\prod_{j=1}^{k+1}\|u_{j}\|_{X_{s,b}}\label{3.02}.
\end{eqnarray}
We define
\begin{eqnarray*}
f(\xi,\tau)=\langle\sigma\rangle^{b_{1}} \mathscr{F}h(\xi,\tau),
f_{j}(\xi_{j},\tau_{j})=\langle\xi_{j}\rangle^{s}\langle\sigma_{j}\rangle^{b}
\mathscr{F}u_{j}(\xi_{j},\tau_{j})(1\leq j\leq k+1).
\end{eqnarray*}
By using the Plancherel identity, to prove (\ref{3.02}), it suffices to prove
\begin{eqnarray}
\int_{\xi=\sum\limits_{j=1}^{k+1}\xi_{j}}\int_{\tau=\sum\limits_{j=1}^{k+1}\tau_{j}}
\frac{|\xi|\langle\xi\rangle^{s}f\prod\limits_{j=1}^{k+1}f_{j}}
{\langle\sigma\rangle^{b_{1}}\prod\limits_{j=1}^{k+1}
\langle\xi_{j}\rangle^{s}\langle\sigma_{j}\rangle^{b}}d\delta
&&\leq C\|f\|_{L^{2}}\prod\limits_{j=1}^{k+1}\|f_{j}\|_{L^{2}},\label{3.03}
\end{eqnarray}
where $d\delta=d\xi_{1}d\xi_{2}\cdot \cdot \cdot d\xi_{k}d\xi d\tau_{1}d\tau_{2}\cdot \cdot \cdot d\tau_{k}d\tau.$

\noindent We define
\begin{eqnarray*}
&&K(\xi_{1},\xi_{2},\cdot\cdot\cdot,\xi_{k},\xi,\tau_{1},\tau_{2},\cdot\cdot\cdot,\tau_{k},\tau)
=\frac{|\xi|\langle\xi\rangle^{s}}
{\langle\sigma\rangle^{b_{1}}\prod\limits_{j=1}^{k+1}
\langle\xi_{j}\rangle^{s}\langle\sigma_{j}\rangle^{b}},\nonumber\\&&
\mathscr{F}F=\frac{f}{\langle\sigma\rangle^{b_{1}}},
\mathscr{F}F_{j}=\frac{f_{j}}{\langle\sigma_{j}\rangle^{b}}
(1\leq j\leq k+1),\\
&&I_{1}=\int_{\xi=\sum\limits_{j=1}^{k+1}\xi_{j}}
\int_{\tau=\sum\limits_{j=1}^{k+1}\tau_{j}}K(\xi_{1},\xi_{2},
\cdot\cdot\cdot,\xi_{k},\xi,\tau_{1},\tau_{2},\cdot\cdot\cdot,\tau_{k},\tau)
f\prod\limits_{j=1}^{k+1}f_{j}d\delta.
\end{eqnarray*}
Without loss of generality, we can assume that $|\xi_{1}|\geq
 |\xi_{2}|\geq\cdot\cdot\cdot\geq |\xi_{k+1}|.$

\noindent Obviously,
\begin{eqnarray*}
\Omega:=\left\{(\xi_{1},\xi_{2},\cdot\cdot\cdot,\xi_{k},\xi,\tau_{1},\tau_{2},
\cdot\cdot\cdot,\tau_{k},\tau)\in \R^{2(k+1)}:|\xi_{1}|\geq |\xi_{2}|\geq
\cdot\cdot\cdot\geq |\xi_{k+1}|\right\}\subset \bigcup\limits_{j=0}^{k+1}\Omega_{j}.
\end{eqnarray*}
Here,
\begin{eqnarray*}
&&\Omega_{0}=\left\{(\xi_{1},\xi_{2},\cdot\cdot\cdot,\xi_{k},\xi,\tau_{1},\tau_{2},
\cdot\cdot\cdot,\tau_{k},\tau)\in \Omega,|\xi_{1}|\leq80k\right\},\\
&&\Omega_{1}=\left\{(\xi_{1},\xi_{2},\cdot\cdot\cdot,\xi_{k},\xi,\tau_{1},\tau_{2},
\cdot\cdot\cdot,\tau_{k},\tau)\in \Omega,|\xi_{1}|\geq80k,|\xi_{1}|\geq80k |\xi_{2}|\right\},\\
&&\Omega_{2}=\left\{(\xi_{1},\xi_{2},\cdot\cdot\cdot,\xi_{k},\xi,\tau_{1},\tau_{2},
\cdot\cdot\cdot,\tau_{k},\tau)\in \Omega,|\xi_{1}|\geq80k, |\xi_{1}|\sim |\xi_{2}|\geq 80k|\xi_{3}|\right\},\\
&&\Omega_{3}=\left\{(\xi_{1},\xi_{2},\cdot\cdot\cdot,\xi_{k},\xi,\tau_{1},\tau_{2},
\cdot\cdot\cdot,\tau_{k},\tau)\in \Omega,|\xi_{1}|\geq80k, |\xi_{1}|\sim
 |\xi_{3}|\geq 80k|\xi_{4}|\right\},\\
&&\Omega_{4}=\left\{(\xi_{1},\xi_{2},\cdot\cdot\cdot,\xi_{k},\xi,\tau_{1},\tau_{2},
\cdot\cdot\cdot,\tau_{k},\tau)\in \Omega,|\xi_{1}|\geq80k, |\xi_{1}|\sim|\xi_{4}|\geq 80k|\xi_{5}|\right\},\\
&&\Omega_{5}=\left\{(\xi_{1},\xi_{2},\cdot\cdot\cdot,\xi_{k},\xi,\tau_{1},\tau_{2},
\cdot\cdot\cdot,\tau_{k},\tau)\in \Omega,|\xi_{1}|\geq80k, |\xi_{1}|\sim  |\xi_{5}|\geq 80k|\xi_{6}|\right\},\\
&&\Omega_{6}=\left\{(\xi_{1},\xi_{2},\cdot\cdot\cdot,\xi_{k},\xi,\tau_{1},\tau_{2},
\cdot\cdot\cdot,\tau_{k},\tau)\in \Omega,|\xi_{1}|\geq80k, |\xi_{1}|\sim |\xi_{l-1}|\geq 80k|\xi_{l}|(7\leq l\leq k+1)\right\},\\
&&\Omega_{7}=\left\{(\xi_{1},\xi_{2},\cdot\cdot\cdot,\xi_{k},\xi,\tau_{1},\tau_{2},
\cdot\cdot\cdot,\tau_{k},\tau)\in \Omega,|\xi_{1}|\geq80k, |\xi_{1}|\sim  |\xi_{k+1}|\right\}.
\end{eqnarray*}
(1)Case $(\xi_{1},\xi_{2},\cdot\cdot\cdot,\xi_{k},\xi,\tau_{1},\tau_{2},\cdot\cdot\cdot,
\tau_{k},\tau)\in \Omega_{0}$,
 by using the Plancherel identity,  the H\"older inequality and  (\ref{2.05}), we have
\begin{eqnarray}
&&I_{1}\leq C\int_{\xi=\sum\limits_{j=1}^{k+1}\xi_{j}}\int_{\tau=\sum\limits_{j=1}^{k+1}\tau_{j}}
\frac{f\prod\limits_{j=1}^{k+1}f_{j}}
{\langle\sigma\rangle^{b_{1}}\prod\limits_{j=1}^{k+1}
\langle\xi_{j}\rangle^{s}\langle\sigma_{j}\rangle^{b}}d\delta\nonumber\\
&&\leq C\|F\|_{L_{xt}^{2}}\prod\limits_{j=1}^{k+1}\|F_{j}\|_{L_{xt}^{2(k+1)}}
\leq C\|f\|_{L^{2}}\prod\limits_{j=1}^{k+1}\|F_{j}\|_{X_{\frac{1}{2}-\frac{2}{k+1},b}}
\leq C\|f\|_{L^{2}}\prod\limits_{j=1}^{k+1}\|f_{j}\|_{L^{2}}\label{3.04}.
\end{eqnarray}
(2)Case  $(\xi_{1},\xi_{2},\cdot\cdot\cdot,\xi_{k},\xi,\tau_{1},\tau_{2},\cdot\cdot\cdot,\tau_{k},\tau)\in \Omega_{1}$,
we have
\begin{eqnarray}
&&K(\xi_{1},\xi_{2},\cdot\cdot\cdot,\xi_{k},\xi,\tau_{1},\tau_{2},\cdot\cdot\cdot,\tau_{k},\tau)\leq C
\frac{|\xi_{1}^{2}-\xi_{2}^{2}|^{\frac{1}{2}}}
{\langle\sigma\rangle^{\frac{1}{2}-\frac{\epsilon}{12}}\prod\limits_{j=2}^{k+1}
\langle\xi_{j}\rangle^{s}\prod\limits_{j=1}^{k+1}\langle\sigma_{j}\rangle^{b}}
.\label{3.05}
\end{eqnarray}
By using  (\ref{3.05}), the Plancherel identity,  the H\"older inequality,
Lemma 2.2, (\ref{2.05}) and (\ref{2.08}),
we have
\begin{eqnarray*}
&&I_{1}\leq C\int_{\xi=\sum\limits_{j=1}^{k+1}\xi_{j}}\int_{\tau=\sum\limits_{j=1}^{k+1}\tau_{j}}
\frac{f(\prod\limits_{j=1}^{k+1}f_{j})|\xi_{1}^{2}-\xi_{2}^{2}|^{\frac{1}{2}}\prod\limits_{j=2}^{k+1}
\langle\xi_{j}\rangle^{-s}}{\langle\sigma\rangle^{\frac{1}{2}-\frac{\epsilon}{12}}
\prod\limits_{j=1}^{k+1}\langle\sigma_{j}\rangle^{b}}d\delta\nonumber\\&&\leq C
\left\|I^{1/2}(F_{1},F_{2})\right\|_{L_{xt}^{2}}
\left(\prod_{j=3}^{k+1}\|J^{-\frac{ks}{k-1}}F_{j}\|_{L_{xt}^{\frac{8(k-1)}{3-\epsilon}}}\right)
\|F\|_{L_{xt}^{\frac{8}{1+\epsilon}}}\nonumber\\
&&\leq C\|F\|_{X_{0,b_{1}}}\prod\limits_{j=1}^{k+1}\|F_{j}\|_{X_{0,b}}\leq
C\left(\prod_{j=1}^{k+1}\|f_{j}\|_{L_{xt}^{2}}\right)\|f\|_{L_{xt}^{2}}.
\end{eqnarray*}
(3)Case $(\xi_{1},\xi_{2},\cdot\cdot\cdot,\xi_{k},\xi,\tau_{1},\tau_{2},
\cdot\cdot\cdot,\tau_{k},\tau)\in \Omega_{2}$,
 since $s\geq \frac{1}{2}-\frac{2}{k}+88\epsilon,$ we  have
\begin{eqnarray}
&&K(\xi_{1},\xi_{2},\cdot\cdot\cdot,\xi_{k},\xi,\tau_{1},\tau_{2},\cdot\cdot\cdot,\tau_{k},\tau)\leq
C\frac{|\xi_{2}^{2}-\xi_{3}^{2}|^{\frac{1}{2}}}
{\langle\sigma\rangle^{b_{1}}\prod\limits_{j=2}^{k+1}
\langle\xi_{j}\rangle^{s}\prod\limits_{j=1}^{k+1}\langle\sigma_{j}\rangle^{b}}
\nonumber\\&&\leq C\frac{|\xi_{2}^{2}-\xi_{3}^{2}|^{\frac{1}{2}}
|\xi_{1}|^{\frac{1}{6}}\prod\limits_{j=4}^{k+1}
\langle\xi_{j}\rangle^{-\frac{3k-11}{6(k-2)}-2k\epsilon}}{\langle\sigma\rangle^{b_{1}}
\prod\limits_{j=1}^{k+1}\langle\sigma_{j}\rangle^{b}}.\label{3.06}
\end{eqnarray}
By using (\ref{3.06}), the Plancherel identity, the H\"older inequality,
  Lemma 2.2, (\ref{2.05})-(\ref{2.06}),
(\ref{2.08}),  we have
\begin{eqnarray*}
&&I_{1}\leq C\int_{\xi=\sum\limits_{j=1}^{k+1}\xi_{j}}\int_{\tau=\sum\limits_{j=1}^{k+1}\tau_{j}}
\frac{f(\prod\limits_{j=1}^{k+1}f_{j})|\xi_{2}^{2}-\xi_{3}^{2}|^{\frac{1}{2}}|\xi_{1}|^{\frac{1}{6}}
\prod\limits_{j=4}^{k+1}
\langle\xi_{j}\rangle^{-\frac{3k-11}{6(k-2)}-2\epsilon}}{\langle\sigma\rangle^{b_{1}}
\prod\limits_{j=1}^{k+1}\langle\sigma_{j}\rangle^{b}}d\delta\nonumber\\&&\leq C
\left\|I^{1/2}(F_{2},F_{3})\right\|_{L_{xt}^{2}}
\left(\prod_{j=4}^{k+1}\|J^{-\frac{3k-11}{6(k-2)}-2\epsilon}F_{j}\|_{L_{xt}^{\frac{24(k-2)}{5-3\epsilon}}}\right)
\|D_{x}^{\frac{1}{6}}F_{1}\|_{L_{xt}^{6}}\|F\|_{L_{xt}^{\frac{8}{1+\epsilon}}}\nonumber\\
&&\leq C\left(\prod_{j=1}^{k+1}\|F_{j}\|_{X_{0,b}}\right)
\|F\|_{X_{0,b_{1}}}\leq C\left(\prod_{j=1}^{k+1}\|f_{j}\|_{L_{xt}^{2}}\right)\|f\|_{L_{xt}^{2}}.
\end{eqnarray*}
(4)Case $(\xi_{1},\xi_{2},\cdot\cdot\cdot,\xi_{k},\xi,\tau_{1},
\tau_{2},\cdot\cdot\cdot,\tau_{k},\tau)\in \Omega_{3}$, we have
\begin{eqnarray}
&&K(\xi_{1},\xi_{2},\cdot\cdot\cdot,\xi_{k},\xi,\tau_{1},\tau_{2},\cdot\cdot\cdot,\tau_{k},\tau)\leq C
\frac{|\xi_{3}^{2}-\xi_{4}^{2}|^{\frac{1}{2}}}
{\langle\sigma\rangle^{b_{1}}\prod\limits_{j=2}^{k+1}
\langle\xi_{j}\rangle^{s}\prod\limits_{j=1}^{k+1}\langle\sigma_{j}\rangle^{b}}
.\label{3.07}
\end{eqnarray}
 By using (\ref{3.07}),  the Plancherel identity, the  H\"older inequality and
  (\ref{2.05}), (\ref{2.07})-(\ref{2.08}) as well as Lemma 2.2,
we have
\begin{eqnarray*}
&&I_{1}\leq C\int_{\xi=\sum\limits_{j=1}^{k+1}\xi_{j}}\int_{\tau=\sum\limits_{j=1}^{k+1}\tau_{j}}
\frac{f|\xi_{3}^{2}-\xi_{4}^{2}|^{\frac{1}{2}}\prod\limits_{j=2}^{k+1}\langle\xi_{j}\rangle^{-s}
(\prod\limits_{j=1}^{k+1}f_{j})}{\langle\sigma\rangle^{b_{1}}
\prod\limits_{j=1}^{k+1}\langle\sigma_{j}\rangle^{b}}d\delta\nonumber\\&&\leq C
\left\|F_{1}\right\|_{L_{xt}^{8}}\left\|F_{2}\right\|_{L_{xt}^{8}}\|F\|_{L_{xt}^{\frac{8}{1+\epsilon}}}
\|I^{\frac{1}{2}}(F_{3},F_{4})\|_{L_{xt}^{2}}
\prod\limits_{j=5}^{k+1}\|J^{-\frac{k-4}{2(k-3)}-\epsilon}F_{j}\|_{L_{xt}^{\frac{8(k-3)}{1-\epsilon}}}\nonumber\\
&&\leq C\left(\prod_{j=1}^{k+1}\|F_{j}\|_{X_{0,b}}\right)
\|F\|_{X_{0,b_{1}}}\leq C\left(\prod_{j=1}^{k+1}\|f_{j}\|_{L_{xt}^{2}}\right)\|f\|_{L_{xt}^{2}}.
\end{eqnarray*}
(5)Case$(\xi_{1},\xi_{2},\cdot\cdot\cdot,\xi_{k},\xi,\tau_{1},\tau_{2},\cdot\cdot\cdot,\tau_{k},\tau)\in \Omega_{4}$, we have
\begin{eqnarray}
&&\hspace{-1cm}K(\xi_{1},\xi_{2},\cdot\cdot\cdot,\xi_{k},\xi,\tau_{1},\tau_{2},\cdot\cdot\cdot,\tau_{k},\tau)\leq
\frac{|\xi_{4}^{2}-\xi_{5}^{2}|^{\frac{1}{2}}\langle\xi_{3}\rangle^{-\frac{k-4}{2(k-3)}-k\epsilon}
\prod\limits_{j=6}^{k+1}\langle\xi_{j}\rangle^{-\frac{k-4}{2(k-3)}-k\epsilon}}
{\langle\sigma\rangle^{b_{1}}
\prod\limits_{j=1}^{k+1}\langle\sigma_{j}\rangle^{b}}
.\label{3.08}
\end{eqnarray}
This case can be proved similarly to Case (4).

\noindent
(6)Case $(\xi_{1},\xi_{2},\cdot\cdot\cdot,\xi_{k},\xi,\tau_{1},\tau_{2},\cdot\cdot\cdot,\tau_{k},\tau)\in \Omega_{5}$,
we consider $k=5$, $k\geq6,$  respectively.

\noindent When $k=5,$  we have
\begin{eqnarray}
K(\xi_{1},\xi_{2},\cdot\cdot\cdot,\xi_{k},\xi,\tau_{1},\tau_{2},\cdot\cdot\cdot,\tau_{k},\tau)\leq C
\frac{\langle\xi_{6}\rangle^{-\frac{1}{3}-3\epsilon}\prod\limits_{j=1}^{5}|\xi_{j}|^{\frac{1}{6}}}
{\langle\sigma\rangle^{b_{1}}
\prod\limits_{j=1}^{6}\langle\sigma_{j}\rangle^{b}}
.\label{3.0888888}
\end{eqnarray}
By using (\ref{3.0888888}), the Plancherel identity, the  H\"older inequality
and   (\ref{2.05}), (\ref{2.06}) and (\ref{2.08}),
we have
\begin{eqnarray*}
&&I_{1}\leq C\int_{\xi=\sum\limits_{j=1}^{6}\xi_{j}}\int_{\tau=\sum\limits_{j=1}^{6}\tau_{j}}
\frac{ff_{6}(\prod\limits_{j=1}^{5}|\xi_{j}|^{\frac{1}{6}}f_{j})\langle\xi_{6}\rangle^{-\frac{1}{3}-3\epsilon}}
{\langle\sigma\rangle^{b_{1}}
\prod\limits_{j=1}^{6}\langle\sigma_{j}\rangle^{b}}d\delta\nonumber\\&&\leq C
\|F\|_{L_{xt}^{\frac{8}{1+\epsilon}}}\|J^{-\frac{1}{3}-3\epsilon}F_{6}\|_{L_{xt}^{\frac{24}{1-3\epsilon}}}
\prod\limits_{j=1}^{5}\left\|D_{x}^{\frac{1}{6}}F_{j}\right\|_{L_{xt}^{6}}
\nonumber\\&&\leq
C\left(\prod_{j=1}^{6}\|f_{j}\|_{L_{xt}^{2}}\right)\|f\|_{L_{xt}^{2}}
.
\end{eqnarray*}
When $k\geq6,$  we have
\begin{eqnarray}
&&\hspace{-1cm}K(\xi_{1},\xi_{2},\cdot\cdot\cdot,\xi_{k},\xi,\tau_{1},\tau_{2},\cdot\cdot\cdot,\tau_{k},\tau)\leq C
\frac{|\xi_{5}^{2}-\xi_{6}^{2}|^{\frac{1}{2}}\prod\limits_{j=1}^{4}\langle\xi_{j}
\rangle^{-s}\prod\limits_{j=7}^{k+1}\langle\xi_{j}\rangle^{-\frac{(k-4)^{2}}{2k(k-5)}-2\epsilon}}
{\langle\sigma\rangle^{b_{1}}\prod\limits_{j=1}^{k+1}\langle\sigma_{j}\rangle^{b}}
.\label{3.09}
\end{eqnarray}
By using (\ref{3.09}), the Plancherel identity, the  H\"older inequality
and   (\ref{2.05})-(\ref{2.08}),  Lemma 2.5,
we have
\begin{eqnarray*}
&&I_{1}\leq C\int_{\xi=\sum\limits_{j=1}^{k+1}\xi_{j}}\int_{\tau=\sum\limits_{j=1}^{k+1}\tau_{j}}
\frac{|\xi_{5}^{2}-\xi_{6}^{2}|^{\frac{1}{2}}f(\prod\limits_{j=1}^{k+1}f_{j})\prod\limits_{j=1}^{4}\langle\xi_{j}
\rangle^{-s}\prod\limits_{j=7}^{k+1}\langle\xi_{j}\rangle^{-\frac{(k-4)^{2}}{2k(k-5)}-2\epsilon}}
{\langle\sigma\rangle^{b_{1}}\prod\limits_{j=1}^{k+1}\langle\sigma_{j}\rangle^{b}}d\delta\nonumber\\&&\leq C
\|F\|_{L_{xt}^{\frac{8}{1+\epsilon}}}
\left(\prod\limits_{j=1}^{4}\left\|J^{-s}F_{j}\right\|_{L_{xt}^{2k}}\right)
\|I^{\frac{1}{2}}(F_{5},F_{6})\|_{L_{xt}^{2}}
\prod\limits_{j=7}^{k+1}\|J^{-\frac{(k-4)^{2}}{2k(k-5)}-2\epsilon}F_{j}\|_{L_{xt}^{\frac{8k(k-5)}{3k-16-k\epsilon}}}
\nonumber\\
&&\leq C\left(\prod_{j=1}^{k+1}\|F_{j}\|_{X_{0,b}}\right)
\|F\|_{X_{0,b_{1}}}\leq C\left(\prod_{j=1}^{k+1}\|f_{j}\|_{L_{xt}^{2}}\right)\|f\|_{L_{xt}^{2}}.
\end{eqnarray*}
(7)Case $(\xi_{1},\xi_{2},\cdot\cdot\cdot,\xi_{k},\xi,\tau_{1},\tau_{2},\cdot\cdot\cdot,\tau_{k},\tau)\in \Omega_{6}$,
this case can be proved similarly to
$(\xi_{1},\xi_{2},\cdot\cdot\cdot,\xi_{k},\xi,\tau_{1},\tau_{2},\cdot\cdot\cdot,\tau_{k},\tau)\in \Omega_{5}$.

\noindent(8)Case $(\xi_{1},\xi_{2},\cdot\cdot\cdot,\xi_{k},\xi,\tau_{1},\tau_{2},\cdot\cdot\cdot,\tau_{k},\tau)\in \Omega_{7}$, we have
\begin{eqnarray}
&&K(\xi_{1},\xi_{2},\cdot\cdot\cdot,\xi_{k},\xi,\tau_{1},\tau_{2},\cdot\cdot\cdot,\tau_{k},\tau)\leq C
\frac{\prod\limits_{j=1}^{k+1}\langle\xi_{j}\rangle^{\frac{1-ks}{k+1}}}
{\langle\sigma\rangle^{b_{1}}\prod\limits_{j=1}^{k+1}\langle\sigma_{j}\rangle^{b}}
.\label{3.010}
\end{eqnarray}
By using (\ref{3.010}), the Plancherel identity, the  H\"older inequality
and   (\ref{2.05}), (\ref{2.08}),
we have
\begin{eqnarray*}
&&I_{1}\leq C\int_{\xi=\sum\limits_{j=1}^{k+1}\xi_{j}}\int_{\tau=\sum\limits_{j=1}^{k+1}\tau_{j}}
\frac{f(\prod\limits_{j=1}^{k+1}f_{j})\prod\limits_{j=1}^{k+1}\langle\xi_{j}\rangle^{\frac{1-ks}{k+1}}}
{\langle\sigma\rangle^{b_{1}}
\prod\limits_{j=1}^{k+1}\langle\sigma_{j}\rangle^{b}}d\delta\nonumber\\&&\leq C
\|D_{x}^{\frac{3-\epsilon}{18}}F\|_{L_{xt}^{\frac{6}{1+\epsilon}}}
\left[\prod\limits_{j=1}^{2}\|D_{x}^{\frac{1}{6}}F_{j}\|_{L_{xt}^{6}}\right]
\left[\prod\limits_{j=3}^{k+1}\|D_{x}^{\frac{1-2ks}{2(k-1)}}F_{j}\|_{L_{xt}^{\frac{2(k-1)}{1-\epsilon}}}\right]
\nonumber\\
&&\leq C\left(\prod_{j=1}^{k+1}\|F_{j}\|_{X_{0,b}}\right)
\|F\|_{X_{0,\frac{1}{2}-\frac{\epsilon}{12}}}\leq C\left(\prod_{j=1}^{k+1}\|f_{j}\|_{L_{xt}^{2}}\right)\|f\|_{L_{xt}^{2}}.
\end{eqnarray*}

We have completed the proof of Lemma 3.1.

\begin{Lemma}\label{Lemma3.2}
Let $s\geq\frac{1}{2}-\frac{2}{k+1}+88\epsilon,k\geq5$ and $b=\frac{1}{2}+\frac{\epsilon}{24},b_{1}=\frac{1}{2}-\frac{\epsilon}{12}$.
Then, we have
\begin{eqnarray}
\left\|\partial_{x}(\prod\limits_{j=1}^{k+1}u_{j})\right\|_{X_{\frac{1}{2}+\epsilon,b_{1}}}\leq C\prod_{j=1}^{k+1}\|u_{j}\|_{X_{s,b}}.\label{3.001}
\end{eqnarray}
\end{Lemma}
\noindent {\bf Proof.}To prove (\ref{3.001}), by duality,  it suffices to prove
\begin{eqnarray}
\left|\int_{\SR^{2}}J^{\frac{1}{2}+\epsilon}\partial_{x}\left(\prod_{j=1}^{k+1}u_{j}\right)\bar{h}dxdt\right|\leq C\|h\|_{X_{0,b_{1}}}\prod_{j=1}^{k+1}\|u_{j}\|_{X_{s,b}}\label{3.002}.
\end{eqnarray}
We define
\begin{eqnarray*}
f(\xi,\tau)=\langle\sigma\rangle^{b_{1}} \mathscr{F}h(\xi,\tau),
f_{j}(\xi_{j},\tau_{j})=\langle\xi_{j}\rangle^{s}\langle\sigma_{j}\rangle^{b}
\mathscr{F}u_{j}(\xi_{j},\tau_{j})(1\leq j\leq k+1).
\end{eqnarray*}
By using the Plancherel identity, to prove (\ref{3.002}), it suffices to prove
\begin{eqnarray}
\int_{\xi=\sum\limits_{j=1}^{k+1}\xi_{j}}\int_{\tau=\sum\limits_{j=1}^{k+1}\tau_{j}}
\frac{|\xi|\langle\xi\rangle^{s}f\prod\limits_{j=1}^{k+1}f_{j}}
{\langle\sigma\rangle^{b_{1}}\prod\limits_{j=1}^{k+1}
\langle\xi_{j}\rangle^{s}\langle\sigma_{j}\rangle^{b}}d\delta
&&\leq C\|f\|_{L^{2}}\prod\limits_{j=1}^{k+1}\|f_{j}\|_{L^{2}},\label{3.003}
\end{eqnarray}
where $d\delta=d\xi_{1}d\xi_{2}\cdot \cdot \cdot d\xi_{k}d\xi d\tau_{1}d\tau_{2}\cdot \cdot \cdot d\tau_{k}d\tau.$

\noindent We define
\begin{eqnarray*}
&&K_{1}(\xi_{1},\xi_{2},\cdot\cdot\cdot,\xi_{k},\xi,\tau_{1},\tau_{2},\cdot\cdot\cdot,\tau_{k},\tau)
=\frac{|\xi|\langle\xi\rangle^{s}}
{\langle\sigma\rangle^{b_{1}}\prod\limits_{j=1}^{k+1}
\langle\xi_{j}\rangle^{s}\langle\sigma_{j}\rangle^{b}},\nonumber\\&&
\mathscr{F}F=\frac{f}{\langle\sigma\rangle^{b_{1}}},
\mathscr{F}F_{j}=\frac{f_{j}}{\langle\sigma_{j}\rangle^{b}}
(1\leq j\leq k+1),\\
&&I_{1}=\int_{\xi=\sum\limits_{j=1}^{k+1}\xi_{j}}
\int_{\tau=\sum\limits_{j=1}^{k+1}\tau_{j}}K(\xi_{1},\xi_{2},
\cdot\cdot\cdot,\xi_{k},\xi,\tau_{1},\tau_{2},\cdot\cdot\cdot,\tau_{k},\tau)
f\prod\limits_{j=1}^{k+1}f_{j}d\delta.
\end{eqnarray*}
Without loss of generality, we can assume that $|\xi_{1}|\geq
 |\xi_{2}|\geq\cdot\cdot\cdot\geq |\xi_{k+1}|.$

\noindent Obviously,
\begin{eqnarray*}
\Omega:=\left\{(\xi_{1},\xi_{2},\cdot\cdot\cdot,\xi_{k},\xi,\tau_{1},\tau_{2},
\cdot\cdot\cdot,\tau_{k},\tau)\in \R^{2(k+1)}:|\xi_{1}|\geq |\xi_{2}|\geq
\cdot\cdot\cdot\geq |\xi_{k+1}|\right\}\subset \bigcup\limits_{j=0}^{k+1}\Omega_{j}.
\end{eqnarray*}
Here, $\Omega_{j}(0\leq j\leq k+1)$ is defined as in Lemma 3.1.

\noindent(1)Case $(\xi_{1},\xi_{2},\cdot\cdot\cdot,\xi_{k},\xi,\tau_{1},\tau_{2},\cdot\cdot\cdot,
\tau_{k},\tau)\in \Omega_{0}$,
 since
 \begin{eqnarray*}
&&K_{1}(\xi_{1},\xi_{2},\cdot\cdot\cdot,\xi_{k},\xi,\tau_{1},\tau_{2},\cdot\cdot\cdot,\tau_{k},\tau)
\leq \frac{C}
{\langle\sigma\rangle^{b_{1}}\prod\limits_{j=1}^{k+1}
\langle\xi_{j}\rangle^{s}\langle\sigma_{j}\rangle^{b}},
 \end{eqnarray*}
 by using a proof similar to Case (1) of Lemma 3.1, we have
 \begin{eqnarray}
&&I_{1}\leq C\|f\|_{L^{2}}\prod\limits_{j=1}^{k+1}\|f_{j}\|_{L^{2}}\label{3.004}.
\end{eqnarray}
(2)Case  $(\xi_{1},\xi_{2},\cdot\cdot\cdot,\xi_{k},\xi,\tau_{1},\tau_{2},\cdot\cdot\cdot,\tau_{k},\tau)\in \Omega_{1}$,
we have
\begin{eqnarray}
&&K_{1}(\xi_{1},\xi_{2},\cdot\cdot\cdot,\xi_{k},\xi,\tau_{1},\tau_{2},\cdot\cdot\cdot,\tau_{k},\tau)\leq C
\frac{|\xi_{1}^{2}-\xi_{2}^{2}|^{\frac{1}{2}}|\xi_{1}|^{\frac{2}{k+1}+\epsilon}}
{\langle\sigma\rangle^{b_{1}}\prod\limits_{j=2}^{k+1}
\langle\xi_{j}\rangle^{s}\prod\limits_{j=1}^{k+1}\langle\sigma_{j}\rangle^{b}}
.\label{3.005}
\end{eqnarray}
By using  (\ref{3.005}), the Plancherel identity,  the H\"older inequality,
Lemma 2.2, (\ref{2.05}) and (\ref{2.08}), since $s\geq \frac{1}{2}-\frac{2}{k+1}+88\epsilon$,
we have
\begin{eqnarray*}
&&I_{1}\leq C\int_{\xi=\sum\limits_{j=1}^{k+1}\xi_{j}}\int_{\tau=\sum\limits_{j=1}^{k+1}\tau_{j}}
\frac{f(\prod\limits_{j=1}^{k+1}f_{j})|\xi_{1}|^{\frac{2}{k+1}+\epsilon}|\xi_{1}^{2}-\xi_{2}^{2}|^{\frac{1}{2}}\prod\limits_{j=2}^{k+1}
\langle\xi_{j}\rangle^{-s}}{\langle\sigma\rangle^{b_{1}}
\prod\limits_{j=1}^{k+1}\langle\sigma_{j}\rangle^{b}}d\delta\nonumber\\&&\leq C
\left\|I^{1/2}(F_{1},F_{2})\right\|_{L_{xt}^{2}}
\left(\prod_{j=3}^{k+1}\|J^{\frac{\frac{2}{k+1}+\epsilon-ks}{k-1}}F_{j}\|_{L_{xt}^{\frac{8(k-1)}{3-\epsilon}}}\right)
\|F\|_{L_{xt}^{\frac{8}{1+\epsilon}}}\nonumber\\
&&\leq C\|F\|_{X_{0,b_{1}}}\prod\limits_{j=1}^{k+1}\|F_{j}\|_{X_{0,b}}\leq
C\left(\prod_{j=1}^{k+1}\|f_{j}\|_{L_{xt}^{2}}\right)\|f\|_{L_{xt}^{2}}.
\end{eqnarray*}
(3)Case $(\xi_{1},\xi_{2},\cdot\cdot\cdot,\xi_{k},\xi,\tau_{1},\tau_{2},
\cdot\cdot\cdot,\tau_{k},\tau)\in \Omega_{2}$,
 since $s\geq \frac{1}{2}-\frac{2}{k+1}+88\epsilon,$ we  have
\begin{eqnarray}
&&K_{1}(\xi_{1},\xi_{2},\cdot\cdot\cdot,\xi_{k},\xi,\tau_{1},\tau_{2},\cdot\cdot\cdot,\tau_{k},\tau)\leq
C\frac{|\xi_{2}^{2}-\xi_{3}^{2}|^{\frac{1}{2}}}
{\langle\sigma\rangle^{b_{1}}\prod\limits_{j=2}^{k+1}
\langle\xi_{j}\rangle^{s}\prod\limits_{j=1}^{k+1}\langle\sigma_{j}\rangle^{b}}
\nonumber\\&&\leq C\frac{|\xi_{2}^{2}-\xi_{3}^{2}|^{\frac{1}{2}}
|\xi_{1}|^{\frac{1}{6}}\prod\limits_{j=4}^{k+1}
\langle\xi_{j}\rangle^{-\frac{ks+\frac{1}{6}-\frac{2}{k+1}-\epsilon}{6(k-2)}}}{\langle\sigma\rangle^{b_{1}}
\prod\limits_{j=1}^{k+1}\langle\sigma_{j}\rangle^{b}}.\label{3.006}
\end{eqnarray}
By using (\ref{3.006}), the Plancherel identity, the H\"older inequality,
  Lemma 2.2, (\ref{2.05})-(\ref{2.06}),
(\ref{2.08}),  we have
\begin{eqnarray*}
&&I_{1}\leq C\int_{\xi=\sum\limits_{j=1}^{k+1}\xi_{j}}\int_{\tau=\sum\limits_{j=1}^{k+1}\tau_{j}}
\frac{f(\prod\limits_{j=1}^{k+1}f_{j})|\xi_{2}^{2}-\xi_{3}^{2}|^{\frac{1}{2}}|\xi_{1}|^{\frac{1}{6}}
\prod\limits_{j=4}^{k+1}
\langle\xi_{j}\rangle^{-\frac{ks+\frac{1}{6}-\frac{2}{k+1}-\epsilon}{6(k-2)}}}{\langle\sigma\rangle^{b_{1}}
\prod\limits_{j=1}^{k+1}\langle\sigma_{j}\rangle^{b}}d\delta\nonumber\\&&\leq C
\left\|I^{1/2}(F_{2},F_{3})\right\|_{L_{xt}^{2}}
\left(\prod_{j=4}^{k+1}\|J^{-\frac{ks+\frac{1}{6}-\frac{2}{k+1}-\epsilon}{6(k-2)}}F_{j}\|_{L_{xt}^{\frac{24(k-2)}{5-3\epsilon}}}\right)
\|D_{x}^{\frac{1}{6}}F_{1}\|_{L_{xt}^{6}}\|F\|_{L_{xt}^{\frac{8}{1+\epsilon}}}\nonumber\\
&&\leq C\left(\prod_{j=1}^{k+1}\|F_{j}\|_{X_{0,b}}\right)
\|F\|_{X_{0,b_{1}}}\leq C\left(\prod_{j=1}^{k+1}\|f_{j}\|_{L_{xt}^{2}}\right)\|f\|_{L_{xt}^{2}}.
\end{eqnarray*}
(4)Case $(\xi_{1},\xi_{2},\cdot\cdot\cdot,\xi_{k},\xi,\tau_{1},
\tau_{2},\cdot\cdot\cdot,\tau_{k},\tau)\in \Omega_{3}$, we have
\begin{eqnarray}
&&K_{1}(\xi_{1},\xi_{2},\cdot\cdot\cdot,\xi_{k},\xi,\tau_{1},\tau_{2},\cdot\cdot\cdot,\tau_{k},\tau)\leq C
\frac{|\xi_{1}|^{\frac{2}{k+1}+\epsilon}|\xi_{3}^{2}-\xi_{4}^{2}|^{\frac{1}{2}}}
{\langle\sigma\rangle^{b_{1}}\prod\limits_{j=2}^{k+1}
\langle\xi_{j}\rangle^{s}\prod\limits_{j=1}^{k+1}\langle\sigma_{j}\rangle^{b}}
.\label{3.07}
\end{eqnarray}
 By using (\ref{3.07}),  the Plancherel identity, the  H\"older inequality and
  (\ref{2.05}), (\ref{2.07})-(\ref{2.08}) as well as Lemma 2.2,
we have
\begin{eqnarray*}
&&I_{1}\leq C\int_{\xi=\sum\limits_{j=1}^{k+1}\xi_{j}}\int_{\tau=\sum\limits_{j=1}^{k+1}\tau_{j}}
\frac{|\xi_{1}|^{\frac{2}{k+1}+\epsilon}f|\xi_{3}^{2}-\xi_{4}^{2}|^{\frac{1}{2}}\prod\limits_{j=2}^{k+1}\langle\xi_{j}\rangle^{-s}
(\prod\limits_{j=1}^{k+1}f_{j})}{\langle\sigma\rangle^{b_{1}}
\prod\limits_{j=1}^{k+1}\langle\sigma_{j}\rangle^{b}}d\delta\nonumber\\&&\leq C
\left\|F_{1}\right\|_{L_{xt}^{8}}\left\|F_{2}\right\|_{L_{xt}^{8}}\|F\|_{L_{xt}^{\frac{8}{1+\epsilon}}}
\|I^{\frac{1}{2}}(F_{3},F_{4})\|_{L_{xt}^{2}}
\prod\limits_{j=5}^{k+1}\|J^{-\frac{ks-\frac{2}{k+1}-\epsilon}{k-3}}F_{j}\|_{L_{xt}^{\frac{8(k-3)}{1-\epsilon}}}\nonumber\\
&&\leq C\left(\prod_{j=1}^{k+1}\|F_{j}\|_{X_{0,b}}\right)
\|F\|_{X_{0,b_{1}}}\leq C\left(\prod_{j=1}^{k+1}\|f_{j}\|_{L_{xt}^{2}}\right)\|f\|_{L_{xt}^{2}}.
\end{eqnarray*}
(5)Case$(\xi_{1},\xi_{2},\cdot\cdot\cdot,\xi_{k},\xi,\tau_{1},\tau_{2},\cdot\cdot\cdot,\tau_{k},\tau)\in \Omega_{4}$, we have
\begin{eqnarray}
&&\hspace{-1cm}K_{1}(\xi_{1},\xi_{2},\cdot\cdot\cdot,\xi_{k},\xi,\tau_{1},\tau_{2},\cdot\cdot\cdot,\tau_{k},\tau)\leq
\frac{|\xi_{4}^{2}-\xi_{5}^{2}|^{\frac{1}{2}}\langle\xi_{3}\rangle^{-\frac{ks-\frac{2}{k+1}-\epsilon}{k-3}}
\prod\limits_{j=6}^{k+1}\langle\xi_{j}\rangle^{-\frac{ks-\frac{2}{k+1}-\epsilon}{k-3}}}
{\langle\sigma\rangle^{b_{1}}
\prod\limits_{j=1}^{k+1}\langle\sigma_{j}\rangle^{b}}
.\label{3.08}
\end{eqnarray}
This case can be proved similarly to Case (4).

\noindent
(6)Case $(\xi_{1},\xi_{2},\cdot\cdot\cdot,\xi_{k},\xi,\tau_{1},\tau_{2},\cdot\cdot\cdot,\tau_{k},\tau)\in \Omega_{5}$,
we consider $k=5$, $k\geq6,$  respectively.

\noindent When $k=5,$  we have
\begin{eqnarray}
K_{1}(\xi_{1},\xi_{2},\cdot\cdot\cdot,\xi_{k},\xi,\tau_{1},\tau_{2},\cdot\cdot\cdot,\tau_{k},\tau)\leq C
\frac{\langle\xi_{6}\rangle^{-\frac{1}{3}-3\epsilon}\prod\limits_{j=1}^{5}|\xi_{j}|^{\frac{1}{6}}}
{\langle\sigma\rangle^{b_{1}}
\prod\limits_{j=1}^{6}\langle\sigma_{j}\rangle^{b}}
.\label{3.00888888}
\end{eqnarray}
By using (\ref{3.00888888}) and the proof similar to $k=5$ of Case (6),
we have
\begin{eqnarray*}
&&I_{1}\leq
C\left(\prod_{j=1}^{6}\|f_{j}\|_{L_{xt}^{2}}\right)\|f\|_{L_{xt}^{2}}
.
\end{eqnarray*}
When $k\geq6,$  we have
\begin{eqnarray}
&&\hspace{-1cm}K_{1}(\xi_{1},\xi_{2},\cdot\cdot\cdot,\xi_{k},\xi,\tau_{1},\tau_{2},\cdot\cdot\cdot,\tau_{k},\tau)\leq C
\frac{|\xi_{1}|^{\frac{2}{k+1}+\epsilon}|\xi_{3}^{2}-\xi_{4}^{2}|^{\frac{1}{2}}}
{\langle\sigma\rangle^{b_{1}}\prod\limits_{j=2}^{k+1}
\langle\xi_{j}\rangle^{s}\prod\limits_{j=1}^{k+1}\langle\sigma_{j}\rangle^{b}}
.\label{3.009}
\end{eqnarray}
By using (\ref{3.009}) and the proof similar to Case (4),
we have
\begin{eqnarray*}
&&I_{1}\leq C\left(\prod_{j=1}^{k+1}\|f_{j}\|_{L_{xt}^{2}}\right)\|f\|_{L_{xt}^{2}}.
\end{eqnarray*}
(7)Case $(\xi_{1},\xi_{2},\cdot\cdot\cdot,\xi_{k},\xi,\tau_{1},\tau_{2},\cdot\cdot\cdot,\tau_{k},\tau)\in \Omega_{6}$,
this case can be proved similarly to
$(\xi_{1},\xi_{2},\cdot\cdot\cdot,\xi_{k},\xi,\tau_{1},\tau_{2},\cdot\cdot\cdot,\tau_{k},\tau)\in \Omega_{5}$.

\noindent(8)Case $(\xi_{1},\xi_{2},\cdot\cdot\cdot,\xi_{k},\xi,\tau_{1},\tau_{2},\cdot\cdot\cdot,\tau_{k},\tau)\in \Omega_{7}$, we have
\begin{eqnarray}
&&K_{1}(\xi_{1},\xi_{2},\cdot\cdot\cdot,\xi_{k},\xi,\tau_{1},\tau_{2},\cdot\cdot\cdot,\tau_{k},\tau)\leq C
\frac{\prod\limits_{j=1}^{k+1}\langle\xi_{j}\rangle^{\frac{1-ks}{k+1}}}
{\langle\sigma\rangle^{b_{1}}\prod\limits_{j=1}^{k+1}\langle\sigma_{j}\rangle^{b}}
.\label{3.0010}
\end{eqnarray}
By using (\ref{3.0010}) and  the proof similar to Case (8) of Lemma 3.1,
we have
\begin{eqnarray*}
&&I_{1}\leq  C\left(\prod_{j=1}^{k+1}\|f_{j}\|_{L_{xt}^{2}}\right)\|f\|_{L_{xt}^{2}}.
\end{eqnarray*}

We have completed the proof of Lemma 3.2.

\begin{Lemma}\label{Lemma3.3}
Let
  $s\geq\frac{1}{2}-\frac{2}{k}+88\epsilon$ and $s_{1}\geq{\rm max}\left\{\frac{s}{k+1},s-\frac{1}{3},\frac{s-\frac{1}{3}}{k}\right\}+88\epsilon$,
   $s_{2}=-1-s_{1}+s-8\epsilon(<0)$, $s_{3}=s-(k+1)s_{1}-8\epsilon(<0)$,  $s_{4}=s-ks_{1}-\frac{1}{3}-8\epsilon(<0)$ and
   $z_{j}(t)=\phi(t)U(t)f^{\omega}(1\leq j\leq k+1,j\in N)$ and $b=\frac{1}{2}+\frac{\epsilon}{24},b_{1}=\frac{1}{2}-\frac{\epsilon}{12}$.
   We  define $\sum=\sum\limits_{N_{1}, \cdot\cdot\cdot,  N_{k+1},N} $.
 Then, we have
\begin{eqnarray}
\left\|\partial_{x}\left(\prod_{j=1}^{k+1}z_{j}\right)\right\|_{X_{s,-\frac{1}{2}
+\frac{\epsilon}{12}}}\leq C\lambda^{k+1}\label{3.011}
\end{eqnarray}
outside a set of probability
$
\leq \epsilon.
$
\end{Lemma}
\noindent{\bf Proof.}To prove (\ref{3.011}), by duality,  it suffices to prove
\begin{eqnarray}
\left|\int_{\SR^{2}}J^{s}\partial_{x}\left(\prod_{j=1}^{k+1}z_{j}\right)\bar{h}dxdt\right|\leq C\lambda^{k+1}\|h\|_{X_{0,\frac{1}{2}-\frac{\epsilon}{12}}}\label{3.012}
\end{eqnarray}
outside a set of probability
$
\leq \epsilon.
$
We define
   \begin{eqnarray*}
   I_{1}=\left|\int_{\SR^{2}}J^{s}\partial_{x}\left(\prod_{j=1}^{k+1}z_{j}\right)\bar{h}dxdt\right|
=C\sum
\left|\int_{\SR^{2}}J^{s}\partial_{x}\left(\prod_{j=1}^{k+1}P_{N_{j}}z_{j}\right)P_{N}\bar{h}dxdt\right|.
   \end{eqnarray*}
 We dyadically decompose $z_{j}$ with $(1\leq j\leq k+1,j\in Z)$  and $h$ such that frequency supports are
$\left\{|\xi_{j}|\sim N_{j}\right\}(1\leq j\leq N)$ and $\left\{|\xi|\sim N\right\}$ for some dyadically $N_{j}$.
Without loss of generality, we can assume that $N_{1}\geq N_{2}\geq \cdot\cdot\cdot\geq N_{k}\geq N_{k+1}$.

\noindent (1).Case $N_{1}\leq 80k.$ By using the  H\"older inequality,   Lemma 2.12  and   the Sobolev embeddings Theorem
as well as (\ref{2.035}),
we have
\begin{eqnarray*}
&&I_{1}\leq C\sum N^{1+s}\|P_{N}h\|_{L_{xt}^{\frac{8}{1+\epsilon}}}
\left[\prod\limits_{j=1}^{k+1}\|P_{N}z_{j}\|_{L_{xt}^{\frac{8(k+1)}{7-\epsilon}}}\right]
\nonumber\\
&&\leq C\sum N^{1+s}
\|P_{N}h\|_{L_{xt}^{\frac{8}{1+\epsilon}}}
\left[\prod\limits_{j=1}^{k+1}\|P_{N}z_{j}\|_{L_{tx}^{\frac{8(k+1)}{7-\epsilon}}}\right]\nonumber\\
&&\leq C\sum N^{1+s}\left(\prod\limits_{j=1}^{k+1}N_{j}^{\epsilon}\right)
\|P_{N}h\|_{L_{xt}^{\frac{8}{1+\epsilon}}}
\left[\prod\limits_{j=1}^{k+1}\|P_{N}z_{j}\|_{L_{t}^{\frac{8(k+1)}{7-\epsilon}}L_{x}^{\frac{8(k+1)}{7+(8k+7)\epsilon}}}\right]\nonumber\\
&&\leq C\|P_{N}h\|_{\ell_{N_{j}}^{\frac{8}{1+\epsilon}}L_{xt}^{\frac{8}{1+\epsilon}}}\prod\limits_{j=1}^{k+1}
\|P_{N_{j}}z_{j}\|_{\ell_{N_{j}}^{\frac{8(k+1)}{7-\epsilon}}L_{t}^{\frac{8(k+1)}{7-\epsilon}}L_{x}^{\frac{8(k+1)}{7+(8k+7)\epsilon}}}\nonumber\\
&&\leq C\|h\|_{L_{xt}^{\frac{8}{1+\epsilon}}}\prod\limits_{j=1}^{k+1}
\|z_{j}\|_{L_{t}^{\frac{8(k+1)}{7-\epsilon}}L_{x}^{\frac{8(k+1)}{7+(8k+7)\epsilon}}}\leq C\lambda^{k+1}
\|h\|_{X_{0,b_{1}}}
\end{eqnarray*}
outside a set of probability
$
\leq C_{1}e^{-C\lambda^{2}\|f\|_{H^{s}}^{-2}}.
$

\noindent (2).Case $N_{1}\geq 80k{\rm max}\left\{1, N_{2}\right\}$ and $N_{4}\leq 1$.
By using the
 H\"older inequality and   Lemmas 2.2, 2.8, 2.12 as well as  (\ref{2.035}),
  since $s_{1}\geq{\rm max}\left\{\frac{s}{k+1},s-\frac{1}{3},\frac{s-\frac{1}{3}}{k}\right\}+88\epsilon,$
we have
\begin{eqnarray*}
&&I_{1}\leq C\sum N^{s_{2}}\left\|I^{\frac{1}{2}}(J^{s_{1}}P_{N_{1}}z_{1},P_{N_{2}}z_{2})\right\|_{L_{xt}^{2}}
\left[\prod\limits_{j=4}^{k+1}\|P_{N_{j}}z_{j}\|_{L_{xt}^{\infty}}\right]
\left\|I^{\frac{1}{2}-2\epsilon}(P_{N_{3}}z_{3},P_{N}h)\right\|_{L_{xt}^{2}}
\nonumber\\&&\leq C\sum N^{s_{2}}
\left(\prod\limits_{j=2}^{k+1}N_{j}^{\epsilon}\right)
\left(\prod\limits_{j=4}^{k+1}N_{j}^{\frac{1}{2}-\frac{k+1}{k-2}\epsilon}\right)
\left[\prod\limits_{j=1}^{k+1}\|P_{N_{j}}z_{j}\|_{X_{s_{1},b}}\right]
\|P_{N}h\|_{X_{0,(1-\epsilon)b}}\nonumber\\
&&\leq C\sum N^{s_{2}}
\left(\prod\limits_{j=2}^{k+1}N_{j}^{\epsilon}\right)
\left[\prod\limits_{j=1}^{k+1}\|P_{N_{j}}z_{j}\|_{X_{s_{1},b}}\right]\|P_{N}h\|_{X_{0,\frac{1}{2}-\frac{\epsilon}{12}}}\nonumber\\
&&\leq C\left[\prod\limits_{j=1}^{k+1}\|P_{N_{j}}z_{j}\|_{\ell_{N_{j}}^{2}X_{s_{1},b}}\right]
\|P_{N}h\|_{\ell_{N}^{2}X_{0,\frac{1}{2}-\frac{\epsilon}{12}}}\nonumber\\
&&\leq C\left[\prod\limits_{j=1}^{k+1}\|z_{j}\|_{X_{s_{1},b}}\right]
\|h\|_{X_{0,b_{1}}}\leq
 C\lambda^{k+1}\|h\|_{X_{0,b_{1}}}
\end{eqnarray*}
outside a set of probability
$
\leq C_{1}e^{-C\lambda^{2}\|f\|_{H^{s}}^{-2}}.
$

\noindent (3).Case  $N_{1}\geq 80k{\rm max}\left\{ 1,N_{2}\right\}$ and
$N_{2}\geq N_{3}\geq\cdot\cdot\cdot\geq N_{l}\geq 1\geq \cdot\cdot\cdot\geq N_{k+1}$.
By using the
 H\"older inequality and   Lemmas 2.2, 2.8, 2.12  as well as  (\ref{2.035}),
   since $s_{1}\geq{\rm max}\left\{\frac{s}{k+1},s-\frac{1}{3},\frac{s-\frac{1}{3}}{k}\right\}+88\epsilon,$
we have
\begin{eqnarray*}
&&I_{1}\leq C\sum N^{s_{2}}\left\|I^{\frac{1}{2}}(J^{s_{1}}P_{N_{1}}z_{1},P_{N_{2}}z_{2})\right\|_{L_{xt}^{2}}
\left[\prod\limits_{j=4}^{k+1}\|P_{N_{j}}z_{j}\|_{L_{xt}^{\infty}}\right]
\left\|I^{\frac{1}{2}-2\epsilon}(P_{N_{3}}z_{3},P_{N}h)\right\|_{L_{xt}^{2}}
\nonumber\\&&\leq C\lambda^{l-3}\sum N^{s_{2}}
\left[\prod\limits_{j=1}^{3}\|P_{N_{j}}z_{j}\|_{X_{s_{1},b}}\right]
\left[\prod\limits_{j=l+1}^{k+1}N_{j}^{\frac{1}{2}}\|P_{N_{j}}z_{j}\|_{X_{0,b}}\right]
\|P_{N}h\|_{X_{0,(1-\epsilon)b}}\nonumber\\
&&\leq C\lambda^{\ell-3}\left[\prod\limits_{j=1}^{3}
\|P_{N_{j}}z_{j}\|_{\ell_{N_{j}}^{2}X_{s_{1},b}}\right]
\left[\prod\limits_{j=l+1}^{k+1}
\|P_{N_{j}}z_{j}\|_{\ell_{N_{j}}^{2}X_{0,b}}\right]
\|P_{N}h\|_{\ell_{N}^{2}X_{0,(1-\epsilon)b}}\nonumber\\
&&\leq C\lambda^{\ell-3}\left[\prod\limits_{j=1}^{k+1}\|z_{j}\|_{X_{s_{1},b}}\right]
\|h\|_{X_{0,(1-\epsilon)b}}\leq C\lambda^{k+1}
\|h\|_{X_{0,b_{1}}}
\end{eqnarray*}
outside a set of probability
$
\leq C_{1}e^{-C\lambda^{2}\|f\|_{H^{s}}^{-2}}.
$

\noindent (4).Case  $N_{1}\geq80k, N_{1}\sim N_{2}\geq 80kN_{3},N_{3}\leq 1.$ By using the proof
 similar to Case (3) of Lemma 3.2,
 we have
\begin{eqnarray*}
&&\left|\int_{\SR^{2}}J^{s}\partial_{x}\left(\prod_{j=1}^{k+1}z_{j}\right)\bar{h}dxdt\right|\leq C\lambda^{k+1}\|h\|_{X_{0,b_{1}}}
\end{eqnarray*}
outside a set of probability
$
\leq C_{1}e^{-C\lambda^{2}\|f\|_{H^{s}}^{-2}}.
$

\noindent (5). Case  $N_{1}\geq80k, N_{1}\sim N_{l}\geq 80kN_{l+1}(3\leq l\leq k-1).$
By using the proof similar to Case (3) of Lemma 3.2,
 we have
\begin{eqnarray*}
&&\left|\int_{\SR^{2}}J^{s}\partial_{x}\left(\prod_{j=1}^{k+1}z_{j}\right)\bar{h}dxdt\right|
\leq C\lambda^{k+1}\|h\|_{X_{0,b_{1}}}
\end{eqnarray*}
outside a set of probability
$
\leq C_{1}e^{-C\lambda^{2}\|f\|_{H^{s}}^{-2}}.
$

\noindent (6). Case  $N_{1}\geq80k, N_{1}\sim N_{k}\geq 80kN_{k+1}\geq 80k.$
By using the
 H\"older inequality and   Lemmas 2.2, 2.12, 2.13  as well as  (\ref{2.035}),
   since $s_{1}\geq{\rm max}\left\{\frac{s}{k+1},s-\frac{1}{3},\frac{s-\frac{1}{3}}{k}\right\}+88\epsilon,$
we have
\begin{eqnarray*}
&&I_{1}\leq C\sum N_{1}^{s_{4}}N^{\epsilon}\left\|I^{\frac{1}{2}}(J^{s_{1}}P_{N_{k}}z_{k},P_{N_{k+1}}z_{k+1})\right\|_{L_{xt}^{2}}
\left[\prod\limits_{j=1}^{2}\|J^{s_{1}}D_{x}^{\frac{1}{6}}P_{N_{j}}z_{j}\|_{L_{xt}^{6}}\right]\nonumber\\&&
\qquad\qquad\qquad\qquad\qquad\qquad\times\left[\prod\limits_{j=3}^{k-1}
\|J^{s_{1}}P_{N_{j}}z_{j}\|_{L_{xt}^{\frac{24(k-3)}{1-3\epsilon}}}\right]
\left\|P_{N}h\right\|_{L_{xt}^{\frac{8}{1+\epsilon}}}
\nonumber\\&&\leq C\sum N_{1}^{s_{4}}N^{\epsilon}
\left[\prod\limits_{j=1}^{2}\|P_{N_{j}}z_{j}\|_{X_{s_{1},b}}\right]
\left[\prod\limits_{j=3}^{k-1}
\|J^{s_{1}}P_{N_{j}}z_{j}\|_{L_{xt}^{\frac{24(k-3)}{1-3\epsilon}}}\right]\nonumber\\&&
\qquad\qquad\qquad\qquad\qquad\qquad\times
\left[\prod\limits_{j=k}^{k+1}\|P_{N_{j}}z_{j}\|_{X_{s_{1},b}}\right]
\|P_{N}h\|_{X_{0,\frac{1}{2}-\frac{\epsilon}{12}}}\nonumber\\
&&\leq C\prod\limits_{j=1}^{2}
\|P_{N_{j}}z_{j}\|_{\ell_{N_{j}}^{2}X_{s_{1},b}}\prod\limits_{j=3}^{k-1}
\|J^{s_{1}}P_{N_{j}}z_{j}\|_{\ell_{N_{j}}^{\frac{24(k-3)}{1-3\epsilon}}
L_{xt}^{\frac{24(k-3)}{1-3\epsilon}}}\nonumber\\
&&\qquad\qquad\qquad\qquad\qquad\qquad\times
\left[\prod\limits_{j=k}^{k+1}
\|P_{N_{j}}z_{j}\|_{\ell_{N_{j}}^{2}X_{s_{1},b}}\right]
\|P_{N}h\|_{\ell_{N}^{2}X_{0,\frac{1}{2}-\frac{\epsilon}{12}}}\nonumber\\
&&\leq C\left[\prod\limits_{j=1}^{2}\|z_{j}\|_{X_{s_{1},b}}\right]
\prod\limits_{j=3}^{k-1}\left[\|J^{s_{1}}z_{j}\|_{L_{xt}^{\frac{24(k-3)}{1-3\epsilon}}}\right]
\left[\prod\limits_{j=k}^{k+1}\|z_{j}\|_{X_{0,b}}\right]
\|h\|_{X_{0,b_{1}}}\nonumber\\&&\leq C\lambda^{k+1}
\|h\|_{X_{0,b_{1}}}
\end{eqnarray*}
outside a set of probability
$
\leq C_{1}e^{-C\lambda^{2}\|f\|_{H^{s}}^{-2}}.
$

\noindent (7). Case  $N_{1}\geq80k, N_{1}\sim N_{k}\geq 80kN_{k+1},N_{k+1}\leq1,N_{1}N_{k+1}\geq1.$
By using the
 H\"older inequality,    Lemmas 2.2, 2.8, 2.13 and   (\ref{2.035}),
   since $s_{1}\geq{\rm max}\left\{\frac{s}{k+1},s-\frac{1}{3},\frac{s-\frac{1}{3}}{k}\right\}+88\epsilon,$
we have
\begin{eqnarray*}
&&I_{1}\leq C\sum N_{1}^{s_{4}}(N_{k+1}N_{1})^{\epsilon}N^{\epsilon}
\left[\prod\limits_{j=1}^{2}\|D_{x}^{s_{1}+\frac{1-8\epsilon}{6}}P_{N_{j}}z_{j}\|_{L_{xt}^{\frac{6}{1-2\epsilon}}}\right]
\left[\prod\limits_{j=3}^{k-1}\|J^{s_{1}}P_{N_{j}}z_{j}\|_{L_{xt}^{\frac{3(k-3)}{\epsilon}}}\right]\nonumber\\&&
\qquad\qquad\qquad\qquad\times\left\|I^{\frac{1}{2}}\left(P_{N_{k+1}}z_{k+1},J^{s_{1}}P_{N_{k}}z_{k}\right)\right\|_{L_{xt}^{2}}
\|D_{x}^{\frac{3-\epsilon}{18}}P_{N}h\|_{L_{xt}^{\frac{6}{1+\epsilon}}}\nonumber\\
&&\leq C\sum N_{1}^{s_{4}+\epsilon}(NN_{k+1})^{\epsilon} \left[\prod\limits_{j=1}^{2}\|z_{j}\|_{X_{s_{1},b}}\right]
\left[\prod\limits_{j=3}^{k-1}\|J^{s_{1}}P_{N_{j}}z_{j}\|_{L_{xt}^{\frac{3(k-3)}{\epsilon}}}\right]\nonumber\\&&
\qquad\qquad\qquad\qquad\times\left[\prod\limits_{j=k}^{k+1}\|z_{j}\|_{X_{s_{1},b}}\right]
\|P_{N}h\|_{X_{0,\frac{1}{2}-\frac{\epsilon}{12}}}\nonumber\\
&&\leq C\lambda^{k+1}\sum N_{1}^{s_{4}+\epsilon}(NN_{k+1})^{\epsilon}
\|P_{N}h\|_{X_{0,\frac{1}{2}-\frac{\epsilon}{12}}}
\leq
 C\lambda^{k+1}\|h\|_{X_{0,\frac{1}{2}-\frac{\epsilon}{12}}}
\end{eqnarray*}
outside a set of probability
$
\leq C_{1}e^{-C\lambda^{2}\|f\|_{H^{s}}^{-2}}.
$

\noindent (8). Case  $N_{1}\geq80k, N_{1}\sim N_{k}\geq 80kN_{k+1},N_{k+1}\leq1,N_{1}N_{k+1}\leq1.$
By using the
 H\"older inequality,    Lemmas 2.2, 2.8, 2.13 and   (\ref{2.035}),
   since $s_{1}\geq{\rm max}\left\{\frac{s}{k+1},s-\frac{1}{3},\frac{s-\frac{1}{3}}{k}\right\}+88\epsilon,$
we have
\begin{eqnarray*}
&&I_{1}\leq C\sum N_{1}^{s-ks_{1}+9\epsilon}N^{\epsilon}
\left(\prod\limits_{j=1}^{4}\|D_{x}^{s_{1}+\frac{1}{6}}P_{N_{j}}z_{j}\|_{L_{xt}^{6}}\right)
\|D_{x}^{s_{1}+\frac{1-8\epsilon}{6}}P_{N_{5}}z_{5}\|_{L_{xt}^{\frac{6}{1-2\epsilon}}}
\nonumber\\&&\qquad\qquad\qquad\qquad\times
\left[\prod\limits_{j=6}^{k}\|D_{x}^{s_{1}}P_{N_{j}}z_{j}\|_{L_{xt}^{\frac{6(k-5)}{\epsilon}}}\right]
\|P_{N_{k+1}}z_{k+1}\|_{L_{xt}^{\infty}}
\left\|D_{x}^{\frac{3-\epsilon}{18}}P_{N}h\right\|_{L_{xt}^{\frac{6}{1+\epsilon}}}
\nonumber\\&&\leq C\sum N_{1}^{s-ks_{1}+9\epsilon}N^{\epsilon}N_{k+1}^{\frac{1}{2}}
\left(\prod\limits_{j=1}^{5}\|P_{N_{j}}z_{j}\|_{X_{s_{1},b}}\right)
\left[\prod\limits_{j=6}^{k}\|D_{x}^{s_{1}}P_{N_{j}}z_{j}\|_{L_{xt}^{\frac{6(k-5)}{\epsilon}}}\right]
\nonumber\\&&\qquad\qquad\qquad\qquad\times\|P_{N_{k+1}}z_{k+1}\|_{X_{0,b}}
\left\|P_{N}h\right\|_{X_{0,b_{1}}}\nonumber\\
&&\leq C\lambda^{k+1}\sum N_{1}^{s_{4}}N^{\epsilon}N_{k+1}^{\epsilon}
\left\|P_{N}h\right\|_{X_{0,b_{1}}}
\leq
 C\lambda^{k+1}\|h\|_{X_{0,\frac{1}{2}-\frac{\epsilon}{12}}}
\end{eqnarray*}
outside a set of probability
$
\leq C_{1}e^{-C\lambda^{2}\|f\|_{H^{s}}^{-2}}.
$

\noindent (9).Case $N_{1}\geq80k, N_{1}\sim N_{k+1}.$
By using the H\"older inequality and the Cauchy-Schwarz
inequality as well as (\ref{2.035}), since $s_{1}\geq{\rm max}\left\{\frac{s}{k+1},s-\frac{1}{3},\frac{s-\frac{1}{3}}{k}\right\}+88\epsilon,$ we have
\begin{eqnarray*}
&&I_{1}\leq C\sum N_{1}^{s_{3}}N^{\epsilon}
\left(\prod\limits_{j=1}^{4}\|D_{x}^{s_{1}+\frac{1}{6}}P_{N_{j}}z_{j}\|_{L_{xt}^{6}}\right)
\|D_{x}^{s_{1}+\frac{1-8\epsilon}{6}}P_{N_{5}}z_{5}\|_{L_{xt}^{\frac{6}{1-2\epsilon}}}
\nonumber\\&&\qquad\qquad\qquad\qquad\times\left[\prod\limits_{j=6}^{k+1}\|D_{x}^{s_{1}}P_{N_{j}}z_{j}\|_{L_{xt}^{\frac{6(k-4)}{\epsilon}}}\right]
\left\|D_{x}^{\frac{3-\epsilon}{18}}P_{N}h\right\|_{L_{xt}^{\frac{6}{1+\epsilon}}}
\nonumber\\&&\leq C\sum N_{1}^{s_{3}}N^{\epsilon}
\left(\prod\limits_{j=1}^{5}\|P_{N_{j}}z_{j}\|_{X_{s_{1},b}}\right)
\prod\limits_{j=6}^{k+1}\|D_{x}^{s_{1}}P_{N_{j}}z_{j}\|_{L_{xt}^{\frac{6(k-4)}{\epsilon}}}
\left\|P_{N}h\right\|_{X_{0,b_{1}}}\nonumber\\
&&\leq C\lambda^{k+1}\sum N_{1}^{s_{3}}N^{\epsilon}\left\|P_{N}h\right\|_{X_{0,b_{1}}}
\leq
 C\lambda^{k+1}\|h\|_{X_{0,b_{1}}}
\end{eqnarray*}
outside a set of probability
$
\leq C_{1}e^{-C\lambda^{2}\|f\|_{H^{s}}^{-2}}.
$

We have completed the proof of Lemma 3.3.

\noindent {\bf Remark 8:}
Case (10) of Lemma 3.2 requires  $(k+1)s_{1}\geq s+88(k+1)\epsilon.$
\begin{Lemma}\label{Lemma3.3}
Let $s\geq \frac{1}{2}-\frac{2}{k}+88\epsilon$
 and $s_{1}\geq{\rm max}
\left\{\frac{s}{k+1},s-\frac{1}{3},\frac{s-\frac{1}{3}}{k}\right\}+88\epsilon,$
$s_{2}=-1-s_{1}+s-8\epsilon(<0),$ $s_{3}=s-(k+1)s_{1}-8\epsilon(<0)$,
 $s_{4}=s-s_{1}-\frac{1}{3}-8\epsilon(<0),$
 $z_{k+1}=\phi(t)U(t)f^{\omega}(x)$
 and   $b=\frac{1}{2}+\frac{\epsilon}{24},b_{1}=\frac{1}{2}-\frac{\epsilon}{12}$.
  We  define $\sum=\sum\limits_{N_{1}, \cdot\cdot\cdot,  N_{k+1},N} $.
Then, we have
\begin{eqnarray}
\left\|\partial_{x}\left(z_{k+1}\prod_{j=1}^{k}v_{j}\right)
\right\|_{X_{s,-\frac{1}{2}+\frac{\epsilon}{12}}}\leq C
\left(\prod_{j=1}^{k}\|v_{j}\|_{X_{s,b}}\right)\lambda
\label{3.013}
\end{eqnarray}
outside a set of probability
$
\leq C_{1}e^{-C\lambda^{2}\|f\|_{H^{s}}^{-2}}.
$
\end{Lemma}
\noindent {\bf Proof.}
To prove (\ref{3.013}), by duality, it suffices to prove
\begin{eqnarray}
\left|\int_{\SR^{2}}J^{s}\partial_{x}\left(z_{k+1}\prod_{j=1}^{k}v_{j}\right)\bar{h}dxdt\right|\leq
C\lambda \left(\prod_{j=1}^{k}\|v_{j}\|_{X_{s,b}}\right) \|h\|_{X_{0,\frac{1}{2}-\frac{\epsilon}{12}}}\label{3.014}
\end{eqnarray}
outside a set of probability
$
\leq C_{1}e^{-C\lambda^{2}\|f\|_{H^{s}}^{-2}}.
$
We define
   \begin{eqnarray*}
   I_{2}=\left|\int_{\SR^{2}}J^{s}\partial_{x}\left(z_{k+1}\prod_{j=1}^{k}v_{j}\right)\bar{h}dxdt\right|
=C\sum
\left|\int_{\SR^{2}}J^{s}\partial_{x}\left(P_{N_{k+1}}z_{k+1}\prod_{j=1}^{k}P_{N_{j}}v_{j}\right)P_{N}\bar{h}dxdt\right|.
   \end{eqnarray*}
We dyadically decompose $v_{j}$ with $(1\leq j\leq k,j\in Z)$, $z_{k+1}$
  and $h$ such that frequency supports are
$\left\{|\xi_{j}|\sim N_{j}\right\}(1\leq j\leq k+1)$  and $\left\{|\xi|\sim N\right\}$
for some dyadically $N_{j}(1\leq j\leq k+1),N$.
Without loss of generality, we can assume that $N_{1}\geq N_{2}\geq N_{3}\geq N_{4}
\geq \cdot\cdot\cdot \geq N_{k}$.

\noindent (1). Case  ${\rm max}\left\{N_{1},N_{k+1}\right\}\leq 80k.$
By using the  H\"older inequality,   (\ref{2.035}), (\ref{2.040}) and the Sobolev embeddings Theorem,
we have
\begin{eqnarray*}
&&I_{1}\leq C\sum N^{1+s}\|P_{N}h\|_{L_{xt}^{2}}
\left[\prod\limits_{j=1}^{k}\|P_{N}v_{j}\|_{L_{xt}^{2(k+1)}}\right]\|P_{N_{k+1}}z_{k+1}\|_{L_{xt}^{2(k+1)}}
\nonumber\\
&&\leq C\sum N^{1+s}(\prod\limits_{j=1}^{k+1}N_{j}^{\epsilon})
\|P_{N}h\|_{L_{xt}^{2}}
\left[\prod\limits_{j=1}^{k}\|P_{N}v_{j}\|_{L_{xt}^{2(k+1)}}\right]
\|P_{N_{k+1}}z_{k+1}\|_{L_{t}^{2(k+1)}L_{x}^{\frac{2(k+1)}{1+(2k+1)\epsilon}}}\nonumber\\
&&\leq C\sum N^{1+s}(\prod\limits_{j=1}^{k+1}N_{j}^{\epsilon})
\|P_{N}h\|_{L_{xt}^{2}}
\left[\prod\limits_{j=1}^{k}\|P_{N_{j}}v_{j}\|_{X_{s,b}}\right]
\|P_{N_{k+1}}z_{k+1}\|_{L_{t}^{2(k+1)}L_{x}^{\frac{2(k+1)}{1+(2k+1)\epsilon}}}\nonumber\\
&&\leq C\|P_{N}h\|_{\ell_{N}^{2}L_{xt}^{2}}
\left[\prod\limits_{j=1}^{k}\|P_{N_{j}}z_{j}\|_{\ell_{N_{j}}^{2}X_{s,b}}\right]
\|P_{N_{k+1}}z_{k+1}\|_{\ell_{N_{j}}^{2(k+1)}L_{t}^{2(k+1)}L_{x}^{\frac{2(k+1)}{1+(2k+1)\epsilon}}}\nonumber\\
&&\leq C\|h\|_{L_{xt}^{2}}\left[\prod\limits_{j=1}^{k}
\|v_{j}\|_{X_{s,b}}\right]\|z_{k+1}\|_{L_{t}^{2(k+1)}L_{x}^{\frac{2(k+1)}{1+(2k+1)\epsilon}}}\leq C\lambda\left(\prod_{j=1}^{k}\|v_{j}\|_{X_{s,b}}\right)
\|h\|_{X_{0,b_{1}}}
\end{eqnarray*}
outside a set of probability
$
\leq C_{1}e^{-C\lambda^{2}\|f\|_{H^{s}}^{-2}}.
$

\noindent
(2). Case $N_{1}\geq 80k{\rm max}\left\{1,N_{2}\right\},N_{1}\geq \frac{N_{k+1}}{80k}$.
By using the  H\"older inequality,   Lemma 2.2,  (\ref{2.05}), (\ref{2.09}),
the Sobolev embeddings Theorem  and Lemma 2.13 as well as (\ref{2.035}),
we have
\begin{eqnarray*}
&&I_{2}\leq C\sum N^{-\epsilon}
\left\|I^{\frac{1}{2}}(J^{s}P_{N_{1}}v_{1},J^{s}P_{N_{2}}v_{2})\right\|_{L_{xt}^{2}}
\nonumber\\&&\qquad\qquad\times
\left(\prod_{j=3}^{k}\|J_{x}^{-\frac{s+\frac{1}{6}}{k-2}}P_{N_{j}}v_{j}\|_{L_{xt}^{\frac{6(k-2)}{2-7\epsilon}}}\right)
\|P_{N_{k+1}}z_{k+1}\|_{L_{xt}^{\frac{1}{\epsilon}}}
\|D_{x}^{\frac{3-\epsilon}{18}}P_{N}h\|_{L_{xt}^{\frac{6}{1+\epsilon}}}\nonumber\\
&&\leq C\sum N^{-\epsilon}\left(\prod\limits_{j=2}^{k+1}N_{j}^{\epsilon}\right)
\left(\prod\limits_{j=1}^{k}\|P_{N_{j}}v_{j}\|_{X_{s,b}}\right)
\|P_{N_{k+1}}z_{k+1}\|_{L_{t}^{\frac{1}{\epsilon}}L_{x}^{\frac{1}{2\epsilon}}}
\|P_{N}h\|_{X_{0,b_{1}}}\nonumber\\
&&\leq C\left[\prod\limits_{j=1}^{k}\|P_{N_{j}}v_{j}\|_{\ell_{N_{j}}^{2}X_{s,b}}\right]
\left[\|P_{N_{k+1}}z_{k+1}\|_{\ell_{N_{k+1}}^{\frac{1}{\epsilon}}L_{t}^{\frac{1}{\epsilon}}L_{x}^{\frac{1}{2\epsilon}}}\right]
\|P_{N}h\|_{\ell_{N}^{2}X_{0,b_{1}}}
\nonumber\\
&&\leq
C\left(\prod_{j=1}^{k}\|v_{j}\|_{X_{s,b}}\right)
\left[\|z_{k+1}\|_{L_{t}^{\frac{1}{\epsilon}}L_{x}^{\frac{1}{2\epsilon}}}\right]
\|h\|_{X_{0,b_{1}}}\leq
 C\lambda\left(\prod_{j=1}^{k}\|v_{j}\|_{X_{s,b}}\right)
\|h\|_{X_{0,b_{1}}}
\end{eqnarray*}
outside a set of probability
$
\leq C_{1}e^{-C\lambda^{2}\|f\|_{H^{s}}^{-2}}.
$

\noindent (3). Case $N_{1}\geq80k, N_{1}\sim N_{l}\geq 80kN_{l+1}(2\leq l\leq k-1), N_{1}\geq\frac{N_{k+1}}{80k}$.
This case can be proved similarly to Case (2).

\noindent (4). Case $N_{1}\geq80k,  N_{1}\sim N_{k},N_{1}\geq\frac{N_{k+1}}{80k}, N_{k+1}\leq1$.
 By using the H\"older inequality and  Lemmas 2.2, 2.7, 2.8, since $s\geq \frac{1}{2}-\frac{2}{k}+88\epsilon,$
we have
\begin{eqnarray*}
&&I_{2}\leq C\sum N^{\epsilon}N_{1}^{-2\epsilon}\|D_{x}^{\frac{3-\epsilon}{18}}P_{N}h\|_{L_{xt}^{\frac{6}{1+\epsilon}}}
\left[\prod\limits_{j=1}^{k}\|D_{x}^{\frac{5+6s}{6k}+2\epsilon}P_{N_{j}}v_{j}\|_{L_{xt}^{\frac{6k}{5-\epsilon}}}\right]
\|P_{N_{k+1}}z_{k}\|_{L_{xt}^{\infty}}
\nonumber\\&&\leq
C\sum N^{\epsilon} N_{1}^{-2\epsilon}N_{k+1}^{\frac{1}{2}}
\left(\prod_{j=1}^{k}\|P_{N_{j}}v_{j}\|_{X_{s,b}}\right)
\|P_{N_{k+1}}z_{k+1}\|_{X_{0,b}}
\|P_{N}h\|_{X_{0,b_{1}}}\nonumber\\
&&\leq C\sum N^{\epsilon}N_{1}^{-2\epsilon}N_{k+1}^{\frac{1}{2}}
\left(\prod_{j=1}^{k}\|P_{N_{j}}v_{j}\|_{X_{s,b}}\right)
\|P_{N_{k+1}}z_{k+1}\|_{X_{s_{1},b}}
\|P_{N}h\|_{X_{0,b_{1}}}\nonumber\\
&&\leq C\left[\prod_{j=1}^{k}\|P_{N_{j}}v_{j}\|_{\ell_{N_{j}}^{2}X_{s,b}}\right]
\|P_{N_{k+1}}z_{k+1}\|_{\ell_{N_{k+1}}^{2}X_{0,b}}
\|P_{N}h\|_{\ell_{N}^{2}X_{0,b_{1}}}\nonumber\\
&&\leq C\left[\prod_{j=1}^{k}
\|v_{j}\|_{X_{s,b}}\right]\|z_{k+1}\|_{X_{0,b}}
\|h\|_{X_{0,b_{1}}}\leq
C\lambda\left(\prod_{j=1}^{k}
\|v_{j}\|_{X_{s,b}}\right)
\end{eqnarray*}
outside a set of probability
$
\leq C_{1}e^{-C\lambda^{2}\|f\|_{H^{s}}^{-2}}.
$

\noindent (5). Case $N_{1}\geq80k, N_{1}\sim N_{k},N_{1}\geq \frac{N_{k+1}}{80k}, N_{k+1}\geq1$.
 By using the H\"older inequality,  (\ref{2.05}), (\ref{2.09}) and  Lemma 2.14,
 since $s\geq \frac{1}{2}-\frac{2}{k}+88\epsilon,$
we have
\begin{eqnarray*}
&&I_{2}\leq C\sum N^{\epsilon}N_{1}^{-2\epsilon}\|D_{x}^{\frac{3-\epsilon}{18}}P_{N}h\|_{L_{xt}^{\frac{6}{1+\epsilon}}}
\left[\prod\limits_{j=1}^{k}\|D_{x}^{\frac{5+6s}{6k}+3\epsilon}P_{N_{j}}v_{j}\|_{L_{xt}^{\frac{6k}{5-\epsilon}}}\right]
\|P_{N_{k+1}}z_{k}\|_{L_{xt}^{\infty}}
\nonumber\\&&\leq C\lambda
\sum N^{\epsilon}N_{1}^{-2\epsilon}
\left(\prod_{j=1}^{k}\|P_{N_{j}}v_{j}\|_{X_{s,b}}\right)
\|P_{N}h\|_{X_{0,b_{1}}}\nonumber\\
&&\leq C\lambda\sum N^{\epsilon}N_{1}^{-2\epsilon}
\left(\prod_{j=1}^{k}\|P_{N_{j}}v_{j}\|_{X_{s,b}}\right)
\|P_{N}h\|_{X_{0,b_{1}}}\nonumber\\
&&\leq C\lambda\left[\prod_{j=1}^{k}\|P_{N_{j}}v_{j}\|_{\ell_{N_{j}}^{2}X_{s,b}}\right]
\|P_{N}h\|_{\ell_{N}^{2}X_{0,b_{1}}}\leq C\lambda\left[\prod_{j=1}^{k}
\|v_{j}\|_{X_{s,b}}\right]
\|h\|_{X_{0,b_{1}}}
\end{eqnarray*}
outside a set of probability
$
\leq C_{1}e^{-C\lambda^{2}\|f\|_{H^{s}}^{-2}}.
$

\noindent (6). Case  $N_{k+1}\geq80k, \frac{N_{k+1}}{80k}\geq N_{1}\geq 80k N_{2},N_{2}\leq1$.
 By using the H\"older inequality and  Lemma 2.2,  2.8, since $s\geq \frac{1}{2}-\frac{2}{k}+88\epsilon,$
we have
\begin{eqnarray*}
&&\hspace{-0.8cm}I_{2}\leq C\sum N_{k+1}^{s_{4}}
\left\|I^{\frac{1}{2}}(J^{s_{1}}P_{N_{k+1}}z_{k+1},J^{s}P_{N_{1}}v_{1})\right\|_{L_{xt}^{2}}
\left\|I^{\frac{1}{2}-2\epsilon}(P_{N}h,J^{s}P_{N_{2}}v_{2})\right\|_{L_{xt}^{2}}
\prod_{j=3}^{k}\|P_{N_{j}}v_{j}\|_{L_{xt}^{\infty}}
\nonumber\\&&\hspace{-0.8cm}\leq
C\sum N_{k+1}^{s_{4}}\left(\prod_{j=3}^{k}N_{j}^{\frac{1}{2}}\right)
\left(\prod_{j=1}^{k}\|P_{N_{j}}v_{j}\|_{X_{s,b}}\right)
\|P_{N_{k+1}}z_{k+1}\|_{X_{s_{1},b}}
\|P_{N}h\|_{X_{0,b_{1}}}\nonumber\\
&&\hspace{-0.8cm}\leq C\lambda\sum N_{k+1}^{s_{4}}\left(\prod_{j=2}^{k}N_{j}^{\frac{k-2}{2(k-1)}}\right)
\left(\prod_{j=1}^{k}\|P_{N_{j}}v_{j}\|_{X_{s,b}}\right)
\|P_{N}h\|_{X_{0,b_{1}}}\nonumber\\
&&\hspace{-0.8cm}\leq C\lambda\left[\prod_{j=1}^{k}\|P_{N_{j}}v_{j}\|_{\ell_{N_{j}}^{2}X_{s,b}}\right]
\|P_{N}h\|_{\ell_{N}^{2}X_{0,b_{1}}}\leq
C\lambda\left(\prod_{j=1}^{k}
\|v_{j}\|_{X_{s,b}}\right)\|h\|_{X_{0,b_{1}}}
\end{eqnarray*}
outside a set of probability
$
\leq C_{1}e^{-C\lambda^{2}\|f\|_{H^{s}}^{-2}}.
$

\noindent (7). Case  $N_{k+1}\geq 80k, \frac{N_{k+1}}{80k}\geq N_{1}\geq 80k N_{2},N_{2}\geq1\geq N_{3}$.
 By using the H\"older inequality and  Lemma 2.2, 2.8, since $s\geq \frac{1}{2}-\frac{2}{k}+88\epsilon,$
we have
\begin{eqnarray*}
&&\hspace{-0.8cm}I_{2}\leq C\sum N_{k+1}^{s_{4}}
\left\|I^{\frac{1}{2}}(J^{s_{1}}P_{N_{k+1}}z_{k+1},J^{s}P_{N_{1}}v_{1})\right\|_{L_{xt}^{2}}
\left\|I^{\frac{1}{2}-2\epsilon}(P_{N}h,J^{s}P_{N_{2}}v_{2})\right\|_{L_{xt}^{2}}
\prod_{j=3}^{k}\|P_{N_{j}}v_{j}\|_{L_{xt}^{\infty}}
\nonumber\\&&\hspace{-0.8cm}\leq
C\sum N_{k+1}^{s_{4}}\left(\prod_{j=3}^{k}N_{j}^{\frac{1}{2}}\right)
\left(\prod_{j=1}^{k}\|P_{N_{j}}v_{j}\|_{X_{s,b}}\right)
\|P_{N_{k+1}}z_{k+1}\|_{X_{s_{1},b}}
\|P_{N}h\|_{X_{0,b_{1}}}\nonumber\\
&&\hspace{-0.8cm}\leq C\lambda\sum N_{k+1}^{s_{4}}\left(\prod_{j=3}^{k}N_{j}^{\frac{1}{2}}\right)
\left(\prod_{j=1}^{k}\|P_{N_{j}}v_{j}\|_{X_{s,b}}\right)
\|P_{N}h\|_{X_{0,b_{1}}}\nonumber\\
&&\hspace{-0.8cm}\leq C\lambda\left[\prod_{j=1}^{k}\|P_{N_{j}}v_{j}\|_{\ell_{N_{j}}^{2}X_{s,b}}\right]
\|P_{N}h\|_{\ell_{N}^{2}X_{0,b_{1}}}\leq C\lambda\left[\prod_{j=1}^{k}
\|v_{j}\|_{X_{s,b}}\right]
\|h\|_{X_{0,b_{1}}}
\end{eqnarray*}
outside a set of probability
$
\leq C_{1}e^{-C\lambda^{2}\|f\|_{H^{s}}^{-2}}.
$

\noindent (8). Case  $N_{k+1}\geq 80k, \frac{N_{k+1}}{80k}\geq N_{1}\geq80k N_{2},N_{1}\geq N_{2}\geq \cdot\cdot\cdot
\geq N_{l}\geq 1\geq\cdot\cdot\cdot\geq N_{k-1}\geq N_{k},N_{1}\leq N_{k+1}^{\frac{2}{3}}$.
 By using the H\"older inequality and  Lemmas 2.2, 2.8, since $s_{1}\geq{\rm max}\left\{\frac{s}{k+1},s-\frac{1}{3},\frac{s-\frac{1}{3}}{k}\right\}+88\epsilon$,
we have
\begin{eqnarray*}
&&I_{2}\leq C\sum N_{k+1}^{s_{4}}
\left\|I^{\frac{1}{2}}(J^{s_{1}}P_{N_{k+1}}z_{k+1},J^{s}P_{N_{1}}v_{1})\right\|_{L_{xt}^{2}}
\left\|I^{\frac{1}{2}-2\epsilon}(P_{N}h,J^{s}P_{N_{2}}v_{2})\right\|_{L_{xt}^{2}}\nonumber\\&&\qquad\qquad\qquad\times
\left(\prod_{j=3}^{l}\|D_{x}^{-\frac{2}{k}}P_{N_{j}}v_{j}\|_{L_{xt}^{\infty}}\right)
\left(\prod_{j=l+1}^{k}\|P_{N_{j}}v_{j}\|_{L_{xt}^{\infty}}\right)
\nonumber\\&&\leq
C\sum N_{k+1}^{s_{4}}\left(\prod_{j=l+1}^{k}N_{j}^{\frac{1}{2}}\right)
\left(\prod_{j=1}^{k}\|P_{N_{j}}v_{j}\|_{X_{s,b}}\right)
\|P_{N_{k+1}}z_{k+1}\|_{X_{s_{1},b}}
\|P_{N}h\|_{X_{0,b_{1}}}\nonumber\\
&&\leq C\lambda\sum N_{k+1}^{s_{4}}\left(\prod_{j=l+1}^{k}N_{j}^{\frac{1}{2}}\right)
\left(\prod_{j=1}^{k}\|P_{N_{j}}v_{j}\|_{X_{s,b}}\right)
\|P_{N}h\|_{X_{0,b_{1}}}\nonumber\\
&&\leq C\lambda\left[\prod_{j=1}^{k}\|P_{N_{j}}v_{j}\|_{\ell_{N_{j}}^{2}X_{s,b}}\right]
\|P_{N}h\|_{\ell_{N}^{2}X_{0,b_{1}}}\leq
C\lambda\left(\prod_{j=1}^{k}
\|v_{j}\|_{X_{s,b}}\right)\|h\|_{X_{0,b_{1}}}
\end{eqnarray*}
outside a set of probability
$
\leq C_{1}e^{-C\lambda^{2}\|f\|_{H^{s}}^{-2}}.
$

\noindent
(9). Case  $N_{k+1}\geq 80k,\frac{N_{k+1}}{80k}\geq N_{1}\geq80k N_{2},N_{1}
\geq N_{k+1}^{\frac{2}{3}},N_{2}\geq N_{k+1}^{\frac{2}{3}}$.
 By using the H\"older inequality,  Lemmas 2.2,  2.7, 2.14 and (\ref{2.011}) as well as (\ref{2.021}),
we have
\begin{eqnarray*}
&&\hspace{-0.8cm}I_{2}\leq C\sum N_{k+1}^{s_{4}}
\left\|I^{\frac{1}{2}}(J^{s}P_{N_{1}}v_{1},J^{s}P_{N_{2}}v_{2})\right\|_{L_{xt}^{2}}\left(\prod_{j=3}^{k}N_{j}^{\epsilon}\right)\left(\prod_{j=3}^{k}
\|D^{\frac{1-2\epsilon}{2(k-2)}-\frac{2}{k}-\epsilon}P_{N_{j}}v_{j}\|_{L_{x}^{\frac{2(k-2)}{(1-2\epsilon)}}L_{t}^{\infty}}\right)
\nonumber\\&&\qquad\qquad\qquad\times
\|D_{x}^{s_{1}}P_{N_{k+1}}z_{k+1}\|_{L_{xt}^{\infty}}
\|D_{x}^{1-2\epsilon}P_{N}h\|_{L_{x}^{\frac{1}{\epsilon}}L_{t}^{2}}\nonumber\\
&&\hspace{-0.8cm}\leq C\lambda\sum N_{k+1}^{s_{4}}\left(\prod_{j=3}^{k}N_{j}^{\epsilon}\right)\left[\prod_{j=1}^{k}\|P_{N_{j}}v_{j}\|_{X_{s,b}}\right]
\|h\|_{X_{0,b_{1}}}\nonumber\\
&&\hspace{-0.8cm}\leq
C\lambda\left[\prod_{j=1}^{k}\|P_{N_{j}}v_{j}\|_{\ell_{N_{j}}^{2}X_{s,b}}\right]
\|P_{N}h\|_{\ell_{N}^{2}X_{0,b_{1}}}\leq C\lambda\left(\prod_{j=1}^{k}\|v_{j}\|_{X_{s,b}}\right)
\|h\|_{X_{0,b_{1}}}
\end{eqnarray*}
outside a set of probability
$
\leq C_{1}e^{-C\lambda^{2}\|f\|_{H^{s}}^{-2}}.
$

\noindent (10). Case $N_{k+1}\geq 80k,\frac{N_{k+1}}{80k}\geq N_{1}\geq80k N_{2},
 N_{1}\geq N_{k+1}^{\frac{2}{3}},N_{2}\leq N_{k+1}^{\frac{2}{3}}$. This case
can be proved similarly to  Case (9).

\noindent (11). Case $N_{k+1}\geq 80k,\frac{N_{k+1}}{80k}\geq N_{1}\sim N_{l}\geq 80kN_{l+1}(2\leq l\leq k-1),
 N_{1}\leq N_{k+1}^{\frac{2}{3}}$. This case
can be proved similarly to  Case (8).

\noindent (12). Case $N_{k+1}\geq 80k,\frac{N_{k+1}}{80k}\geq N_{1}\sim N_{l}\geq 80kN_{l+1}(2\leq l\leq k-1),
 N_{1}\geq N_{k+1}^{\frac{2}{3}}$. This case
can be proved similarly to  Case (9).

\noindent (13). Case $N_{k+1}\geq 80k,\frac{N_{k+1}}{80k}\geq N_{1}\sim N_{k},
 N_{1}\leq N_{k+1}^{\frac{2}{3}}$. This case
can be proved similarly to  Case (10).

\noindent (14). Case $N_{k+1}\geq 80k,\frac{N_{k+1}}{80k}\geq N_{1}\sim N_{k},
 N_{1}\geq N_{k+1}^{\frac{2}{3}}$.
 By using the H\"older inequality,  Lemmas 2.2, 2.14 and (\ref{2.05}) as well as (\ref{2.035}),
we have
\begin{eqnarray*}
&&I_{2}\leq C\sum N_{k+1}^{s_{4}}
\left\|I^{\frac{1}{2}-2\epsilon}(J^{s}P_{N_{1}}v_{1},P_{N}h)\right\|_{L_{xt}^{2}}\left(\prod_{j=2}^{k}
\|J^{\frac{2}{k(k-1)}+2\epsilon}P_{N_{j}}v_{j}\|_{L_{xt}^{2(k-1)}}\right)\nonumber\\&&\qquad\qquad\qquad\qquad\times
\|D_{x}^{s_{1}}P_{N_{k+1}}z_{k+1}\|_{L_{xt}^{\infty}}\nonumber\\
&&\leq C\lambda\sum N_{k+1}^{s_{4}}\left[\prod_{j=1}^{k}\|P_{N_{j}}v_{j}\|_{X_{s,b}}\right]
\|P_{N}h\|_{X_{0,b_{1}}}\nonumber\\
&&\leq
C\lambda\left[\prod_{j=1}^{k}\|P_{N_{j}}v_{j}\|_{\ell_{N_{j}}^{2}X_{s,b}}\right]
\|P_{N}h\|_{\ell_{N}^{2}X_{0,b_{1}}}
\leq C\lambda\left(\prod_{j=1}^{k}\|v_{j}\|_{X_{s,b}}\right)
\|h\|_{X_{0,b_{1}}}
\end{eqnarray*}
outside a set of probability
$
\leq C_{1}e^{-C\lambda^{2}\|f\|_{H^{s}}^{-2}}.
$

We have completed the proof of Lemma 3.4.

\noindent {\bf Remark 9:}
Cases(8)-(10) of Lemma 3.4 requires  $s_{1}\geq s-\frac{1}{3}+88\epsilon.$

\begin{Lemma}\label{Lemma3.5}
Let $s\geq \frac{1}{2}-\frac{2}{k}+88\epsilon$ and $s_{1}\geq{\rm max}
\left\{\frac{s}{k+1},
s-\frac{1}{3},\frac{s-\frac{1}{3}}{k}\right\}+88\epsilon$,
$s_{2}=-1-s_{1}+s-8\epsilon(<0)$, $s_{3}=s-(k+1)s_{1}-8\epsilon$, $s_{4}=s-s_{1}-\frac{1}{3}-8\epsilon$
 and $z_{j}=\phi(t)U(t)f^{\omega}(k\leq j\leq k+1,j\in Z)$
 and   $b=\frac{1}{2}+\frac{\epsilon}{24},b_{1}=\frac{1}{2}-\frac{\epsilon}{12}$.
 We  define $\sum=\sum\limits_{N_{1}, \cdot\cdot\cdot,  N_{k+1},N}$.
  Then, we have
\begin{eqnarray}
\left\|\partial_{x}\left[\left(\prod_{j=1}^{k-1}v_{j}\right)\left(\prod\limits_{j=k}^{k+1}z_{j}
\right)\right]\right\|_{X_{s,-b_{1}}}\leq C
\lambda^{2}\left(\prod_{j=1}^{k-1}\|v_{j}\|_{X_{s,b}}\right)\label{3.015}
\end{eqnarray}
outside a set of probability
$
\leq C_{1}e^{-C\lambda^{2}\|f\|_{H^{s}}^{-2}}.
$

\end{Lemma}
\noindent {\bf Proof.} To prove (\ref{3.015}), by duality,  it suffices to prove
\begin{eqnarray}
\left|\int_{\SR^{2}}J^{s}\partial_{x}\left[\left(\prod_{j=1}^{k-1}v_{j}\right)
\left(\prod\limits_{j=k}^{k+1}z_{j}\right)\right]\bar{h}dxdt\right|\leq
 C\lambda^{2}\left(\prod_{j=1}^{k-1}\|v_{j}\|_{X_{s,b}}\right)
\|h\|_{X_{0,b_{1}}}\label{3.016}
\end{eqnarray}
outside a set of probability
$
\leq C_{1}e^{-C\lambda^{2}\|f\|_{H^{s}}^{-2}}.
$

We define
  \begin{eqnarray*}
  && I_{3}=\left|\int_{\SR^{2}}J^{s}\partial_{x}
  \left[\left(\prod_{j=1}^{k-1}v_{j}\right)\left(\prod\limits_{j=k}^{k+1}z_{j}
\right)\right]\bar{h}dxdt\right|\nonumber\\&&
=C\sum
\left|\int_{\SR^{2}}J^{s}\partial_{x}
\left[\left(\prod_{j=1}^{k-1}P_{N_{j}}v_{j}\right)\left(\prod\limits_{j=k}^{k+1}P_{N_{j}}z_{j}
\right)\right]P_{N}\bar{h}dxdt\right|.
   \end{eqnarray*}
We dyadically decompose $\mathscr{F}z_{j}$
with $(k\leq j\leq k+1,j\in Z)$ and $\mathscr{F}v_{j}$
with $(1\leq j\leq k-1,j\in Z)$ and $h$ such that frequency supports are
$\left\{|\xi_{j}|\sim N_{j}\right\}(1\leq j\leq k+1)$ and  $\left\{|\xi|\sim N\right\}$ for
 some dyadically $N_{j},N$.
Without loss of generality, we can assume that
 $N_{1}\geq N_{2}\geq\cdot\cdot\cdot\geq N_{k-1}$ and $N_{k}\geq  N_{k+1}$.

\noindent (1). Case ${\rm max}\left\{N_{1},N_{k}\right\}\leq 80k.$
By using the proof similar to Case (1) of Lemma 3.3,
we have
\begin{eqnarray*}
I_{3}
\leq C\lambda^{2}\left(\prod_{j=1}^{k-1}\|v_{j}\|_{X_{s,b}}\right)
\|h\|_{X_{0,b_{1}}}
\end{eqnarray*}
outside a set of probability
$
\leq C_{1}e^{-C\lambda^{2}\|f\|_{H^{s}}^{-2}}.
$

\noindent (2). Case $N_{1}\geq80 k{\rm max}\left\{1,N_{2},N_{k}\right\}.$
By using the proof similar to Case (2) of Lemma 3.3,
we have
\begin{eqnarray*}
I_{3}
\leq C\lambda^{2}\left(\prod_{j=1}^{k-1}\|v_{j}\|_{X_{s,b}}\right)
\|h\|_{X_{0,b_{1}}}
\end{eqnarray*}
outside a set of probability
$
\leq C_{1}e^{-C\lambda^{2}\|f\|_{H^{s}}^{-2}}.
$

\noindent (3). Case $N_{1}\geq80k{\rm max}\left\{N_{2},1\right\},
\frac{N_{k}}{80k}\leq N_{1}\leq80kN_{k}.$
By using the proof similar to Case (3) of Lemma 3.3,
we have
\begin{eqnarray*}
I_{3}
\leq C\lambda^{2}\left(\prod_{j=1}^{k-1}\|v_{j}\|_{X_{s,b}}\right)
\|h\|_{X_{0,b_{1}}}
\end{eqnarray*}
outside a set of probability
$
\leq C_{1}e^{-C\lambda^{2}\|f\|_{H^{s}}^{-2}}.
$

\noindent (4). Case $N_{1}\geq80k,N_{1}\sim N_{l}\geq 80kN_{l+1}(2\leq l\leq k-2),
\frac{N_{k}}{80k}\leq N_{1}\leq80kN_{k}.$
By using the proof similar to Case (2) of Lemma 3.3,
we have
\begin{eqnarray*}
I_{3}
\leq C\lambda^{2}\left(\prod_{j=1}^{k-1}\|v_{j}\|_{X_{s,b}}\right)
\|h\|_{X_{0,b_{1}}}
\end{eqnarray*}
outside a set of probability
$
\leq C_{1}e^{-C\lambda^{2}\|f\|_{H^{s}}^{-2}}.
$

\noindent (5). Case $N_{1}\geq80k,N_{1}\sim N_{k-1},
\frac{N_{k}}{80k}\leq N_{1}\leq80kN_{k}, N_{1}\geq80kN_{k+1}.$
By using the proof similar to Case (7) of Lemma 3.2,
we have
\begin{eqnarray*}
I_{3}
\leq C\lambda^{2}\left(\prod_{j=1}^{k-1}\|v_{j}\|_{X_{s,b}}\right)
\|h\|_{X_{0,b_{1}}}
\end{eqnarray*}
outside a set of probability
$
\leq C{\rm exp}\left(-c\frac{\lambda^{2}}{\|f\|_{H^{s_{1}}}^{2}}\right).
$

\noindent (6). Case $N_{k}\geq 80k{\rm max}\left\{N_{k+1},1\right\}$, $\frac{N_{k}}{80k}\geq N_{1}\geq80 N_{2}.$
By using the proof similar to Case (3) of Lemma 3.3,
we have
\begin{eqnarray*}
I_{3}
\leq C\lambda^{2}\left(\prod_{j=1}^{k-1}\|v_{j}\|_{X_{s,b}}\right)
\|h\|_{X_{0,b_{1}}}
\end{eqnarray*}
outside a set of probability
$
\leq C_{1}e^{-C\lambda^{2}\|f\|_{H^{s}}^{-2}}.
$

\noindent (7). Case $N_{k}\geq80k,\frac{N_{k}}{80k}\geq N_{1}$.
By using the proof similar to Cases (6)-(14) of Lemma 3.3,
we have
\begin{eqnarray*}
I_{3}
\leq C\lambda^{2}\left(\prod_{j=1}^{k-1}\|v_{j}\|_{X_{s,b}}\right)
\|h\|_{X_{0,b_{1}}}
\end{eqnarray*}
outside a set of probability
$
\leq C_{1}e^{-C\lambda^{2}\|f\|_{H^{s}}^{-2}}.
$

We have completed the proof of Lemma 3.5.

\begin{Lemma}\label{Lemma3.6}
Let $s\geq \frac{1}{2}-\frac{2}{k}+88\epsilon$ and $s_{1}\geq{\rm max}
\left\{\frac{s}{k+1},
s-\frac{1}{3},\frac{s-\frac{1}{3}}{k}\right\}+88\epsilon$,
$s_{2}=-1-s_{1}+s-8\epsilon(<0)$, $s_{3}=s-(k+1)s_{1}-8\epsilon$,
 $s_{4}=s-s_{1}-\frac{1}{3}-8\epsilon$ and $z_{j}=
\phi(t)U(t)f^{\omega}(2\leq j\leq k+1,j\in Z)$
 and   $b=\frac{1}{2}+\frac{\epsilon}{24},b_{1}=\frac{1}{2}-\frac{\epsilon}{12}$. We  define
  $\sum=\sum\limits_{N_{1}, \cdot\cdot\cdot,  N_{k+1},N}$.
  Then, we have
\begin{eqnarray}
\left\|\partial_{x}\left[v_{1}
\left(\prod\limits_{j=2}^{k+1}z_{j}
\right)\right]\right\|_{X_{s,-b_{1}}}\leq C
\lambda^{k}\|v_{1}\|_{X_{s,b}}\label{3.017}
\end{eqnarray}
outside a set of probability
$
\leq C_{1}e^{-C\lambda^{2}\|f\|_{H^{s}}^{-2}}.
$
\end{Lemma}
\noindent {\bf Proof.} To prove (\ref{3.017}), by duality,
it suffices to prove
\begin{eqnarray}
\left|\int_{\SR^{2}}J^{s}\partial_{x}\left[v_{1}
\left(\prod\limits_{j=2}^{k+1}z_{j}\right)\right]\bar{h}dxdt\right|\leq
 C\lambda^{k}\|v_{1}\|_{X_{s,b}}
\|h\|_{X_{0,b_{1}}}\label{3.018}
\end{eqnarray}
outside a set of probability
$
\leq C_{1}e^{-C\lambda^{2}\|f\|_{H^{s}}^{-2}}.
$

 We dfine
  \begin{eqnarray*}
  && I_{4}=\left|\int_{\SR^{2}}J^{s}\partial_{x}
   \left[v_{1}\left(\prod\limits_{j=2}^{k+1}z_{j}
\right)\right]\bar{h}dxdt\right|
=C\sum
\left|\int_{\SR^{2}}J^{s}\partial_{x}\left[P_{N_{1}}v_{1}
\left(\prod\limits_{j=2}^{k+1}z_{j}
\right)\right]P_{N}\bar{h}dxdt\right|.
   \end{eqnarray*}
We dyadically decompose $\mathscr{F}z_{j}$
with $(2\leq j\leq k+1,j\in Z)$ and $\mathscr{F}v_{1}$ and $h$ such that frequency supports are
$\left\{|\xi_{j}|\sim N_{j}\right\}(1\leq j\leq k+1)$ and  $\left\{|\xi|\sim N\right\}$ for
 some dyadically $N_{j},N$.
Without loss of generality, we can assume that
 $N_{1}$ and $N_{2}\geq
  \cdot\cdot\cdot\geq N_{k}\geq  N_{k+1}$.

\noindent (1). Case ${\rm max}\left\{N_{1},N_{2}\right\}\leq 80k.$
This case can be proved similarly to Case (1) of Lemma 3.2.

\noindent (2). Case $N_{1}\geq80 k{\rm max}\left\{1,N_{2}\right\}.$
By using the proof similar to Case (2) of Lemma 3.2,
we have
\begin{eqnarray*}
&&I_{4}
\leq C\lambda^{k}\|v_{1}\|_{X_{s,b}}
\|h\|_{X_{0,b_{1}}}
\end{eqnarray*}
outside a set of probability
$
\leq C_{1}e^{-C\lambda^{2}\|f\|_{H^{s}}^{-2}}.
$

\noindent (3). Case $N_{2}\geq80k{\rm max}\left\{N_{3},1\right\},\frac{N_{2}}{80k}\leq N_{1}\leq80kN_{2}.$
By using  Lemma 2.2, (\ref{2.05}), (\ref{2.06}) and (\ref{2.08}) as well as (\ref{2.035}),
we have
\begin{eqnarray*}
&&I_{4}
\leq C\sum N_{1}^{-\frac{1}{6}-s_{1}}\|D_{x}^{s+\frac{1}{6}}P_{N_{1}}v_{1}\|_{L_{xt}^{6}}
\|I^{\frac{1}{2}}(J^{s_{1}}P_{N_{2}}z_{2},
J^{s_{1}}P_{N_{3}}z_{3})\|_{L_{xt}^{2}}
\left[\prod\limits_{j=4}^{k+1}\|P_{N_{j}}z_{j}\|_{L_{xt}^{\frac{24(k-2)}{5-3\epsilon}}}\right]\nonumber\\
&&\qquad\qquad\qquad\qquad\qquad\times\|P_{N}h\|_{L_{xt}^{\frac{8}{1+\epsilon}}}\nonumber\\
&&\leq C\sum N_{1}^{-\frac{1}{6}-s_{1}}\|P_{N_{1}}v_{1}\|_{X_{s_{1},b}}
\left(\prod\limits_{j=2}^{3}\|P_{N_{j}}z_{j}\|_{X_{s_{1},b}}\right)
\prod\limits_{j=3}^{k+1}N_{j}^{\epsilon}
\nonumber\\
&&\qquad\qquad\qquad\qquad\qquad
\times\left[\prod\limits_{j=4}^{k+1}\|P_{N_{j}}z_{j}\|_{L_{t}^{\frac{24(k-2)}{5-3\epsilon}}
L_{x}^{\frac{24(k-2)}{5+(24k-27)\epsilon}}}\right]
\|P_{N}h\|_{X_{0,b_{1}}}\nonumber\\
&&\leq C\lambda^{2}\sum N_{1}^{-\frac{1}{6}-s_{1}}\|v_{1}\|_{X_{s_{1},b}}
\prod\limits_{j=2}^{k+1}N_{j}^{\epsilon}
\left[\prod\limits_{j=4}^{k+1}\|P_{N_{j}}z_{j}\|_{L_{t}^{\frac{24(k-2)}{5-3\epsilon}}
L_{x}^{\frac{24(k-2)}{5+(24k-27)\epsilon}}}\right]\nonumber\\
&&\qquad\qquad\qquad\qquad\qquad\times
\|P_{N}h\|_{X_{0,b_{1}}}\nonumber\\
&&\leq C\|P_{N_{1}}v_{1}\|_{\ell_{N_{1}}^{2}X_{s_{1},b}}
\left(\prod\limits_{j=2}^{3}\|P_{N_{j}}z_{j}\|_{\ell_{N_{j}}^{2}X_{s_{1},b}}\right)
\left[\prod\limits_{j=4}^{k+1}\|P_{N_{j}}z_{j}\|_{\ell_{N_{j}}^{\frac{24(k-2)}{5-3\epsilon}}L_{t}^{\frac{24(k-2)}{5-3\epsilon}}
L_{x}^{\frac{24(k-2)}{5+(24k-27)\epsilon}}}\right]\nonumber\\
&&\qquad\qquad\qquad\qquad\qquad\times
\|P_{N}h\|_{\ell_{N}^{2}X_{0,b_{1}}}\nonumber\\
&&\leq C\|v_{1}\|_{X_{s_{1},b}}\left(\prod\limits_{j=2}^{3}\|z_{j}\|_{X_{s_{1},b}}\right)
\left[\prod\limits_{j=4}^{k+1}\|z_{j}\|_{L_{t}^{\frac{24(k-2)}{5-3\epsilon}}L_{x}^{\frac{24(k-2)}{5+(24k-27)\epsilon}}}\right]
\|h\|_{X_{0,b_{1}}}\nonumber\\
&&\leq C\lambda^{k}\|v_{1}\|_{X_{s,b}}
\|h\|_{X_{0,b_{1}}}
\end{eqnarray*}
outside a set of probability
$
\leq C_{1}e^{-C\lambda^{2}\|f\|_{H^{s}}^{-2}}.
$

\noindent (4). Case $N_{2}\geq80k,N_{2}\sim N_{l}\geq 80kN_{l+1}(3\leq l\leq k),\frac{N_{2}}{80k}\leq N_{1}\leq80kN_{2}.$
By using a  proof   similarly to Case (3) of Lemma 3.5,  we have
\begin{eqnarray*}
I_{4}
\leq C\lambda^{k}\|v_{1}\|_{X_{s,b}}
\|h\|_{X_{0,b_{1}}}
\end{eqnarray*}
outside a set of probability
$
\leq C_{1}e^{-C\lambda^{2}\|f\|_{H^{s}}^{-2}}.
$

\noindent (5). Case $N_{2}\geq80k,N_{2}\sim N_{k+1},\frac{N_{2}}{80k}\leq N_{1}\leq80kN_{2}.$
By using a  proof   similarly to Case (11) of Lemma 3.2,  we have
\begin{eqnarray*}
I_{4}
\leq C\lambda^{k}\|v_{1}\|_{X_{s,b}}
\|h\|_{X_{0,b_{1}}}
\end{eqnarray*}
outside a set of probability
$
\leq C_{1}e^{-C\lambda^{2}\|f\|_{H^{s}}^{-2}}.
$

\noindent (6). Case $N_{2}\geq80k{\rm max}\left\{N_{3},1\right\},\frac{N_{2}}{80k}\geq N_{1}\geq1.$
By using  Lemmas 2.2, 2.12, 2.14 and (\ref{2.035}),
we have
\begin{eqnarray*}
&&I_{4}
\leq C\sum N_{2}^{s_{2}}\left\|I^{\frac{1}{2}}(J^{s}P_{N_{1}}v_{1},J^{s_{1}}P_{N_{2}}z_{2})\right\|_{L_{xt}^{2}}
\left\|I^{\frac{1}{2}-2\epsilon}(J^{s}P_{N_{3}}z_{3},P_{N}h)\right\|_{L_{xt}^{2}}
\left(\prod_{j=4}^{k+1}\|P_{N_{j}}z_{j}\|_{L_{xt}^{\infty}}\right)
\nonumber\\
&&\leq C\sum N_{2}^{s_{2}}\|P_{N_{1}}v_{1}\|_{X_{s,b}}\left[\prod\limits_{j=2}^{3}\|P_{N_{j}}z_{j}\|_{X_{s,b}}\right]
\left(\prod_{j=4}^{k+1}\|P_{N_{j}}z_{j}\|_{L_{tx}^{\infty}}\right)\nonumber\\
&&\leq C\lambda^{k}\sum N_{2}^{s_{2}}\|P_{N_{1}}v_{1}\|_{X_{s,b}}\left(\prod_{j=3}^{k+1}N_{j}^{\frac{2(k-2)}{k-1}\epsilon}\right)
\|P_{N}h\|_{X_{0,b_{1}}}
\nonumber\\&&\leq C\lambda^{k}\|P_{N_{1}}v_{1}\|_{\ell_{N_{1}}^{2}X_{s,b}}\|P_{N}h\|_{\ell_{N}^{2}X_{0,b_{1}}}
\leq C\lambda^{k}\|v_{1}\|_{X_{s,b}}
\|h\|_{X_{0,b_{1}}}
\end{eqnarray*}
outside a set of probability
$
\leq C_{1}e^{-C\lambda^{2}\|f\|_{H^{s}}^{-2}}.
$

\noindent (7). Case $\frac{N_{2}}{80k}\geq N_{1}\geq1,N_{2}\sim N_{l}(3\leq l\leq k)\geq 80k N_{l+1}.$
By using  a   proof  similarly to Case (6), we have
\begin{eqnarray*}
&&I_{4}
\leq C\lambda^{k}\|v_{1}\|_{X_{s,b}}
\|h\|_{X_{0,b_{1}}}
\end{eqnarray*}
outside a set of probability
$
\leq C_{1}e^{-C\lambda^{2}\|f\|_{H^{s}}^{-2}}.
$

\noindent (8). Case $\frac{N_{2}}{80k}\geq N_{1}\geq1,N_{2}\sim N_{k+1}.$
By using  Lemmas 2.2, 2.13 and (\ref{2.035}),
we have
\begin{eqnarray*}
&&I_{4}
\leq C\sum N_{2}^{s-ks_{1}-\frac{1}{3}}N^{\epsilon}
\left\|I^{\frac{1}{2}}(J^{s}P_{N_{1}}v_{1},J^{s_{1}}P_{N_{2}}z_{2})\right\|_{L_{xt}^{2}}
\|J^{s_{1}+\frac{1}{6}}P_{N_{3}}z_{3}\|_{L_{xt}^{6}}
\nonumber\\&&\qquad\qquad\qquad \times \left(\prod_{j=4}^{k+1}
\|J^{s_{1}}P_{N_{j}}z_{j}\|_{L_{xt}^{\frac{6(k-2)}{1-\epsilon}}}\right)
\left\|D_{x}^{\frac{3-\epsilon}{18}}P_{N}h\right\|_{L_{xt}^{\frac{6}{1+\epsilon}}}\nonumber\\
&&\leq C\sum N_{2}^{s-ks_{1}-\frac{1}{3}}N^{\epsilon}\|P_{N_{1}}v_{1}\|_{X_{s,b}}
\left[\prod\limits_{j=2}^{3}\|P_{N_{j}}z_{j}\|_{X_{s_{1},b}}\right]
\left(\prod_{j=4}^{k+1}\|P_{N_{j}}z_{j}\|_{L_{xt}^{\frac{6(k-2)}{1-\epsilon}}}\right)
\|P_{N}h\|_{X_{0,b_{1}}}\nonumber\\
&&\leq C\lambda^{k}\sum N_{2}^{s-ks_{1}-\frac{1}{3}}N^{\epsilon}\|P_{N_{1}}v_{1}\|_{X_{s,b}}
\|P_{N}h\|_{X_{0,b_{1}}}
\nonumber\\&&\leq C\lambda^{k}\|P_{N_{1}}v_{1}\|_{\ell_{N_{1}}^{2}X_{s,b}}
\|P_{N}h\|_{\ell_{N}^{2}X_{0,b_{1}}}
\leq C\lambda^{k}\|v_{1}\|_{X_{s,b}}
\|h\|_{X_{0,b_{1}}}
\end{eqnarray*}
outside a set of probability
$
\leq C_{1}e^{-C\lambda^{2}\|f\|_{H^{s}}^{-2}}.
$

\noindent (9). Case $N_{2}\geq80k{\rm max}\left\{N_{3},1\right\}, \frac{N_{2}}{80k}\geq N_{1},N_{1}\leq1.$
By using  Lemmas 2.2, 2.15 and (\ref{2.035}),
we have
\begin{eqnarray*}
&&I_{4}
\leq C\sum N_{2}^{s_{2}}\|P_{N_{1}}v\|_{L_{xt}^{\infty}}
\left\|I^{\frac{1}{2}}(J^{s_{1}}P_{N_{2}}z_{2},J^{s_{1}}P_{N_{3}}z_{3})\right\|_{L_{xt}^{2}}
\nonumber\\&&\qquad\qquad\qquad\times
\left\|I^{\frac{1}{2}-2\epsilon}(J^{s}P_{N_{4}}z_{4},P_{N}h)\right\|_{L_{xt}^{2}}
\left(\prod_{j=5}^{k+1}\|P_{N_{j}}z_{j}\|_{L_{xt}^{\infty}}\right)
\nonumber\\
&&\leq C\sum N_{2}^{s_{2}}N_{1}^{\frac{1}{2}}\|P_{N_{1}}v_{1}\|_{X_{0,b}}
\left[\prod\limits_{j=2}^{4}\|P_{N_{j}}z_{j}\|_{X_{s,b}}\right]
\left(\prod_{j=5}^{k+1}N_{j}^{2\epsilon}\right)
\left(\prod_{j=5}^{k+1}\|P_{N_{j}}z_{j}\|_{L_{t}^{\infty}L_{x}^{\frac{1}{\epsilon}}}\right)\nonumber\\
&&\leq C\lambda^{k}\sum N_{2}^{s_{2}}\|P_{N_{1}}v_{1}\|_{X_{s,b}}
\left(\prod_{j=3}^{k+1}N_{j}^{\frac{2(k-2)}{k-1}\epsilon}\right)
\|P_{N}h\|_{X_{0,b_{1}}}
\nonumber\\&&\leq C\lambda^{k}\|P_{N_{1}}v_{1}\|_{\ell_{N_{1}}^{2}X_{s,b}}
\|P_{N}h\|_{\ell_{N}^{2}X_{0,b_{1}}}
\leq C\lambda^{k}\|v_{1}\|_{X_{s,b}}
\|h\|_{X_{0,b_{1}}}
\end{eqnarray*}
outside a set of probability
$
\leq C_{1}e^{-C\lambda^{2}\|f\|_{H^{s}}^{-2}}.
$

\noindent (10). Case $N_{2}\geq80k,\frac{N_{2}}{80k}\geq N_{1},N_{1}\leq1,N_{2}\sim N_{l}(3\leq l\leq k)
\geq 80k N_{l+1}.$
By using  a   proof  similarly to Case (9), we have
\begin{eqnarray*}
&&I_{4}
\leq C\lambda^{k}\|v_{1}\|_{X_{s,b}}
\|h\|_{X_{0,b_{1}}}
\end{eqnarray*}
outside a set of probability
$
\leq C_{1}e^{-C\lambda^{2}\|f\|_{H^{s}}^{-2}}.
$

\noindent (11). Case $N_{2}\geq80k,\frac{N_{2}}{80k}\geq N_{1},N_{1}\leq1,N_{2}\sim N_{k+1},N_{1}N_{2}\geq1.$
By using  a   proof  similarly to Case (7) of Lemma 3.2, we have
\begin{eqnarray*}
&&I_{4}
\leq C\lambda^{k}\|v_{1}\|_{X_{s,b}}
\|h\|_{X_{0,b_{1}}}
\end{eqnarray*}
outside a set of probability
$
\leq C_{1}e^{-C\lambda^{2}\|f\|_{H^{s}}^{-2}}.
$

\noindent (12). Case $N_{2}\geq80k,\frac{N_{2}}{80k}\geq N_{1},N_{1}\leq1,N_{2}\sim N_{k+1},N_{1}N_{2}\leq1.$
By using  a   proof  similarly to Case (8) of Lemma 3.2, we have
\begin{eqnarray*}
&&I_{4}
\leq C\lambda^{k}\|v_{1}\|_{X_{s,b}}
\|h\|_{X_{0,b_{1}}}
\end{eqnarray*}
outside a set of probability
$
\leq C_{1}e^{-C\lambda^{2}\|f\|_{H^{s}}^{-2}}.
$

We have completed the proof of Lemma 3.6.

\noindent {\bf Remark 10:} Case (8) of Lemma 3.5 requires  $s_{1}\geq \frac{s-\frac{1}{3}}{k}+88\epsilon.$

\begin{Lemma}\label{Lemma3.7}
Let $s\geq \frac{1}{2}-\frac{2}{k}+88\epsilon$ and $s_{1}\geq{\rm max}
\left\{\frac{s}{k+1},
s-\frac{1}{3},\frac{s-\frac{1}{3}}{k}\right\}+88\epsilon$,
$s_{2}=-1-s_{1}+s-8\epsilon(<0)$, $s_{3}=s-(k+1)s_{1}-8\epsilon$,
 $s_{4}=s-s_{1}-\frac{1}{3}-8\epsilon$ and $z_{j}=
\phi(t)U(t)f^{\omega}(2\leq j\leq k+1,j\in Z)$
 and   $b=\frac{1}{2}+\frac{\epsilon}{24},b_{1}=\frac{1}{2}-\frac{\epsilon}{12}$. We  define
  $\sum=\sum\limits_{N_{1}, \cdot\cdot\cdot,  N_{k+1},N}$.
  Then, we have
\begin{eqnarray}
\left\|\partial_{x}\left[v_{1}v_{2}
\left(\prod\limits_{j=3}^{k+1}z_{j}
\right)\right]\right\|_{X_{s,-b_{1}}}\leq C
\lambda^{k-1}\prod\limits_{j=1}^{2}\|v_{j}\|_{X_{s,b}}\label{3.019}
\end{eqnarray}
outside a set of probability
$
\leq C_{1}e^{-C\lambda^{2}\|f\|_{H^{s}}^{-2}}.
$

\end{Lemma}
\noindent {\bf Proof.} To prove (\ref{3.019}), by duality,
it suffices to prove
\begin{eqnarray}
\left|\int_{\SR^{2}}J^{s}\partial_{x}\left[v_{1}v_{2}
\left(\prod\limits_{j=3}^{k+1}z_{j}\right)\right]\bar{h}dxdt\right|\leq
 C\lambda^{k-1}\|h\|_{X_{0,b_{1}}}\prod\limits_{j=1}^{2}\|v_{j}\|_{X_{s,b}}
\label{3.020}
\end{eqnarray}
outside a set of probability
$
\leq C_{1}e^{-C\lambda^{2}\|f\|_{H^{s}}^{-2}}.
$

 We define
  \begin{eqnarray*}
  && I_{5}=\left|\int_{\SR^{2}}J^{s}\partial_{x}
   \left[v_{1}v_{2}\left(\prod\limits_{j=3}^{k+1}z_{j}
\right)\right]\bar{h}dxdt\right|\nonumber\\&&
=C\sum
\left|\int_{\SR^{2}}J^{s}\partial_{x}\left[P_{N_{1}}v_{1}P_{N_{2}}v_{2}
\left(\prod\limits_{j=3}^{k+1}z_{j}
\right)\right]P_{N}\bar{h}dxdt\right|.
   \end{eqnarray*}
We dyadically decompose $\mathscr{F}z_{j}$
with $(3\leq j\leq k+1,j\in Z)$ and $\mathscr{F}v_{j}(j=1,2)$ and $h$ such that frequency supports are
$\left\{|\xi_{j}|\sim N_{j}\right\}(1\leq j\leq k+1)$ and  $\left\{|\xi|\sim N\right\}$ for
 some dyadically $N_{j},N$.
Without loss of generality, we can assume that
 $N_{1}\geq N_{2}$ and $N_{3}\geq
  \cdot\cdot\cdot\geq N_{k}\geq  N_{k+1}$.

\noindent (1). Case ${\rm max}\left\{N_{1},N_{3}\right\}\leq 80k.$
By using the proof similar to Case (1) of Lemma 3.2,
we have
\begin{eqnarray*}
&&I_{5}
\leq C
 \lambda^{k-1}\|h\|_{X_{0,b_{1}}}\prod\limits_{j=1}^{2}\|v_{j}\|_{X_{s,b}}
\end{eqnarray*}
outside a set of probability
$
\leq C_{1}e^{-C\lambda^{2}\|f\|_{H^{s}}^{-2}}.
$

\noindent (2). Case $N_{1}\geq80 k{\rm max}\left\{1,N_{3}\right\}.$
By using the proof similar to Case (2) of Lemma 3.2,
we have
\begin{eqnarray*}
&&I_{5}
\leq C
 \lambda^{k-1}\|h\|_{X_{0,b_{1}}}\prod\limits_{j=1}^{2}\|v_{j}\|_{X_{s,b}}
\end{eqnarray*}
outside a set of probability
$
\leq C_{1}e^{-C\lambda^{2}\|f\|_{H^{s}}^{-2}}.
$

\noindent (3). Case $N_{1}\geq80k{\rm max}\left\{N_{2},1\right\},\frac{N_{3}}{80k}\leq N_{1}\leq80kN_{3}.$
By using  Lemmas 2.2, 2.13 and (\ref{2.08}) as well as (\ref{2.035}),
we have
\begin{eqnarray*}
&&I_{5}
\leq C\sum N_{1}^{-s_{1}}
\|I^{\frac{1}{2}}(J^{s}P_{N_{1}}v_{1},J^{s}P_{N_{2}}v_{2})\|_{L_{xt}^{2}}
\left[\prod\limits_{j=3}^{k+1}\|J^{s_{1}}P_{N_{j}}z_{j}\|_{L_{xt}^{\frac{8(k-1)}{3-\epsilon}}}\right]
\|P_{N}h\|_{L_{xt}^{\frac{8}{1+\epsilon}}}\nonumber\\
&&\leq C\sum N_{1}^{-s_{1}}\left[\prod\limits_{j=1}^{2}\|P_{N_{j}}v_{j}\|_{X_{s,b}}\right]
\left(\prod\limits_{j=3}^{k+1}\|J^{s_{1}}P_{N_{j}}z_{j}\|_{L_{xt}^{\frac{8(k-1)}{3+(2k-1)\epsilon}}}\right)
\left(\prod\limits_{j=2}^{k+1}N_{j}^{\epsilon}\right)
\|P_{N}h\|_{X_{0,b_{1}}}\nonumber\\
&&\leq C\sum N_{1}^{-s_{1}}\prod\limits_{j=1}^{2}\|v_{j}\|_{X_{s,b}}
\left(\prod\limits_{j=3}^{k+1}\|J^{s_{1}}P_{N_{j}}z_{j}\|_{L_{xt}^{\frac{8(k-1)}{3+(2k-1)\epsilon}}}\right)
\left[\prod\limits_{j=2}^{k+1}N_{j}^{\epsilon}\right]
\|P_{N}h\|_{X_{0,b_{1}}}\nonumber\\
&&\leq C\prod\limits_{j=1}^{2}\|P_{N_{j}}v_{j}\|_{\ell_{N_{j}}^{2}X_{s_{1},b}}
\left[\prod\limits_{j=3}^{k+1}\|J^{s_{1}}P_{N_{j}}z_{j}\|_{\ell_{N_{j}}^{\frac{8(k-1)}
{3-\epsilon}}L_{xt}^{\frac{8(k-1)}{3+(2k-1)\epsilon}}}\right]
\|P_{N}h\|_{\ell_{N}^{2}X_{0,b_{1}}}\nonumber\\
&&\leq C\left[\prod\limits_{j=1}^{2}\|v_{j}\|_{X_{s,b}}\right]
\left[\prod\limits_{j=3}^{k+1}\|J^{s_{1}}z_{j}\|_{L_{xt}^{\frac{8(k-1)}{3+(2k-1)\epsilon}}}\right]
\|h\|_{X_{0,b_{1}}}\nonumber\\
&&\leq C\lambda^{k-1}\|h\|_{X_{0,b_{1}}}\prod\limits_{j=1}^{2}\|v_{j}\|_{X_{s,b}}
\end{eqnarray*}
outside a set of probability
$
\leq C_{1}e^{-C\lambda^{2}\|f\|_{H^{s}}^{-2}}.
$

\noindent (4). Case $N_{1}\geq80k,N_{1}\sim N_{2},N_{3}\geq80k{\rm max}\left\{N_{4},1\right\},
\frac{N_{3}}{80k}\leq N_{1}\leq80kN_{3}.$
By using a  proof  similarly to Case (3) of Lemma 3.5, we have
\begin{eqnarray*}
&&I_{5}
\leq  C\lambda^{k-1}\|h\|_{X_{0,b_{1}}}\prod\limits_{j=1}^{2}\|v_{j}\|_{X_{s,b}}
\end{eqnarray*}
outside a set of probability
$
\leq C_{1}e^{-C\lambda^{2}\|f\|_{H^{s}}^{-2}}.
$

\noindent (5). Case $N_{1}\geq80k,N_{1}\sim N_{2},\frac{N_{3}}{80k}\leq N_{1}\leq80kN_{3},
N_{3}\sim  N_{m}\geq 80kN_{m+1}(4\leq m\leq k).$
By using a proof  similarly to Case (3) of Lemma 3.3,  we have
\begin{eqnarray*}
&&I_{5}
\leq  C\lambda^{k-1}\|h\|_{X_{0,b_{1}}}\prod\limits_{j=1}^{2}\|v_{j}\|_{X_{s,b}}
\end{eqnarray*}
outside a set of probability
$
\leq C_{1}e^{-C\lambda^{2}\|f\|_{H^{s}}^{-2}}.
$

\noindent (6). Case $N_{1}\geq80k,N_{1}\sim N_{2},\frac{N_{3}}{80k}\leq N_{1}\leq80kN_{3},N_{3}\sim N_{k+1}.$
By using a proof similarly   to Case (10) of Lemma 3.2, we have
\begin{eqnarray*}
&&I_{5}
\leq  C\lambda^{k-1}\|h\|_{X_{0,b_{1}}}\prod\limits_{j=1}^{2}\|v_{j}\|_{X_{s,b}}
\end{eqnarray*}
outside a set of probability
$
\leq C_{1}e^{-C\lambda^{2}\|f\|_{H^{s}}^{-2}}.
$

\noindent (7). Case $N_{1}\geq80k,\frac{N_{3}}{80k}\geq N_{1}, N_{3}\geq80k.$
By using a proof similarly   to Cases (6)-(9) of Lemma 3.3, we have
\begin{eqnarray*}
&&I_{5}
\leq  C\lambda^{k-1}\|h\|_{X_{0,b_{1}}}\prod\limits_{j=1}^{2}\|v_{j}\|_{X_{s,b}}
\end{eqnarray*}
outside a set of probability
$
\leq C_{1}e^{-C\lambda^{2}\|f\|_{H^{s}}^{-2}}.
$

We have completed the proof of Lemma 3.7.

\begin{Lemma}\label{Lemma3.8}
Let $3\leq l\leq k-2$,  $s\geq \frac{1}{2}-\frac{2}{k}+88\epsilon$,  $s_{1}\geq{\rm max}
\left\{\frac{s}{k+1},
s-\frac{1}{3},\frac{s-\frac{1}{3}}{k}\right\}+88\epsilon$,
$s_{2}=-1-s_{1}+s-8\epsilon(<0)$, $s_{3}=s-(k+1)s_{1}-8\epsilon$,
 $s_{4}=s-s_{1}-\frac{1}{3}-8\epsilon$ and $z_{j}=
\phi(t)U(t)f^{\omega}(l+1\leq j\leq k+1,j\in Z)$
 and   $b=\frac{1}{2}+\frac{\epsilon}{24},,b_{1}=\frac{1}{2}-\frac{\epsilon}{12}$. We  define
  $\sum=\sum\limits_{N_{1}, \cdot\cdot\cdot,  N_{k+1},N}$.
  Then, we have
\begin{eqnarray}
\left\|\partial_{x}\left[\prod\limits_{j=1}^{l}v_{j}
\left(\prod\limits_{j=l+1}^{k+1}z_{j}
\right)\right]\right\|_{X_{s,-b_{1}}}
\leq C\lambda^{k-l+1}\prod\limits_{j=1}^{l}\|v_{j}\|_{X_{s,b}}\label{3.021}
\end{eqnarray}
outside a set of probability
$
\leq C_{1}e^{-C\lambda^{2}\|f\|_{H^{s}}^{-2}}.
$
\end{Lemma}
\noindent {\bf Proof.} To prove (\ref{3.021}), by duality,
it suffices to prove
\begin{eqnarray}
\left|\int_{\SR^{2}}J^{s}\partial_{x}\left[\prod\limits_{j=1}^{l}v_{j}
\left(\prod\limits_{j=l+1}^{k+1}z_{j}\right)\right]\bar{h}dxdt\right|\leq
 C\lambda^{k-l+1}\|h\|_{X_{0,b_{1}}}\prod\limits_{j=1}^{l}\|v_{j}\|_{X_{s,b}}\label{3.022}
\end{eqnarray}
outside a set of probability
$
\leq \epsilon.
$
 We define
  \begin{eqnarray*}
  && I_{6}=\left|\int_{\SR^{2}}J^{s}\partial_{x}
   \left[\prod\limits_{j=1}^{l}v_{j}\left(\prod\limits_{j=l+1}^{k+1}z_{j}
\right)\right]\bar{h}dxdt\right|
\nonumber\\&&=C\sum
\left|\int_{\SR^{2}}J^{s}\partial_{x}\left[\prod\limits_{j=1}^{l}P_{N_{j}}v_{j}
\left(\prod\limits_{j=l+1}^{k+1}z_{j}
\right)\right]P_{N}\bar{h}dxdt\right|.
   \end{eqnarray*}
We dyadically decompose $\mathscr{F}z_{j}$
with $(l+1\leq j\leq k+1,j\in Z)$ and $\mathscr{F}v_{j}(1\leq j\leq l)$ and $h$ such that frequency supports are
$\left\{|\xi_{j}|\sim N_{j}\right\}(1\leq j\leq k+1)$ and  $\left\{|\xi|\sim N\right\}$ for
 some dyadically $N_{j},N$.
Without loss of generality, we can assume that
 $N_{1}\geq N_{2}\geq\cdot\cdot\cdot\geq N_{l}$ and $N_{l+1}\geq
  \cdot\cdot\cdot\geq N_{k}\geq  N_{k+1}$.

\noindent (1). Case ${\rm max}\left\{N_{1},N_{l+1}\right\}\leq 80k.$
By using a proof  similarly to Case (1) of Lemma 3.2,  we have
\begin{eqnarray*}
&&I_{6}
\leq C\lambda^{k-l+1}\|h\|_{X_{0,b_{1}}}\prod\limits_{j=1}^{l}\|v_{j}\|_{X_{s,b}}
\end{eqnarray*}
outside a set of probability
$
\leq C_{1}e^{-C\lambda^{2}\|f\|_{H^{s}}^{-2}}.
$

\noindent (2). Case $N_{1}\geq80 k{\rm max}\left\{1,N_{l+1}\right\}.$
By using the proof similar to Case (2) of Lemma 3.2,
we have
\begin{eqnarray*}
&&I_{6}
\leq C\lambda^{k-l+1}\|h\|_{X_{0,b_{1}}}\prod\limits_{j=1}^{l}\|v_{j}\|_{X_{s,b}}
\end{eqnarray*}
outside a set of probability
$
\leq C_{1}e^{-C\lambda^{2}\|f\|_{H^{s}}^{-2}}.
$

\noindent (3). Case $N_{1}\geq80k{\rm max}\left\{N_{2},1\right\}, N_{1}\geq\frac{N_{l+1}}{80k}.$
By using a proof  similarly to Case (2) of Lemma 3.3, we have
\begin{eqnarray*}
&&I_{6}
\leq C\lambda^{k-l+1}\|h\|_{X_{0,b_{1}}}\prod\limits_{j=1}^{l}\|v_{j}\|_{X_{s,b}}
\end{eqnarray*}
outside a set of probability
$
\leq C_{1}e^{-C\lambda^{2}\|f\|_{H^{s}}^{-2}}.
$

\noindent (4). Case $N_{1}\geq80k,N_{1}\sim N_{m}\geq 80kN_{m+1}(2\leq m\leq l-1),N_{1}\geq\frac{N_{l+1}}{80k}.$
By using a proof  similarly to Case (2) of Lemma 3.3, we have
\begin{eqnarray*}
&&I_{6}
\leq C\lambda^{k-l+1}\|h\|_{X_{0,b_{1}}}\prod\limits_{j=1}^{l}\|v_{j}\|_{X_{s,b}}
\end{eqnarray*}
outside a set of probability
$
\leq C_{1}e^{-C\lambda^{2}\|f\|_{H^{s}}^{-2}}.
$

\noindent (5). Case $N_{1}\geq80k,N_{1}\sim N_{l},N_{1}\geq 80kN_{l+1}.$
By using a  proof  similarly to Case (2) of Lemma 3.3,  we have
\begin{eqnarray*}
&&I_{6}
\leq  C\lambda^{k-l+1}\|h\|_{X_{0,b_{1}}}\prod\limits_{j=1}^{l}\|v_{j}\|_{X_{s,b}}
\end{eqnarray*}
outside a set of probability
$
\leq C_{1}e^{-C\lambda^{2}\|f\|_{H^{s}}^{-2}}.
$

\noindent (6). Case $N_{1}\geq80k,N_{1}\sim N_{l},\frac{N_{l+1}}{80k}\leq N_{1}\leq 80kN_{l+1},N_{l+1}\geq 800kN_{l+2}.$
By using a  proof  similarly to Case (2) of Lemma 3.3,  we have
\begin{eqnarray*}
&&I_{6}
\leq  C\lambda^{k-l+1}\|h\|_{X_{0,b_{1}}}\prod\limits_{j=1}^{l}\|v_{j}\|_{X_{s,b}}
\end{eqnarray*}
outside a set of probability
$
\leq C_{1}e^{-C\lambda^{2}\|f\|_{H^{s}}^{-2}}.
$

\noindent (7). Case $N_{1}\geq80k,N_{1}\sim N_{l},\frac{N_{l+1}}{80k}
\leq N_{1}\leq 80kN_{l+1},N_{l+1}\sim N_{m}\geq 80kN_{m+1}(l+2\leq m\leq k-1).$
By using a  proof  similarly to Case (2) of Lemma 3.3,  we have
\begin{eqnarray*}
&&I_{6}
\leq  C\lambda^{k-l+1}\|h\|_{X_{0,b_{1}}}\prod\limits_{j=1}^{l}\|v_{j}\|_{X_{s,b}}
\end{eqnarray*}
outside a set of probability
$
\leq C_{1}e^{-C\lambda^{2}\|f\|_{H^{s}}^{-2}}.
$

\noindent (8). Case $N_{1}\geq80k, N_{1}\sim N_{l},\frac{N_{l+1}}{80k}
\leq N_{1}\leq 80kN_{l+1},N_{l+1}\sim N_{k+1}.$
By using a  proof  similarly to Case (9) of Lemma 3.2,  we have
\begin{eqnarray*}
&&I_{6}
\leq  C\lambda^{k-l+1}\|h\|_{X_{0,b_{1}}}\prod\limits_{j=1}^{l}\|v_{j}\|_{X_{s,b}}
\end{eqnarray*}
outside a set of probability
$
\leq C_{1}e^{-C\lambda^{2}\|f\|_{H^{s}}^{-2}}.
$

\noindent (9). Case $N_{k+1}\geq80k,\frac{N_{l+1}}{80k}\geq N_{1}.$
By using the proof similar to Cases (6)-(14) of Lemma 3.3,
we have
\begin{eqnarray*}
&&I_{6}
\leq  C\lambda^{k-l+1}\|h\|_{X_{0,b_{1}}}\prod\limits_{j=1}^{l}\|v_{j}\|_{X_{s,b}}
\end{eqnarray*}
outside a set of probability
$
\leq C_{1}e^{-C\lambda^{2}\|f\|_{H^{s}}^{-2}}.
$

We have completed the proof of Lemma 3.8.

\bigskip

\noindent{\large\bf 4. Proof of Theorem 1.1}
\setcounter{equation}{0}
\setcounter{Theorem}{0}

\setcounter{Lemma}{0}

\setcounter{section}{4}

To present the proof of Theorem  1.1, as in \cite{CLS},
we use frequency truncated technique to prove Theorem 1.1.
Let $u^{N}$ be the solution to
\begin{eqnarray}
&&(u^{N})_{t}+(u^{N})_{xxx}+\frac{1}{k+1}\partial_{x}P_{\leq N}((u^{N})^{k+1})=0\label{4.01},\\
&&
u^{N}(x,0)=P_{\leq N}u(x,0)=P_{\leq N}f.\label{4.02}
\end{eqnarray}
Here,  $P_{\leq N}u(x,0)=\frac{1}{2\pi}\int_{\SR}\phi\left(\frac{\xi}{N}\right)e^{ix\xi}\mathscr{F}_{x}u(\xi,0)d\xi
=\frac{1}{2\pi}\int_{\SR}\phi\left(\frac{\xi}{N}\right)e^{ix\xi}\mathscr{F}_{x}f(\xi)d\xi.$
We define $u:=\lim\limits_{N\longrightarrow \infty}u^{N}.$

\begin{Lemma}\label{lem4.1}
Suppose $f\in H^{s}(\R)(s>\frac{1}{2}-\frac{2}{k},k\geq8).$ Then, we have
\begin{eqnarray}
\lim\limits_{N\longrightarrow \infty}
\left\|\sup\limits_{0\leq t\leq \delta}\left|u-u^{N}\right|\right\|_{L_{x}^{4}}=0\label{4.08}.
\end{eqnarray}

\end{Lemma}
\noindent {\bf Proof.}By using (\ref{2.022}), we have
\begin{eqnarray}
\left\|\sup\limits_{0\leq t\leq \delta}\left|u-u^{N}\right|\right\|_{L_{x}^{4}}\leq
\left\|\phi\left(\frac{t}{\delta}\right)\left(u-u^{N}\right)\right\|_{X_{s,b}}\label{4.09}.
\end{eqnarray}
We define $\phi\left(\frac{s}{\delta}\right)u=v, \phi\left(\frac{s}{\delta}\right)u^{N}=v^{N}.$
Obviously, for $t\in [0,\delta]$, we have
\begin{eqnarray}
&&\hspace{-1cm}\phi\left(\frac{t}{\delta}\right)\left(u-u^{N}\right)\nonumber\\&&
\hspace{-1cm}=\phi(t)U(t)P_{>N}f+\frac{1}{k+1}\phi\left(\frac{t}{\delta}\right)
\int_{0}^{t}U(t-s)\partial_{x}\left[(v^{k+1})-P_{\leq N}(v^{N})^{k+1})\right]ds\label{4.010}.
\end{eqnarray}
By using Lemma 2.1, we have
\begin{eqnarray}
 &&\left\|\phi\left(\frac{t}{\delta}\right)\left(u-u^{N}\right)\right\|_{X_{s,b}}\nonumber\\
 &&\leq C\left\|P_{>N}f\right\|_{H^{s}(\SR)}+C\delta^{\frac{\epsilon}{24}}
 \left\|\partial_{x}\left[(v^{k+1})-P_{\leq N}((v^{N})^{k+1})\right]\right\|_{X_{s,b_{1}}}\label{4.011}.
\end{eqnarray}
Obviously, we have
\begin{eqnarray}
\partial_{x}\left[(v^{k+1})-P_{\leq N}((v^{N})^{k+1})\right]=
P_{\leq N}\partial_{x}\left[(v^{k+1})-(v^{N})^{k+1})\right]+P_{>N}\partial_{x}\left[v^{k+1}\right].\label{4.012}
\end{eqnarray}
From (\ref{4.011}), (\ref{4.012}) and Lemma 3.1, we have
\begin{eqnarray}
 &&\left\|\phi\left(\frac{t}{\delta}\right)\left(u-u^{N}\right)\right\|_{X_{s,b}}
 =\left\|v-v^{N}\right\|_{X_{s,b}}\nonumber\\
 &&\leq C\left\|P_{>N}f\right\|_{H^{s}(\SR)}+C\delta^{\frac{\epsilon}{24}}
 \left[\left\|P_{\leq N}\partial_{x}\left[v^{k+1}-(v^{N})^{k+1}\right]\right\|_{X_{s,b_{1}}}
  +\left\|P_{>N}\partial_{x}\left[v^{k+1}\right]\right\|_{X_{s,b_{1}}}\right]\nonumber\\
 &&\leq C\left\|P_{>N}f\right\|_{H^{s}(\SR)}+C\delta^{\frac{\epsilon}{24}}
 \left[\left\| v-v^{N}\right\|_{X_{s,b}}^{k+1}+\left\|P_{>N}\partial_{x}\left[v^{k+1}\right]\right\|_{X_{s,b_{1}}}\right]\nonumber\\
 &&\leq C\left\|P_{>N}f\right\|_{H^{s}(\SR)}+\frac{1}{2}\left\|v-v^{N}\right\|_{X_{s,b}}^{k+1}+
 C\delta^{\frac{\epsilon}{24}}\left\|P_{>N}\partial_{x}\left[v^{k+1}\right]\right\|_{X_{s,b_{1}}}.\label{4.013}
\end{eqnarray}
From (\ref{4.013}), we have
\begin{eqnarray}
 &&\left\|\phi\left(\frac{t}{\delta}\right)\left(u-u^{N}\right)\right\|_{X_{s,b}}=\left\|v-v^{N}\right\|_{X_{s,b}}\nonumber\\&&\leq C\left\|P_{>N}f\right\|_{H^{s}(\SR)}+C\delta^{\frac{\epsilon}{24}}\left\|P_{>N}\partial_{x}\left[v^{k+1}\right]\right\|_{X_{s,b_{1}}}.\label{4.014}
\end{eqnarray}
Since $f\in H^{s}(\R),\left\|P_{>N}\partial_{x}\left[v^{k+1}\right]\right\|_{X_{s,b_{1}}}\leq C\|v\|_{X_{s,b}}^{k+1},$
we have
\begin{eqnarray}
&&\lim\limits_{N\longrightarrow \infty}\left\|\sup\limits_{0\leq t\leq \delta}\left|u-u^{N}\right|\right\|_{L_{x}^{4}}\leq \lim\limits_{N\longrightarrow \infty}\left\|\phi\left(\frac{t}{\delta}\right)\left(u-u^{N}\right)\right\|_{X_{s,b}}\nonumber\\&&\leq \lim\limits_{N\longrightarrow \infty}\left[C\left\|P_{>N}f\right\|_{H^{s}(\SR)}+C\delta^{\frac{\epsilon}{24}}\left\|P_{>N}\partial_{x}\left[v^{k+1}\right]\right\|_{X_{s,b_{1}}}\right]=0\label{4.015}.
\end{eqnarray}

Thus, we have completed the proof of Lemma 4.1.

Inspired by Proposition 3.3 of \cite{CLS}, we use Lemma 4.2 to prove Lemma 4.3.
To prove Theorem 1.1, it suffices to prove Lemma 4.2.

\begin{Lemma}\label{lem4.2}Let $k\geq8$ and $f\in H^{s}(\R),s>\frac{1}{2}-\frac{2}{k}$.
 Then, $u\longrightarrow f$ as $t\longrightarrow0$ for a.e. $x\in \R.$
\end{Lemma}
\noindent {\bf Proof.} Since $P_{\leq N}f\in H^{\infty}$, then, we have
\begin{eqnarray}
\lim\limits_{t\longrightarrow0}u^{N}=P_{\leq N}f\label{4.016}
\end{eqnarray}
for a.e. $x\in \R$. Since
\begin{eqnarray}
\left|u-f\right|\leq |u-u^{N}|+|u^{N}-P_{\leq N}f|+|P_{>N}f|,\label{4.017}
\end{eqnarray}
by using (\ref{4.016}), we have
\begin{eqnarray}
&&\lim\limits_{t\longrightarrow0}\sup\left|u-f\right|\leq \lim\limits_{t\longrightarrow0}\sup|u-u^{N}|
+\lim\limits_{t\longrightarrow0}\sup|u^{N}-P_{\leq N}f|+\lim\limits_{t\longrightarrow0}\sup|P_{>N}f|\nonumber\\
&&\leq \lim\limits_{t\longrightarrow0}\sup|u-u^{N}|+\lim\limits_{t\longrightarrow0}\sup|P_{>N}f|.\label{4.018}
\end{eqnarray}
By using the Chebyshev inequality and the Sobolev embeddings theorem $H^{\frac{1}{4}}(\R)\hookrightarrow L^{4}(\R)$
 as well as Lemma 4.1, since $f\in H^{s}(\R)$,  we have
\begin{eqnarray}
&&\left|\left\{x\in B_{1}: \lim\limits_{t\longrightarrow0}\sup\left|u-f\right|>\alpha\right\}\right|\leq
 \left|\left\{x\in B_{1}: \lim\limits_{t\longrightarrow0}\sup\left|u-u^{N}\right|>\frac{\alpha}{2}\right\}\right|
\nonumber\\&&\qquad\qquad+\left|\left\{x\in B_{1}: \lim\limits_{t\longrightarrow0}\sup\left|P_{>N}f\right|>\frac{\alpha}{2}\right\}\right|\nonumber\\
&&\leq C\alpha^{-4}\left\|\sup\limits_{0\leq t\leq \delta}\left|u-u^{N}\right|\right\|_{L_{x}^{4}}+C\alpha^{-4}\|P_{>N}f\|_{L_{x}^{4}}\nonumber\\
&&\leq C\alpha^{-4}\left\|\sup\limits_{0\leq t\leq \delta}\left|u-u^{N}\right|\right\|_{L_{x}^{4}}+C\alpha^{-4}\|P_{>N}f\|_{H_{x}^{\frac{1}{4}}}\nonumber\\
&&\leq C\alpha^{-4}\left\|\sup\limits_{0\leq t\leq \delta}\left|u-u^{N}\right|\right\|_{L_{x}^{4}}+C\alpha^{-4}\|P_{>N}f\|_{H_{x}^{s}}\label{4.019}.
\end{eqnarray}
Here $B_{1}\subset \R.$
From Lemma 4.1 and $f\in H^{s}(\R),$ we have
\begin{eqnarray}
\left|\left\{x\in B_{1}: \lim\limits_{t\longrightarrow0}\sup\left|u-f\right|>\alpha\right\}\right|\leq \left|\left\{x\in B_{1}:
\lim\limits_{t\longrightarrow0}\sup\left|u-u^{N}\right|>\frac{\alpha}{2}\right\}\right|=0.\label{4.020}
\end{eqnarray}

We have completed the proof of Lemma 4.2.

This completes the proof of Theorem  1.1.$\hfill\Box$

\bigskip
\bigskip
\noindent{\large\bf 5. Proof of Theorem 1.2}
\setcounter{equation}{0}
\setcounter{Theorem}{0}

\setcounter{Lemma}{0}

\setcounter{section}{5}

In this section, we prove Theorem 1.2.

\noindent{\bf  Proof of Theorem 1.2.}
By using (\ref{4.020}) and a proof similar to  Lemma 2.3 of \cite{D}, since $f\in H^{s}(\R)(s>\frac{1}{2}-\frac{2}{k}\geq \frac{1}{4},k\geq8)$ which yields
\begin{eqnarray}
\|U(t)f\|_{L_{x}^{4}L_{t}^{\infty}}\leq C\|f\|_{H^{\frac{1}{4}}}\leq C\|f\|_{H^{s}}\label{5.01},
\end{eqnarray}
we have
\begin{eqnarray}
\lim\limits_{t\longrightarrow0}\sup\left|u-U(t)f\right|\leq \lim\limits_{t\longrightarrow0}
\sup|u-f|+\lim\limits_{t\longrightarrow0}\sup|f-U(t)f|=0.\label{5.02}
\end{eqnarray}

\bigskip
\bigskip

\noindent{\large\bf 6. Proof of Theorem 1.3}
\setcounter{equation}{0}
\setcounter{Theorem}{0}

\setcounter{Lemma}{0}

\setcounter{section}{6}

In this section,  by using the idea of \cite{Compaan},  we will prove Theorem 1.3.

We firstly prove Lemma 6.1.

\begin{Lemma}\label{lem6.1}
Suppose $f\in H^{s}(\R)(s>\frac{1}{2}-\frac{2}{k},k\geq5).$ Then, there exists a unique solution $u\in X_{s,b}(b>\frac{1}{2})$ to (\ref{1.01})-(\ref{1.02}).
\end{Lemma}
\noindent{\bf Proof.}We define
\begin{eqnarray}
&&\Phi(u)=\phi(t)U(t)f+\frac{1}{k+1}\phi\left(\frac{t}{T}\right)\int_{0}^{t}U(t-\tau)\partial_{x}(u^{k+1})d\tau,\label{6.003}\\
&&B=\left\{u\in X_{s,b}:\|u\|_{X_{s,b}}\leq 2C\|f\|_{H^{s}}\right\}.\label{6.004}
\end{eqnarray}
By using (\ref{2.01})-(\ref{2.02}) and Lemma 3.1,  for  $T\leq  \left(\frac{1}{C^{k+1}2^{k+2}\|f\|_{H^{s}}^{k}}\right)^{\frac{1}{\epsilon}}$,     we have
\begin{eqnarray}
&&\left\|\Phi(u)\right\|_{X_{s,b}}\leq \left\|\phi(t)U(t)f\right\|_{X_{s,b}}+
\left\|\frac{1}{k+1}\phi\left(\frac{t}{T}\right)\int_{0}^{t}U(t-\tau)\partial_{x}(u^{k+1})d\tau\right\|_{X_{s,b}}\nonumber\\
&& \leq C\|f\|_{H^{s}}+CT^{\epsilon}\left\|\partial_{x}(u^{k+1})\right\|_{X_{s,b_{1}}}\nonumber\\
&&\leq C\|f\|_{H^{s}}+CT^{\epsilon}\|u\|_{X_{s,b}}^{k+1}\nonumber\\
&&\leq C\|f\|_{H^{s}}+CT^{\epsilon}(2C\|f\|_{H^{s}})^{k+1}\leq 2C\|f\|_{H^{s}}\label{6.005}
\end{eqnarray}
and
\begin{eqnarray}
&&\left\|\Phi(u)-\Phi(v)\right\|_{X_{s,b}}\leq
\left\|\frac{1}{k+1}\phi\left(\frac{t}{T}\right)\int_{0}^{t}U(t-\tau)\partial_{x}(u^{k+1}-v^{k+1})d\tau\right\|_{X_{s,b}}\nonumber\\
&&\leq CT^{\epsilon}\|u-v\|_{X_{s,b}}\left[\|u\|_{X_{s,b}}^{k}+\|v\|_{X_{s,b}}^{k}\right]\leq  CT^{\epsilon}(2C\|f\|_{H^{s}})^{k}\|u-v\|_{X_{s,b}}\nonumber\\
&&\leq \frac{1}{2}\|u-v\|_{X_{s,b}}\label{6.006}.
\end{eqnarray}
Thus, $\phi$ is a contraction mapping from $B$ to $B$. Consequently, $\Phi$ has a fixed point. That is $\Phi(u)=u.$ From (\ref{6.006}), we have
\begin{eqnarray}
\left\|u-v\right\|_{X_{s,b}}\leq
\frac{1}{2}\|u-v\|_{X_{s,b}}\label{6.007}.
\end{eqnarray}

We have completed the proof of Lemma 6.1.

\noindent{\bf Proof of Theorem 1.3.}We define
$b:=\frac{1}{2}+\frac{\epsilon}{24},b_{1}=\frac{1}{2}-\frac{\epsilon}{12}.$
By using Duhamel's formula, we can rewrite  the local solution to (\ref{4.01})-(\ref{4.02})  as follows.
\begin{eqnarray}
u=\phi(t)U(t)f+\frac{1}{k+1}\phi\left(\frac{t}{T}\right)\int_{0}^{t}U(t-\tau)\partial_{x}(u^{k+1})d\tau,\label{6.01}
\end{eqnarray}
For $s\geq \frac{1}{2}-\frac{2}{k+1}+88\epsilon,$ for sufficiently small $t>0,$  by using Lemmas 3.2, 6.1 and
 $X_{\frac{1}{2}+\epsilon,b}\hookrightarrow C([0,T];H^{\frac{1}{2}+\epsilon})$ as well as (\ref{6.004}),  we have
\begin{eqnarray}
&&\left\|u-\phi(t)U(t)f\right\|_{L_{x}^{\infty}}=\left\|u-U(t)f\right\|_{L_{x}^{\infty}}
\leq C\left\|u-U(t)f\right\|_{C([0,T];H_{x}^{\frac{1}{2}+\epsilon})}\nonumber\\&&= C\left\|u-U(t)f\right\|_{X_{\frac{1}{2}+\epsilon,b}}\nonumber\\&&=
C\left\|\phi\left(\frac{t}{T}\right)\int_{\SR}U(t-\tau)\partial_{x}\left[(u)^{k+1}\right]\right\|_{X_{\frac{1}{2}+\epsilon,b}}
\nonumber\\&&\leq CT^{\epsilon}\|\partial_{x}\left[u^{k+1}\right]\|_{X_{s,b_{1}}}\nonumber\\&&\leq  CT^{\epsilon}\|u\|_{X_{s,b}}^{k+1}
\leq CT^{\epsilon}\|f\|_{H^{s}}^{k+1}<\infty.\label{6.02}
\end{eqnarray}
From (\ref{6.02}), we have
\begin{eqnarray}
\lim\limits_{t\longrightarrow 0}\left\|u-U(t)f\right\|_{L_{x}^{\infty}}=0.\label{6.03}
\end{eqnarray}

This completes the proof of Theorem  1.3.$\hfill\Box$

\bigskip
\bigskip

\noindent {\large\bf 7. Proof of Theorem  1.4}

\setcounter{equation}{0}

 \setcounter{Theorem}{0}

\setcounter{Lemma}{0}

\setcounter{section}{7}
\noindent {\bf Proof.} For $1\leq l\leq k$, by using Lemmas 3.1-3.7 and the Young inequality, we have
\begin{eqnarray}
&&\left\|\partial_{x}\left[(v+z)^{k+1}\right]\right\|_{X_{s,-\frac{1}{2}+\frac{\epsilon}{12}}}\leq C
\left[\|v\|_{X_{s,b}}^{k+1}+\lambda^{k-l+1}\|v\|_{X_{s,b}}^{l}+\lambda^{k+1}\right]\nonumber\\&&
\leq C\left[\|v\|_{X_{s,b}}^{k+1}+\lambda^{k+1}\right]\label{7.01}
\end{eqnarray}
outside a set of probability
$
\leq C_{1}e^{-C\lambda^{2}\|f\|_{H^{s}}^{-2}}.
$

This completes the proof of Theorem 1.4.

\bigskip
\bigskip

\noindent {\large\bf 8. Proof of Theorem  1.5}

\setcounter{equation}{0}

 \setcounter{Theorem}{0}

\setcounter{Lemma}{0}

\setcounter{section}{8}

\noindent{\bf Proof of Theorem 1.5.} We define
$b:=\frac{1}{2}+\frac{\epsilon}{24},b_{1}=\frac{1}{2}-\frac{\epsilon}{12}.$
Assume that $f^{\omega}$ is defined as in (\ref{1.05}), which
 belongs to $H^{s_{1}}(\R)$ almost surely.
Now we consider the Cauchy problem for (\ref{1.01}) with $u(x,0)=f^{\omega}$.
Let $z(t):=z^{\omega}(t)=U(t)f^{\omega}$ and $v(t)=u(t)-z(t)$, (\ref{1.01})
 can be rewritten as follows:
\begin{eqnarray}
&&v_{t}+v_{xxx}+\frac{1}{k+1}\partial_{x}\left[(v+z)^{k+1}\right]=0,\label{8.01},\\
&&v(x,0)=0\label{8.02}.
\end{eqnarray}
(\ref{8.01})-(\ref{8.02}) are equivalent to the following integral equation
\begin{eqnarray}
v(t)=\frac{1}{k+1}\int_{0}^{t}U(t-\tau)\partial_{x}\left[(v+z)^{k+1}\right]d\tau.\label{8.03}
\end{eqnarray}
Obviously, $v$ satisfies (\ref{8.03}) on $[-T,T]$ if $v$ satisfies
\begin{eqnarray}
v(t)=\frac{1}{k+1}\phi\left(\frac{t}{T}\right)\int_{0}^{t}U(t-\tau)
\partial_{x}\left[(v+\phi(t)z)^{k+1}\right]d\tau.\label{8.04}
\end{eqnarray}
for some $T\ll1$. For $t\in[0,T],$
we define
\begin{eqnarray}
\Gamma v=\frac{1}{k+1}\phi\left(\frac{t}{T}\right)\int_{0}^{t}U(t-\tau)
\partial_{x}\left[(v+\phi(t)z)^{k+1}\right]d\tau.\label{8.05}
\end{eqnarray}
By using Lemma 2.1, Theorem 1.4, we have
\begin{eqnarray}
&&\|\Gamma v\|_{X_{s,b}}\leq CT^{\frac{\epsilon}{24}}\left\|
\partial_{x}\left[(v+\phi(t)z)^{k+1}\right]\right\|_{X_{s,-b_{1}}}\nonumber\\
&&\leq CT^{\frac{\epsilon}{24}}
\left(\|v\|_{X_{s,b}}^{k+1}+\lambda^{k+1}\right)\label{8.06}
\end{eqnarray}
and
\begin{eqnarray}
&&\|\Gamma v_{1}-\Gamma v_{2}\|_{X_{s,b}}\leq C
T^{\frac{\epsilon}{24}}\|v_{1}-v_{2}\|_{X_{s,b}}
\left(\|v_{1}\|_{X_{s,b}}^{k}+\|v_{2}\|_{X_{s,b}}^{k}+\lambda^{k}\right)\label{8.07}
\end{eqnarray}
outside a set of probability  $
\leq C_{1}e^{-C\lambda^{2}\|f\|_{H^{s}}^{-2}}.
$
Let $B_{2}=\left\{v:\|u\|_{X_{s,b}}\leq1\right\}$, for
$
T\leq1$,
we choose
 $\lambda=\lambda(T)\sim T^{-\frac{\epsilon}{48k}}$
such that
\begin{eqnarray}
&&CT^{\frac{\epsilon}{24}}(1+\lambda^{k+1})\leq 1,CT^{\frac{\epsilon}{24}}(2+\lambda^{k})
\leq \frac{1}{2}\label{8.08}.
\end{eqnarray}
Thus, for $v,v_{1},v_{2}\in B_{1}$, we have
\begin{eqnarray}
&&\|\Gamma v\|_{X_{s,b}}\leq1,\|\Gamma v_{1}-\Gamma v_{2}\|_{X_{s,b}}
\leq\frac{1}{2}\|v_{1}-v_{2}\|_{X_{s,b}}\label{8.09}
\end{eqnarray}
outside an  set of probability at most
\begin{eqnarray*}
C_{1}{\rm exp}\left(-C\frac{\lambda^{2}}{\|f\|_{H^{s_{1}}}^{2}}\right)
=C_{1}{\rm exp}\left(-\frac{C}{T^{\frac{\epsilon}{24k}}\|f\|_{H^{s_{1}}}^{2}}\right).
\end{eqnarray*}
Consequently, let $\Omega_{T}$  be the complement of the exceptional set,
for $\omega \in \Omega_{T}$, there exists a unique $v^{\omega}\in B_{2}$ such
that $\Gamma v^{\omega}=v^{\omega}$.

We have completed the proof of Theorem 1.5.

\bigskip
\bigskip

\noindent{\large\bf 9. Proof of Theorem 1.6}
\setcounter{equation}{0}
\setcounter{Theorem}{0}

\setcounter{Lemma}{0}

\setcounter{section}{9}

In this section,  we will prove Theorem 1.6.

\noindent{\bf Proof of Theorem 1.6.}
Let $v=u-U(t)f^{\omega}.$ Then,
 $\forall \omega \in \Omega_{T}$,  by using  $X_{\frac{1}{2}+\epsilon,b}\hookrightarrow C([0,T];H^{\frac{1}{2}+\epsilon})$ and Theorem 1.4,  we have
\begin{eqnarray}
&&\left\|u-U(t)f^{\omega}\right\|_{L_{x}^{\infty}}=\left\|v\right\|_{L_{x}^{\infty}}
\leq C\left\|v\right\|_{C([0,T];H_{x}^{\frac{1}{2}+\epsilon})}\leq
C\left\|v\right\|_{X_{\frac{1}{2}+\epsilon,b}}\nonumber\\&&\leq C
\left\|\int_{\SR}U(t-\tau)\partial_{x}\left[(v+U(t)f^{\omega})^{k+1}\right]\right\|_{X_{\frac{1}{2}+\epsilon,b}}
\leq C\left[\|v\|_{X_{s,b}}^{k+1}+\lambda^{k+1}\right]\nonumber\\
&&\leq C\left[1+\lambda^{k+1}\right].\label{9.01}
\end{eqnarray}
From (\ref{9.01}) and $v(0)=0$, we have
\begin{eqnarray}
\lim\limits_{t\longrightarrow 0}\|v(t)\|_{H^{s}}=\lim\limits_{t\longrightarrow 0}\left|\|v(t)\|_{H^{s}}-\|v(0)\|_{H^{s}}\right|=0.\label{9.02}
\end{eqnarray}
From (\ref{9.02}) and $H^{s}(\R)\hookrightarrow L^{\infty}(\R)(s>\frac{1}{2})$, we have
\begin{eqnarray}
\lim\limits_{t\longrightarrow 0}\left\|u-U(t)f^{\omega}\right\|_{L_{x}^{\infty}}
=\lim\limits_{t\longrightarrow 0}\left\|v\right\|_{L_{x}^{\infty}}
\leq C\lim\limits_{t\longrightarrow 0}\|v(t)\|_{H^{s}}=0.\label{9.03}
\end{eqnarray}
From (\ref{9.03}), we have
\begin{eqnarray}
\lim\limits_{t\longrightarrow 0}\left\|u-U(t)f^{\omega}\right\|_{L_{x}^{\infty}}=0.\label{9.04}
\end{eqnarray}

This completes the proof of Theorem  1.6.$\hfill\Box$

\bigskip
\bigskip
\noindent{\large\bf 10. Proof of Theorem 1.7}
\setcounter{equation}{0}
\setcounter{Theorem}{0}

\setcounter{Lemma}{0}

\setcounter{section}{10}
In this section, we prove Theorem 1.7.

\noindent{\bf Proof of Theorem 1.7.}

Combining Theorem 1.5 with  Theorem 1.1 of \cite{YDLY},  $\forall \omega \in \Omega_{T},$ we have
\begin{eqnarray}
\lim\limits_{t\longrightarrow0}\sup\left|u-f^{\omega}\right|\leq \lim\limits_{t\longrightarrow0}
\sup|u-U(t)f^{\omega}|+\lim\limits_{t\longrightarrow0}\sup|f^{\omega}-U(t)f^{\omega}|=0.\label{10.01}
\end{eqnarray}

We have completed the proof of Theorem  1.7.$\hfill\Box$

\bigskip

\leftline{\large \bf Acknowledgments}

\noindent Wei Yan was supported by NSFC grants (No. 11771127).


\bigskip

  \leftline{\large\bf  References}
\begin{thebibliography}{99}
\bibitem{BOP} A. B$\acute{e}$nyi, T. Oh, O. Pocovnicu, Wiener randomization
on unbounded domains
and an application to almost  sure
well-posedness of NLS, Excursions in Harmonic Analysis.

\bibitem{BOP2015} A. B$\acute{e}$nyi, T. Oh, O. Pocovnicu, On the probabilistic
  Cauchy theory of
the cubic nonlinear Schr\"odinger
equation on $\R^{d},$ $d\geq3$, {\it Tran. Amer. Math. Soc.} 2(2015), 1-50.















\bibitem{Bourgain1992} J. Bourgain,  A remark on Schr\"odinger operators,
{\it Isreal J. Math.}  77(1992), 1-16.

\bibitem{Bourgain-S}J.  Bourgain,  Fourier transform restriction phenomena for certain lattice subsets and applications to nonlinear evolution equations,  I. Schr\"odinger equations,
{\it  Geom. Funct. Anal.} 3(1993),  107-156.

\bibitem{Bourgain-GAFA93}
J. Bourgain,
Fourier transform restriction phenomena for certain
lattice subsets and applications to nonlinear evolution equations,
part II: The KdV equation, {\it Geom.   Funct. Anal.} 3(1993), 209-262.



 \bibitem{B1994CMP} J. Bourgain,    Periodic  Schr\"odinger  equation
 and invariant measures,
  {\it Comm. Math. Phys.} 166(1994),  1-26.








\bibitem{B1994Duke} J. Bourgain, On the Cauchy and invariant measure problem for the
 periodic Zakharov system, {\it Duke. Math. J.} 76(1994), 175-202.



 \bibitem{Bourgain1995}J. Bourgain,  Some new estimates on osillatory integrals,
  In: Essays on Fourier Analysis in Honor of
Elias M. Stein, Princeton, NJ 1991. Princeton Mathematical Series, vol. 42,
 pp. 83.112. Princeton
University Press, New Jersey (1995).









\bibitem{B1996CMP} J. Bourgain, Invariant measures for the 2 D-defocusing
 nonlinear Schr\"odinger equation,  {\it Comm. Math. Phys.} 176(1996), 421-445



\bibitem{Bourgain2013}J. Bourgain,  On the Schr\"odinger maximal
 function in higher dimensions,  {\it Proc. Steklov Inst. Math.}
280(2013), 46-60.



\bibitem{B-2014}
J. Bourgain, A. Bulut,  Almost sure global well-posedness for
the radial nonlinear Schr\"odinger equation on the unit ball II: the 3d case.
{\it J. Eur. Math. Soc. }
16(2014), 1289-1325


\bibitem{BB2014}
J.  Bourgain, A. Bulut,  Invariant Gibbs measure evolution for the
radial nonlinear wave equation on the 3d ball.
{\it J. Funct. Anal.}  266(2014), 2319-2340.




\bibitem{Bourgain2016} J. Bourgain,  A note on the Schr\"odinger
maximal function,  {\it J. Anal. Math.} 130(2016),  393-396.






\bibitem{BT2007}
N. Burq, N. Tzvetkov,  Invariant measure for a three dimensional
nonlinear wave equation, {\it Int. Math. Res. Not.}  2007, no. 22, Art. ID rnm108, 26 pp

\bibitem{BT2008-L}
N. Burq,  N. Tzvetkov,  Random data Cauchy theory for supercritical wave equations,
 I. Local theory,  {\it Invent. Math.}  173(2008),449-475.

\bibitem{BT2008-G}
N. Burq,  N. Tzvetkov, Random data Cauchy theory for supercritical wave equations,
 II. A global existence result,  {\it Invent. Math.} 173(2008), 477-496.

\bibitem{BT2014}
N. Burq, N. Tzvetkov,  Probabilistic well-posedness for the cubic
wave equation, {\it  J. Eur. Math. Soc.}  16(2014),  1-30.


\bibitem{Carleson} L. Carleson,  Some analytical problems related to statistical
 mechanics. Euclidean Harmonic Analysisi.
Lecture Notes in Mathematics, vol. 779, pp. 5.45, Springer, Berlin, (1979).


\bibitem{CG}
Y. Chen, H. J. Gao, The  Cauchy problem for the Hartree equations
  under random influences,
{\it J. Diff. Eqns.}  259(2015), 5192-5219.



\bibitem{CLV}C. Cho, S. Lee and  A. Vargas, Problems on pointwise convergence of solutions
to the Schr\"odinger
equation, {\it J. Fourier Anal. Appl.} 18(2012), 972-994.

\bibitem{CO}
J.  Colliander, T. Oh,  Almost sure well-posedness of
the cubic nonlinear
 Schr\"odinger equation below $L^{2}(\mathbf{T})$,
 {\it Duke Math. J.} 161(2012), 367-414.

\bibitem{Compaan} E. Compaan, A smoothing estimate for the nonlinear Schr\"odinger equation,
{\it UIUC Research Experience for Graduate  Students report, 2013}.

\bibitem{CLS} E. Compaan, R. Luc$\grave{a}$, G. Staffilani,
  Pointwise convergence of the Schr$\ddot{o}$dinger
flow, {\it  Int. Math. Res. Not.} 1(2021),  599-650.


\bibitem{Cowling}M.  Cowling,  Pointwise behavior of solutions to
 Schr\"odinger equations. In: Harmonic Analysis (Cortona,
1982). Lecture Notes in Mathematics, vol. 992, pp. 83.90. Springer, Berlin, (1983).




\bibitem{DJEMS} Y. Deng,  Invariance of the Gibbs measure for the
 Benjamin-Ono equation,
{\it  J. Eur. Math. Soc. } 17(2015),  1107-1198.

\bibitem{DC}C. Deng, S. B. Cui, Random-data Cauchy problem for
 the Navier-Stokes equations
on $\mathbf{T}^3,$ {\it  J. Diff. Eqns.} 251(2011),  902-917.

 \bibitem{DK} B. E. Dahlberg and  C. E. Kenig, A note on the almost everywhere behavior
  of solutions to the Schr\"odinger
equation. In: Proceedings of Italo-American Symposium in Harmonic Analysis,
University of Minnesota.
Lecture Notes in Mathematics, vol. 908, pp. 205-209. Springer, Berlin, (1981).





\bibitem {DG} C. Demeter and   S. Guo,  Schr\"odinger maximal function estimates via
 the pseudoconformal transformation,
 arXiv: 1608.07640.

 \bibitem{ST}
A. de Suzzoni,   N. Tzvetkov,  On the propagation of
 weakly nonlinear random dispersive waves,
  {\it  Arch. Ration. Mech. Anal.}
   212(2014),  849-874.


\bibitem{de2013}A. de Suzzoni, Large data low regularity
  scattering results for
 the wave equaion on the Euclidean space,
{\it Comm. Partial Diff. Eqns.} 38(2012), 1-49.

\bibitem{de2014}
A. de Suzzoni, Consequences  of the Choice of a  particular basis
of $L^{2}(S^{3})$
for the cubic wave equation on the sphere and the Euclidean space,
{\it Comm. Pure Appl. Anal.}  13(2012), 991-1015.


\bibitem{DGL} X. Du, L. Guth and   X. Li,  A sharp Schr\"odinger maximal estimate in $\R^2,$
  {\it Ann. Math.} 188(2017), 607-640.

 \bibitem{D} X. Du,   A sharp Schr\"odinger maximal estimate in $\R^2,$
  {\it Dissertation} 2017.


\bibitem{DZ}X.  Du and  R.  Zhang, Sharp $L^2$ estimates of the Schr\"odinger maximal function
 in higher dimensions, {\it Ann. Math.} 189(2019),
837-861.

\bibitem{DGZ}X.  Du, L. Guth,  X. Li and   R. Zhang, Pointwise convergence of Schr\"{o}dinger solutions
 and multilinear refined Strichartz estimates, {\it Forum Math.Sigma}, 6(2018).






\bibitem{EL} D. Eceizabarrena and R. Luc$\grave{a}$, Convergence over fractals
 for the periodic Schr\"odinger equation,
arXiv:2005.07581.








\bibitem{FS} C. Federico,  A. S. de Suzzoni,  Invariant measure for
 the Schr\"odinger equation
on the real line, {\it J. Funct. Anal.} 269(2015), 271-324.


\bibitem{GS}G. Gigante and  F. Soria,  On the the boundedness in $H^{1/4}$ of the maximal square
 function associated with the Schr\"odinger equation,  {\it J. Lond. Math. Soc.} 77(2008), 51-68.



\bibitem{G}A.  Gr\"unrock,   New applications of the Fourier restriction
norm method to wellposedness problems for nonlinear Evolution Equations,
Ph.D. Universit$\ddot{a}$t Wuppertal,  2002, Germany, Dissertation.





\bibitem{HO} H. Hirayama, M. Okamoto, Random data Cauchy
 theory for the
fourth order nonlinear Schr\"odinger equation with cubic
 nonlinearity, arXiv:1505.06497 .



\bibitem{HO15}   H. Hirayama, M. Okamoto,
Random data Cauchy problem for the
 nonlinear Schr\"odinger
equation with derivative nonlinearity,
{\it Dis. Contin. Dyn. Syst.} 36(2016),
 6943-6974.

\bibitem{HK}
G. Hwang, C. Kwak, Probabilistic well-posedness of generalized KdV,
{\it Proc. Amer. Math. Soc.}  146(2018), 267-280.








\bibitem{LR}  F. Linares, Jo$\tilde{a}$o P. G. Ramos,  Maximal function estimates
 and local well-posedness
 for the generalized Zakharov-Kuznetsov equation,
{\it  SIAM J. Math. Anal.} 53(2021),  914-936.


\bibitem{LR2019}R. Luc$\grave{a}$, K.  M. Rogers,  A note on pointwise convergence
 for the Schr\"odinger equation,
{\it Math. Proc. Cambridge Philos. Soc.} 166(2019),  209-218.



\bibitem{KPV1991}  C. E. Kenig,  G. Ponce and  L. Vega, Oscillatory integrals and
 regularity of dispersive equations.
{\it  Indiana Univ.  Math. J.} 40(1991),  33-69.



\bibitem{KPV1996} C.     Kenig, G. Ponce, L. Vega, A bilinear estimate with
 applications to
   the KdV equation, {\it  J. Amer. Math. Soc.}  9(1996),  573-603.


\bibitem{KPV-Duke} C. E. Kenig, G. Ponce, L. Vega,
The Cauchy problem for the Korteweg-de Vries equation in Sobolev spaces of
 negative indices,
 {\it Duke Math.} J. 71(1993),  1-21.


\bibitem{KPVCPAM}  C. E. Kenig, G. Ponce, L. Vega,  Well-posedness
 and scattering results for the generalized Korteweg-de Vries equation
  via the contraction  principle,  {\it  Comm. Pure Appl. Math.}
  46(1993), 527-620.





\bibitem{KMV}R.  Killip, J.  Murphy, M. Visan,
 Invariance of white noise for KdV on the line,  {\it Invent. Math.}
   222(2020),   203-282.





\bibitem{LRS} J. Lebowitz, H. Rose, E. Speer,
Statistical mechanics of the nonlinear Schr\"odinger equation,
{\it J. Stat. Phys.} 50(1988), 657-687.




\bibitem{LMCPDE}J. L\"uhrmann and  D. Mendelson,  Random data Cauchy
 theory for nonlinear wave equations of power-type on
$\R^3$, {\it Comm. Partial Diff. Eqns.}  39(2014),  2262-2283.






\bibitem{Lee}S.  Lee,  On pointwise convergence of the solutions
to Schr\"odinger equation
 in $\R^2.$ {\it Int. Math. Res. Not.}  2006, 32597.

\bibitem{LR2015}R. Luc$\grave{a}$ and  K. M. Rogers,  An improved neccessary
 condition for Schr\"odinger
 maximal estimate,  arXiv:
1506.05325.

\bibitem{LR2017}R. Luc$\grave{a}$ and K.  M. Rogers,  Coherence on fractals
 versus pointwise convergence
 for the Schr\"odinger equation,
{\it Commun. Math. Phys.} 351(2017), 341-359.




\bibitem{MYZ2015} C. Miao, J. Yang and  J. Zheng, An improved maximal inequality for 2D
 fractional order Schr\"odinger
operators,  Stud. Math. 230(2015), 121-165.

\bibitem{MZZ2015}C.  Miao, J. Zhang and  J. Zheng, Maximal estimates for Schr\"odinger
 equation with inverse-square potential,
{\it Pac. J. Math.} 273(2015), 1-19,


\bibitem{MVV}A.  Moyua, A. Vargas and  L. Vega,  Schr\"odinger maximal function and
restriction properties of the Fourier
transform,  {\it IMRN},  1996(1996), 793-815.

\bibitem{MV} A. Moyua and L. Vega,
Bounds for the maximal function associated to periodic solutions of one-dimensional
 dispersive equations,
{\it Bull.  Lon. Math. Soc.}  40(2008),  117-128.

\bibitem{NPS} A. R. Nahmod, N. Pavlovi$\acute{c}$, G. Staffilani, Almost sure existence of
global weak solutions for supercritical Navier-Stokes equations,  {\it SIAM J. Math. Anal.}
 45(2013),  3431-3452.


\bibitem{NO}A. Nahmod, T.  Oh,  L. Rey-Bellet, G. Staffilani, Invariant
  weighted Wiener measures and almost sure global well-posedness for
   the periodic derivative NLS, {\it  J. Eur. Math. Soc.} 14(2012),  1275-1330.

   \bibitem{Oh} T.  Oh, Invariance of the white noise for KdV,
 {\it  Comm. Math. Phys.}  292(2009),  217-236.




\bibitem{OR}T. Oh, G. Richards, L. Thomann,  On invariant Gibbs measures for the
generalized KdV equations,  {\it Dyn. Partial Diff. Eqns.} 13(2016),  133-153.



 \bibitem{OhCMP} T. Oh,  Invariance of the white noise for KdV,
 {\it  Comm. Math. Phys.} 292(2009),  217-236.



\bibitem{OSIAM} T. Oh,  Invariance of the Gibbs measure for
 the Schr\"odinger-Benjamin-Ono system,
{\it  SIAM J. Math. Anal.} 41(2009),  2207-2225.



\bibitem{OhDIE}   T. Oh, Invariant Gibbs measures and a.s. global well posedness
 for coupled KdV systems,
{\it Diff. Int. Eqns.} 22(2009), 637-668.



\bibitem{OQ2013} T. Oh, J. Quastel, On invariant Gibbs measures conditioned on mass and momentum,
{\it  J. Math. Soc. Japan} 65(2013),  13-35.


\bibitem{OP} T. Oh, O. Pocovnicu, Probabilistic global well-posedness of
 the energy-critical defocusing
quintic nonlinear wave equation on $\R^{3}$,
{\it J.  Math.  Pures  Appl.} 105(2016), 342-366.





\bibitem{ORT} T. Oh, G. Richards, L. Thomann,  On invariant Gibbs measures
 for the generalized KdV  equations,
{\it Dyn.  Partial Diff. Eqns.}
13(2016), 133-153.



   \bibitem{PZ} R. Paley, A. Zygmund, On some series of functions (1), (2), (3),
{\it Proc. Camb. Philos. Soc.} 26(1930), 337-357, 458-474; 28(1932), 190-205.


\bibitem{P} O. Pocovnicu,  Almost sure global well-posedness for the
 energy-critical defocusing nonlinear wave equation on $\R^d,$  $d=4$ and $5$,
  {\it J. Eur. Math. Soc.}   19(2017),  2521-2575.


\bibitem{Poiret12}A. Poiret, Solutions globales pour des $\acute{e}$quation de
Schr\"odinger sur-critiques en toutes dimensions [Global solutions for
supercritical Schr\"odinger equations in all dimensions]. Preprint, arXiv:1207.3519.


\bibitem{PRT}A. Poiret, D. Robert, L. Thomann, Probabilistic global
 well-posedness for the supercritical nonlinear harmonic oscillator,
 {\it  Anal.   Partial  Diff.  Eqns.} 7(2014), 997-1026.






 \bibitem{Porn}D. Pornnopparath,
Small data well-posedness for derivative nonlinear Schr\"odinger equations,
{\it J. Diff. Eqns.} 265(2018),   3792-3840.














   \bibitem{R} G. Richards, Invariance of the gibbs measure for the periodic
 quartic gKdV,   {\it Annales de l'Institut Henri Poincare (C) Nonlinear Analysis,}
33(2016), 699-766.


\bibitem{RVV} K. M.  Rogers, A. Vargas and  L. Vega,  Pointwise convergence of solutions
 to the nonelliptic Schr\"odinger
equation,  {\it Indiana Univ. Math. J.} 55(2006), 1893-1906.

 \bibitem {RV} K. M.  Rogers and  P. Villarroya,  Sharp estimates for maximal operators
  associated to the wave equation,  Ark.
Mat. 46(2008), 143-151.



\bibitem{Stein}
E.M. Stein, Singular integrals and differentiability properties of functions,
{\it Princeton University Press, New Jersey, 1970}.

\bibitem{Stein-new} E.M. Stein,  Harmonic analysis : real-variable methods,
 orthogonality,
and oscillatory integrals,
{\it Princeton University Press, New Jersey,	1993}.










\bibitem{S} P. Sj\"olin,  Regularity of solutions to the Schr\"odinger equation,
{\it  Duke Math. J.}  55(1987), 699-715.

\bibitem{Shao} S. Shao, On localization of the Schr\"odinger maximal operator,
 arXiv: 1006.2787v1.




\bibitem{Strichartz}   R.   Strichartz,
 Restrictions of  Fourier transforms to  surfaces and  decay  of  solutions
of wave equation,  {\it Duke  Math.  J.}  44(1977), 705-714.

\bibitem {Tao} T. Tao,  A sharp bilinear restriction estimate for parabloids,
 Geom. Funct. Anal. 13(2003), 1359-1384.

\bibitem{TV}T. Tao and  A. Vargas, A bilinear approach to cone multipliers,
II. Appl. {\it Geom. Funct. Anal.} 10(2003), 216-258.














\bibitem{TT} L.  Thomann, N.  Tzvetkov, Gibbs measure for the periodic derivative
nonlinear Schr\"odinger equation, {\it Nonlinearity} 23(2010), 2771-2791.


\bibitem{T2009} L. Thomann, Random data  Cauchy problem for supercritical
 Schr\"odinger equations,
{\it Ann. Inst. H. Poincar$\acute{e}$ Anal. Nonlin$\acute{e}$aire,} 26(2009), 2385-2402.


 \bibitem{T}  N. Tzvetkov,  Invariant measures for the nonlinear Schr\"odinger equation
 on the disc, {\it Dyn. Partial Diff. Equ.} 3(2006),  111-160.



 \bibitem{TV2013}
 N. Tzvetkov, N. Visciglia,
 Gaussian measures associated to the  higher order conservation laws of
 the Benjamin-Ono equation, {\it  Ann. Sci. $\acute{E}$c. Norm. Sup$\acute{e}$r.}
   46(2013), 249-299.



\bibitem{TV2014}  N. Tzvetkov, N. Visciglia,  Invariant measures and long-time behavior for
the Benjamin-Ono equation, {\it Int. Math. Res. Not.}  17(2014),   4679-4714.



\bibitem{Vega}L. Vega,  Schr\"odinger equations: pointwise convergence to
the initial data,  Proc. Am. Math. Soc.
102(1988), 874-878.



\bibitem{WYY} J. F. Wang, W. Yan,  X. Q. Yan,  Probabilistic pointwise convergence problem
 of Schr\"odinger equations on manifolds,
 {\it  Proc. Amer. Math. Soc.} 149(2021), 3367-3378.

\bibitem{WZ}X. Wang and  C. J. Zhang,  Pointwise convergence of solutions to the
 Schr\"odinger equation on manifolds,
{\it  Canad. J. Math. }  71(2019),  983-995.








\bibitem{YDLY}   Wei Yan, Jinqiao Duan, Yongsheng Li, Meihua Yang,
Probabilistic pointwise convergence problem of some dispersive equations, {\it  arXiv:2004.01553}




\bibitem{ZF2011} T. Zhang, D. Y. Fang,  Random data Cauchy theory for the incompressible
three dimensional Navier-Stokes equations,
 {\it Proc. Amer. Math. Soc.}  139(2011), 2827-2837.


\bibitem{ZF2012} T. Zhang, D. Y. Fang, Random data Cauchy theory for the generalized
incompressible Navier-Stokes equations, {\it J. Math. Fluid Mech.} 14(2012),  311-324.
\end{thebibliography}
\end{document}